\documentclass[reqno, 11pt]{amsart}
\usepackage{amssymb,amsmath,amsfonts}
\usepackage[margin=1in]{geometry}
\usepackage{dsfont}
\usepackage{color}
\usepackage{graphicx,graphics}
\usepackage{latexsym}
\usepackage[dvips]{epsfig}
\usepackage{mathtools,mathrsfs}
\usepackage{stmaryrd,cite}
\usepackage{cases}
\usepackage{enumerate}
\usepackage[usenames,dvipsnames,svgnames]{xcolor}
\usepackage{todonotes}
\usepackage{marginnote}
\usepackage{listings}
\usepackage{cleveref}
\allowdisplaybreaks
\definecolor{refkey}{rgb}{1,0,0.5}
\definecolor{labelkey}{rgb}{0,0.4,1}
\numberwithin{equation}{section}
\newtheorem{thm}{Theorem}[section]

\newtheorem{lem}[thm]{Lemma}
\newtheorem{prop}[thm]{Proposition}

\newcommand{\ea}{\epsilon}

\newcommand{\al}{\alpha}

\newcommand{\da}{\delta}

\newcommand{\na}{\nabla}

\newcommand{\sa}{\sigma}

\newcommand{\ga}{\gamma}
\newcommand{\ba}{\beta}

\newcommand{\za}{\zeta}

\newcommand{\oa}{\omega}
\newcommand{\Oa}{\Omega}

\newcommand{\iy}{\infty}
\newcommand{\pl}{\partial}

\newcommand{\lt}{\left}
\newcommand{\rt}{\right}
\newcommand{\be}{\begin{equation}}
\newcommand{\ee}{\end{equation}}
\newcommand{\bee}{\begin{equation*}}
\newcommand{\eee}{\end{equation*}}

\newcommand{\ef}{\eqref}

\newcommand{\les}{\lesssim}

\begin{document}
\title[] {Global Solution to the Vacuum Free Boundary Problem with Physical Singularity  of Compressible Euler Equations with Damping and Gravity}
\maketitle
\begin{center}
\author{Huihui Zeng}\\
{\small Department of Mathematical Sciences,
Tsinghua University,
Beijing, 100084, China
}\\
{\small E-mail address: hhzeng@mail.tsinghua.edu.cn}
\end{center}

\noindent {\small\textbf{MSC}: 35Q35, 76N10

\begin{abstract} The global existence of smooth solutions to the  vacuum free boundary problem with physical singularity of compressible Euler equations with damping and gravity is proved  in space dimensions $n=1, 2, 3$, for the initial data being small perturbations of the stationary solution. Moreover, the exponential decay of the velocity is obtained for $n=1, 2,  3$. The exponentially fast convergence of the density and vacuum boundary to those of the stationary solution is shown for $n=1$, and it is proved for $n=2, 3$ that they stay close to those of the stationary solution if they do so initially.
The proof is based on the weighted estimates of both hyperbolic and parabolic types with weights capturing the singular behavior of higher-order normal derivatives near vacuum states, exploring the balance between the physical singularity  which pushes the vacuum boundary outwards and the effect of gravity which pulls it inwards,  and the dissipation of the frictional damping.
The results obtained in this paper are the first ones on the global existence of solutions to the vacuum free boundary problems of inviscid compressible fluids with the non-expanding background solutions.
Exponentially fast convergence when the vacuum state is involved  discovered in this paper is a new feature of the problem studied.
 \end{abstract}


 \tableofcontents

\section{Introduction}

Consider the following vacuum free boundary problem for
the compressible gravity-driven Euler equations with frictional damping:
\begin{subequations}\label{eq-ed}
\begin{align}
& \pl_t \rho   + {\rm div}(\rho   u ) = 0 &  {\rm in}& \ \ \Omega(t)=\mathbb{T}^{n-1}\times\lt(0,\ \Gamma(t,x_*)\rt), \label{eq-ed-a}\\
 &  \pl_t  (\rho   u )   + {\rm div}(\rho   u \otimes   u ) +\na  p = -\rho u -\rho g e_n& {\rm in}& \ \ \Omega(t),\label{eq-ed-b}\\
 &\partial_t\Gamma=u_n-\sum_{1\le i\le n-1}u_i\partial_{x_i}\Gamma   &    {\rm on}& \  \{x_n= \Gamma(t,x_*)\}, \label{eq-ed-c}\\
  & \rho=0   &    {\rm on}& \  \{x_n= \Gamma(t,x_*)\},  \label{eq-ed-d}\\
 & u_n=0  &    {\rm on}& \  \{x_n= 0\}, \label{eq-ed-e}\\
&(\rho,{ u})=(\rho_0, { u}_0) & {\rm on} & \ \   \Omega(0)=\mathbb{T}^{n-1}\times\lt(0, \ \Gamma(0,x_*)\rt).  \label{eq-ed-f}
 \end{align} \end{subequations}
Here $t\ge 0$, $x=(x_*, x_n)\in \mathbb{R}^n$, $\rho$, $u$, and  $p$
denote, respectively, the time variable,  space variable, density, velocity, and pressure;
$\mathbb{T}^{n-1}$, $g$, and $e_n$ represent, respectively, the torus which can be thought of as the unit volume  with periodic boundary conditions,  acceleration of gravity, and   $n$-th unit vector; $\Oa(t)\subset \mathbb{R}^n$, and $\Gamma(t,x_*):  \mathbb{R}^{+}\times \mathbb{T}^{n-1} \to \mathbb{R}$ describe, respectively, the changing volume occupied by the fluid at time $t$, and moving vacuum boundary.
We are concerned with the polytropic gas for which the equation of state is  given by
$$ p(\rho)=\rho^{\gamma}, \  \ {\rm where} \ \  \gamma>1 {~\rm is~the~adiabatic~exponent}. $$
From the statistical consideration and simplification of the complicated microscopic flow picture,  the system of compressible Euler equations with damping can be derived from the  Navier-Stokes equations for compressible fluids  (cf. \cite{whitham}). The system of \ef{eq-ed-a} and \ef{eq-ed-b}  is also related to the shallow water model for which $\gamma=2$  (cf. \cite{stoker, whitham}).

Let $c(\rho)=\sqrt{ p'(\rho)}$ be the sound speed, the condition
\be\label{physical} -\infty<\nabla_\mathcal{N}\lt(c^2(\rho)\rt)<0  \  \ {\rm on} \ \  \{x_n=\Gamma(t,x_*) \}
\ee
defines  a  physical vacuum boundary  (cf. \cite{7,10',16',23,24,25}), which is also called a vacuum boundary with physical singularity in contrast to the case that $ \nabla_\mathcal{N}\lt(c^2(\rho)\rt)=0$ on  $\{x_n=\Gamma  \}$,
 where
$\mathcal{N}$ is the exterior unit normal vector to $\{x_n=\Gamma \}$ defined by
$\mathcal{N}= \frac{(-\pl_{x_1}\Gamma, \cdots, -\pl_{x_{n-1}} \Gamma, 1)}{\sqrt{1+|\pl_{x_1}\Gamma|^2 + \cdots + |\pl_{x_{n-1}}\Gamma|^2}}.$
To  capture the physical singularity, we assume that the initial density satisfies that
\be\label{initial density}\begin{split}
& \rho_0>0 \ \ {\rm in} \ \ \Omega(0), \ \   \int_{\Omega(0)} \rho_0(x) dx =M, \ \
\\
&  \rho_0=0 \ \ {\rm and} \ \ -\infty<\nabla_\mathcal{N}\lt(c^2(\rho_0)\rt)<0  \  \ {\rm on} \ \  \{x_n=\Gamma(0,x_*) \},
\end{split}\ee
where $M\in (0, \iy)$ is the initial total mass.
Restricting the momentum equation \ef{eq-ed-b} on $\{x_n=\Gamma  \}$, we have that
$$D_t u \cdot \mathcal{N} = - (\ga-1)^{-1} \nabla_\mathcal{N}\lt(c^2(\rho)\rt) - u\cdot \mathcal{N}-g e_n\cdot \mathcal{N}, $$
where $D_t u= (\pl_t + u\cdot \nabla)u$ is the acceleration of the moving vacuum boundary,
which is due to  three parts: the physical vacuum singularity, frictional damping and  gravity. The balance of these three effects is crucial to the global-in-time existence of  smooth solutions to the vacuum free boundary problem \ef{eq-ed} and \ef{initial density}, which is the issue we will address in this paper.

The vacuum free boundary problem \ef{eq-ed} admits a stationary solution with total mass $M$ given by
\begin{align}\label{1.25}
\bar\rho(x)=\lt(\nu(\hbar-x_n) \rt)^{{1}/({\ga-1})},  \ \  \bar u=0, \ \
\Oa=\mathbb{T}^{n-1} \times (0,\hbar),
\end{align}
where $\nu$ and $\hbar$ are positive constants determined by
\be\label{1.25-2}
\nu = g \ga^{-1} (\ga-1) \ \  {\rm and  } \ \ \hbar = \ga(\ga-1)^{-1} g^{-1} (M g)^{(\ga-1)/\ga}.
\ee
It should be noted that
$ \na_{e_n} {\bar\rho^{\ga-1}}= -\nu$  on $\{x_n=\hbar\}$, which means that
 the stationary solution satisfies the physical vacuum condition \ef{physical}. This motivates the study of the vacuum free boundary problem with physical singularity in the present work.

The aim of this paper is to prove the global-in-time  existence of smooth solutions to the physical vacuum free boundary problem \ef{eq-ed} and
\ef{initial density}  and the large time asymptotic behavior of solutions in space dimensions $n=1, 2, 3$, which is challenging  due to the strong degeneracy of the system near vacuum states. Indeed,  the problem  is a degenerate and characteristic hyperbolic system violating the uniform Kreiss-Lopatinskii condition (cf. \cite{17}) due to resonant wave speeds at vacuum boundaries, so that the standard approach of symmetric hyperbolic systems (cf. \cite{Friedrichs, Kato,17}) do not apply.
As realized in \cite{10', 16'}, the appearance of density functions as coefficients in a nonlinear
wave equation governing the propagation of  acoustic waves shows that this wave equation loses derivatives
with respect to the non-degenerate case of  compressible fluids.
In fact, a satisfactory theory of local-in-time  well-posedness theory is quite recent \cite{16,10, 7, 10', 16'}. (See also \cite{zhenlei,zhenlei1, LXZ,serre} for related works on the local theory.)   It is of fundamental importance to extend  the local-in-time existence theory to the long-time one  for nonlinear problems.
The current available results for the equations of inviscid compressible fluids with physical vacuum are mainly for the cases of expanding vacuum boundaries.
For instance,  the  global-in-time existence  of  smooth solutions of which the expanding rate of vacuum boundaries is linear, $O(1+t)$, when the initial data are small perturbations of affine motions (cf. \cite{sideris1,sideris2}), is proved in  \cite{HaJa1,ShSi} for the $n$-dimensional ($n=1,2,3$) compressible Euler equations.  (See also \cite{CHJ,HaJa2,PHJ} for related results of this type.) For the $n$-dimensional Euler equations with damping (but without gravity),
the global- or almost global-in-time existence  of  smooth solutions with the sub-linear, $O(1+t)^{1/(n(\gamma-1)+2)}$, expanding rate of vacuum boundaries is shown in \cite{LZ,HZ,HZeng}, when the initial data are small perturbations of Barenballt
self-similar solutions (cf. \cite{ba}) of the corresponding porous media equation simplified via Darcy's law.
Indeed,
a family of particular solutions capturing the physical vacuum singularity
was constructed in \cite{23}, where the time-asymptotical equivalence of those solutions to Barenblatt solutions with the same total masses was shown.
Due to the construction in \cite{23}, it became a natural question to investigate the  long-time existence of smooth solutions and their time-asymptotic equivalence for the vacuum free boundary problem capturing  \ef{physical} for the compressible Euler equations with damping.
This question was answered in \cite{LZ} for the one-dimensional case and in \cite{HZ} for the three-dimensional case with spherical symmetry, and partially answered in \cite{HZeng} for the three-dimensional case without symmetry assumption. (See also \cite{HZ2} for the precise time asymptotics of vacuum boundaries.)

The problem of the global-in-time existence of solutions to \ef{eq-ed} studied in the current work differs from those mentioned above  for expanding vacuum boundaries (cf. \cite{LZ,HZ,HZeng,HaJa1,ShSi}).
Fluid expansions have stabilizing effects due to the dispersion which was used in the proof of the results just mentioned, and also plays the important role in the analysis in other context, for example, in general relativistic cosmological models (cf. \cite{Rodnianski, Oliynyk}).
In the presence of gravity, problem \ef{eq-ed} admits a stationary solution, and the global-in-time solutions to the  problem \ef{eq-ed} and
\ef{initial density} are perturbations of the stationary solution. The mechanism of the stability for this problem is the balance of the physical singularity to pushing the vacuum boundary,  the gravity to deceleration the vacuum boundary, and the dissipation of the frictional damping.

The system of \ef{eq-ed-a} and \ef{eq-ed-b} can be observed from both the hyperbolic and parabolic points of view. When we treat the damping term $-\rho u$ as a lower order term, the system is of hyperbolic nature. However, when we view the acceleration term $ (\pl_t + u\cdot \nabla)u$ as a fast decay  perturbation in large time in the spirit of Darcy's law, ignoring this term, the system becomes
 \be\label{darcy}  \rho u=-\nabla p(\rho)-\rho g e_n,  \ \ \ \
  \rho_t=\Delta p(\rho)+g\pl_{x_n} \rho, \ee
which is of degenerate parabolic nature. The approximation for large time can be seen formally  by the following rescaling ${ x}'=\ea { x}$, $t'=\ea^2 t$, and ${ u}'= { u}/ \ea$.  It   should be noted that system \ef{darcy} has the same stationary solution as problem \ef{eq-ed}. This indicates that the
long time dynamics of the original problem should be governed by that of \ef{darcy}. This point of view is adopted in the proof of main theorems in this paper.

We prove the global-in-time existence of smooth solutions to problem \ef{eq-ed} and \ef{initial density} in space dimensions $n=1, 2, 3$ for initial data being small perturbations of the stationary solution \ef{1.25}, and investigate the large time asymptotic behavior towards the stationary solution.
We start with the one-dimensional case of $n=1$, where we obtain the exponentially fast convergence of the density, velocity, and vacuum boundary to those of the stationary solution. Then we investigate the multi-dimensional case of $n=2$ and $n=3$.
The situation becomes much more involved due to the intricate time evolution of the vacuum boundary geometry and the complicated multi-dimensional motions in both tangential and normal directions.
This is handled by the weighted energy estimates of both the hyperbolic and parabolic types, together with  the elliptic estimates near the bottom and the curl estimates.

 The analysis in this paper is different from that for the case when vacuum boundaries expand in time as studied in  \cite{HaJa1,ShSi, HZeng}.  Let $\bar x(t, y)$ be the particle position at time $t$ starting from the initial position $y$,
the background flows discussed in \cite{HaJa1,ShSi, HZeng} are typically of the form $\bar x(t, y)= \alpha(t) y$ with
$y$ in a ball, and with $\alpha(t)\sim 1+t$
for   affine solutions  in \cite{HaJa1,ShSi} and $\alpha(t)\sim (1+t)^{1/(n(\gamma-1)+2)}$ for the Barenblatt solution in \cite{HZeng}.
Let $x=x(t,y)$ be the flow of the velocity field $u$: $\pl_t x(t,y)= u(t, x(t,y))$ for $t>0$ with
$x(0,y)=x_0(y)$ for $y\in \Omega$, where $\Omega$ is the initial domain of the corresponding background solution.
Due to the form of $\bar x(t,y)=\alpha(t)y$, it is natural to decompose
$$x(t,y)=\bar x(t,y)+\alpha(t) \omega(t,y)=\alpha(t) \lt(y+\omega(t,y)\rt),$$
and the problems are reduced to study the equations for the unknown $\omega(t,y)$.
The higher-order norms defined in \cite{HaJa1,ShSi, HZeng} are of the form
\be\label{a} \mathcal{E}(t)=
\alpha(t)^{\beta}\mathcal{E}_I(t)
+\mathcal{E}_{II}(t)\ee
for some constant $\beta>0$, where $\mathcal{E}_I(t)$ contains the weighted $L^2$-norms of   space-time mixed derivatives  of $\pl_t \omega$,
and $\mathcal{E}_{II}(t)$ contains the  weighted $L^2$-norms only involving  the space derivatives of $\omega$.  The approaches adopted in \cite{HaJa1,ShSi, HZeng} to obtain the bounds for $\mathcal{E}(t)$ are to derive the following inequalities:
$$  \alpha(t)^{\beta}\mathcal{E}_I(t) \le C_1\mathcal{E}(0)+ F_1(\mathcal{E}(t)) \ \ {\rm and} \  \ \mathcal{E}_{II}(t)\le  C_2 \mathcal{E}(0) +F_2(\mathcal{E}(t)) , $$
for some positive constants  $C_1$, $C_2$ and functions $F_1$, $F_2$.
It should be noted that the decay for $\pl_t \omega$  and its derivatives is naturally encoded in the definition of $\mathcal{E}(t)$ due to the growth of $\alpha(t)$ in time for solutions with expanding vacuum boundaries.
However, for the problem studied in the present paper, the background solution is stationary, of the form  of $\bar x(y, t)=y$, so
we do not have a natural growing factor $\alpha(t)$.
Instead, we will have to identify the decay of $\pl_t \oa$
and its derivatives first. By the same notations, we may decompose
$$x(t,y)=y+\omega(t,y) $$
and define the higher-order norm as
 \be\label{b} \mathcal{E}(t)=\mathcal{E}_I(t)+\mathcal{E}_{II}(t).\ee
Compared \ef{b} with \ef{a}, we do not have a growing factor of $\alpha(t)$ in front of $\mathcal{E}_I(t)$.
Therefore, we have to identify the decay property of $\mathcal{E}_I(t)$ first.
Indeed, we  identify a number $\delta>0$ depending only on the adiabatic exponent $\gamma$ and the total mass $M$ such that
$$ e^{\da t}\mathcal{E}_I(t)\le C \mathcal{E}_I(0)  $$
for some positive constant $C$, based on which, we prove the bound  of $\mathcal{E}_{II}(t)$.
The estimates for $\mathcal{E}_I(t)$ and $\mathcal{E}_{II}(t)$ are in order, the estimate for $\mathcal{E}_{II}(t)$ is based on that for $\mathcal{E}_{I}(t)$, so we have to make the distinction between $\mathcal{E}_I(t)$ and $\mathcal{E}_{II}(t)$ in the estimates.
This is different from those for expanding solutions in \cite{HaJa1,ShSi, HZeng}, where $\mathcal{E}_I(t)$ and $\mathcal{E}_{II}(t)$ can be estimated simultaneously.

It should be remarked that the exponential decay mentioned above is a new feature of the problem studied in this paper, by noting that the decay of the velocity field for the problems with expanding vacuum boundaries is algebraic. Indeed, for problems involving the vacuum state, even with the higher-order dissipation such as viscosity, the currently available results on the decay are of algebraic rates, see, for instance, \cite{LXZ1, LXZ2, OZ} for the related results.

Besides the results mentioned above, we review  some works related to the study of vacuum states. The theoretical study of vacuum states of gas dynamics dates back to 1980 when it was shown in \cite{LiuSmoller} that shock waves vanish at the vacuum.  The well-posedness of smooth solutions with sound speed $c(\rho)$ smoother than $C^{ {1}/{2}}$-H$\ddot{\rm o}$lder continuous near vacuum states for compressible inviscid fluids can be found in \cite{chemin1,chemin2,24,25,MUK,Makino,38,39}.
The phenomena of physical vacuum arise additionally in several important situations, for example, the equilibrium and dynamics of boundaries of gaseous stars (cf. \cite{6', cox, HaJa2, LXZ}).
There have been extensive studies on viscous flows with vacuum. Since the ideas and techniques are quite different from the inviscid flow, so we omit the discussions on this topic here.

{\bf Notation} Throughout the rest of this paper, $C$ will denote a positive constant which only depends on the parameters of the problem, the adiabatic exponent $\ga$ and the initial total mass $M$, but does not depend on the data. They are referred as universal and can change from one inequality to another one. We will adopt the notation $ a \les b$ to denote $a \le C b$, and $a \thicksim b$ to denote $C^{-1} b\le a \le Cb$, where $C$ is the universal constant as defined
above. We will use $\int=\int_\Oa$,
$\|\cdot\|_{W^{k,q}}=\|\cdot\|_{W^{k,q}(\Oa)}$
for any constants $k\ge 0$ and $q\ge 1$, and $\|\cdot\|_{H^s}=\|\cdot\|_{H^s(\Oa)}$ for any constant $s\ge 0$. In particular, $\Oa=I=(0,\hbar)$ in one space dimension ($n=1$).

\section{The one-dimensional motions}
In one space dimension, the vacuum free boundary problem \ef{eq-ed} reads
\begin{subequations}\label{eq-ed-1d}
\begin{align}
&\pl_t \rho + \pl_x \lt(\rho u\rt) =0         & {\rm in} &\ \  (0, {\Gamma}(t)),
\label{eq-ed-1d-a}\\
&\rho \lt( \pl_t u  +  u \pl_x u \rt)  +\pl_x \rho^\ga  =-\rho u -\rho g  & {\rm in } &\ \ (0, {\Gamma}(t)),  \label{eq-ed-1d-b}\\
&   \frac{d}{dt} {\Gamma}(t)= u  (t,\Gamma(t)),                & &\\
&\rho(t,\Gamma(t))=0  \  \ {\rm and} \ \   u(t,0)=0  ,       &                       &   \\
&(\rho, u)=(\rho_0, u_0)                 & {\rm on} & \ \ (0,\Gamma(0)),
 \end{align}\end{subequations}
and the statioanry solution \ef{1.25} becomes
\be\label{9-24}
\bar\rho(x)=\lt( \nu(\hbar-x) \rt)^{\frac{1}{\ga-1}},  \ \  \bar u=0, \ \
\Oa={I}=(0,\hbar),
\ee
where $\nu$ and $\hbar$ are positive constants given by \ef{1.25-2}.

We transform the free boundary problem \ef{eq-ed-1d} into Lagrangian variables to fix the boundary, and choose the interval of the stationary solution, $I$, as the reference interval.
For $y\in I$, we define $x$ as the Lagrangian flow of the velocity $u$ by
\begin{align}\label{lag-1d}
\pl_t x(t,y)= u(t, x(t,y)) \   {\rm for} \   t>0 \ ,  \  {\rm and} \   x(0,y)=x_0(y) ,
\end{align}
where $x_0: {I}\rightarrow(0,\Gamma(0))$ is a diffeomorphism defined by
$\int_0^{x_0(y)} \rho_0(z) dz=\int_0^y \bar\rho(z) dz$. Clearly,
$
\rho_0(x_0(y))x_0'(y)=\bar\rho(y)
$.
Then, it follows from
\ef{eq-ed-1d-a} and \ef{lag-1d} that
\be\label{9-25-1}
\rho(t,x(t,y))=\frac{\rho_0(x_0(y)) x_0'(y)}{\pl_y x(t,y)}= \frac{\bar \rho(y)}{\pl_y x(t,y)}.
\ee
Substituting \ef{9-25-1} into \ef{eq-ed-1d-b}, and using \ef{lag-1d}, we  rewrite  problem \ef{eq-ed-1d} as
\begin{subequations}\label{9-25-2}
\begin{align}
&\bar\rho \lt(\pl_t^2 x + \pl_t x\rt) + \pl_y \lt(\frac{\bar\rho}{\pl_y x}\rt)^\ga=- \bar\rho g         & {\rm in} &\   (0, T]\times I ,
\label{9-25-2-a}\\
&\pl_t x=0        & {\rm on} &\ (0, T] \times  \{y=0\},
\label{9-25-2-b}\\
&(x,\pl_t x)=(x_0, u_0(x_0))        & {\rm on} &\   \{t=0\} \times I\label{9-25-2-c}.
\end{align}\end{subequations}
In the setting, the moving vacuum boundary for problem \ef{eq-ed-1d} is given by
$
\Gamma(t)=x(t,\hbar)$ for $t> 0$.

\subsection{Main results}
We define the perturbation $\oa$ by
$\oa(t,y) = x(t,y)- y
$,
and set $\sigma(y) = \bar\rho^{\ga-1}(y)$
and $
\iota= ({\ga-1})^{-1}$.
We introduce that for nonnegative integers $m$ and $i$,
\begin{align*}
\mathcal{E}^{m,i}(t)= \lt\| \sa^{\frac{\iota+i}{2}}  \pl_t^{m+1}  \pl_y^i \oa \rt\|_{L^2}^2+ \lt\| \sa^{\frac{\iota+i}{2}} \pl_t^{m}   \pl_y^i \oa \rt\|^2_{L^2} +   \lt\| \sa^{\frac{\iota+i+1}{2}}
  \pl_t^m   \pl_y^{i+1} \oa \rt\|^2_{L^2} ,
\end{align*}
and define the higher-order weighted Sobolev norm
$\mathcal{E}(t)$ by
\begin{align*}
\mathcal{E}(t) =  \sum_{0\le m+i\le [\iota]+4} \mathcal{E}^{m,i}(t)
.
\end{align*}

The main results in this section are stated as follows.
\begin{thm}\label{thm-1d}
There exists a positive constant $\bar \ea>0$ such that if $\mathcal{E}(0)\le \bar\ea$, then problem \ef{9-25-2} admits a global  smooth solution in $[0,\iy)\times I$ satisfying
\be\label{5.14-1}
\mathcal{E}(t) \le  C e^{-\da t} \mathcal{E}(0) \ \ {\rm for} \ t \ge 0,
\ee
where $C$ and $\da$ are positive constants which only depend on the adiabatic exponent $\ga$ and the initial total mass $M$, but do not depend  on the time $t$.
\end{thm}

As a corollary of Theorem \ref{thm-1d}, we have the following theorem for solutions to the original vacuum free boundary problem \ef{eq-ed-1d} concerning the convergence of the vacuum boundary $\Gamma$, density $\rho$ and velocity $u$ to those of the stationary solution.

\begin{thm}\label{thm-1d'}
There exists a positive constant $\bar \ea>0$ such that if $\mathcal{E}(0)\le \bar\ea$, then problem \ef{eq-ed-1d} admits a global  smooth solution $(\rho,u,\Gamma(t))$ for $t\in [0,\iy)$ satisfying
\begin{subequations}\label{5.14}
\begin{align}
&|\rho(t,x(t,y))-\bar\rho(y)|\le C (\hbar-y)^{{1}/({\ga-1})} \sqrt{e^{-\da t} \mathcal{E}(0)}, \\
&|u(t,x(t,y))|+|\Gamma(t)-\hbar| + \sum_{1\le m\le 3} \lt|\frac{d^m\Gamma(t)}{d t^m } \rt| \le C \sqrt{e^{-\da t} \mathcal{E}(0)},
 \end{align}\end{subequations}
 for all $y\in I$ and $t\ge 0$. Here
 $C$ and $\da$ are positive constants which only depend on the adiabatic exponent $\ga$ and the initial total mass $M$, but do not depend on the time $t$.
\end{thm}

\subsection{Proof of Theorems \ref{thm-1d} and \ref{thm-1d'}}
We first prove  Theorem \ref{thm-1d}.
The proof of the global existence of smooth solutions is based on the local existence of smooth solutions (cf. \cite{10,16}), and the following a priori estimates stated in Proposition \ref{5.13} whose proof will be presented in the next subsection.
\begin{prop}\label{5.13}
Let $x(t,y)=y+\oa(t,y)$ be a solution to problem \ef{9-25-2}
in the time interval $[0,T]$ satisfying the a priori assumption:
\begin{align}\label{4.21-1}
\mathcal{E}(t) \le \ea_0^2, \ \ t\in[0,T],
\end{align}
for some suitably small fixed positive number $\ea_0$ independent of $t$. Then there exist positive constants $C$ and $\da$,  only depending  on $\ga$ and $M$,  but not on $t$,  such that
\bee
\mathcal{E}(t) \le C e^{-\da t} \mathcal{E}(0), \ \ t\in [0, T].
\eee
\end{prop}

To prove Theorem  \ref{thm-1d'}, we recall the following estimate  achieved in Lemma 3.7 of \cite{LZ}.
\begin{align}\label{4.21-2'}
\sum_{0\le m+i \le [\iota]+4}\lt\|\sa^{ \max\lt\{0,\ \frac{m+2i-3}{2} \rt\}} \pl_t^m\pl_y^i \oa \rt\|_{L^\iy}^2
 \les \mathcal{E}(t),
\end{align}
provided that $\mathcal{E}(t)$ is finite. Indeed, the proof of  \ef{4.21-2'} is based on the weighted Sobolev embedding \ef{wsv} and the Hardy inequality \ef{hard}. Therefore,  \ef{5.14} follows from  \ef{5.14-1}, \ef{4.21-2'}, and
\begin{align*}
&\rho(t,x(t,y))-\bar\rho(y)=-\bar\rho(y)
\frac{\pl_y \oa(t,y)}{1+\pl_y \oa(t,y)},\\
&u(t,x(t,y))=\pl_t \oa(t,y), \ \ \Gamma(t)-\hbar=\oa (t,\hbar).
\end{align*}

\subsection{A priori estimates}
This subsection devotes to proving Proposition \ref{5.13}. To this end,
we introduce the following Sobolev norms for the  dissipation:
\begin{align*}
\mathcal{D}^{m,i}(t)= & \lt\| \sa^{\frac{\iota+i}{2}}  \pl_t^{m+1}  \pl_y^i \oa \rt\|_{L^2}^2 +   \lt\| \sa^{\frac{\iota+i+1}{2}}
  \pl_t^m   \pl_y^{i+1} \oa \rt\|^2_{L^2},
\end{align*}
where $m$ and $i$ are nonnegative integers. To specify the behavior of solutions near the bottom and the top, we divide
$\mathcal{D}^{m,i}$
into two parts $\mathcal{D}_1^{m,i}$ and $\mathcal{D}_2^{m,i}$ as follows:
\begin{align*}
 &  \mathcal{D}_1^{m,i}(t)=  \lt\| \za_1  \pl_t^{m+1}  \pl_y^i \oa \rt\|_{L^2}^2 +   \lt\|\za_1
  \pl_t^m   \pl_y^{i+1} \oa \rt\|^2_{L^2},\\
 & \mathcal{D}_2^{m,i}(t)=  \lt\| \za_2 \sa^{\frac{\iota+i}{2}}  \pl_t^{m+1}  \pl_y^i \oa \rt\|_{L^2}^2 +   \lt\| \za_2 \sa^{\frac{\iota+i+1}{2}}
  \pl_t^m   \pl_y^{i+1} \oa \rt\|^2_{L^2},
\end{align*}
where $\zeta_1=\zeta_1(y)$ and $\zeta_2=\zeta_2(y)$ are  smooth cut-off functions satisfying
\begin{subequations}\label{4.13}\begin{align}
&\zeta_1=1 \ \ {\rm on} \ \ [0, \hbar/{2}], \ \  \zeta_1=0 \ \ {\rm on} \ \  [3\hbar/{4}, \hbar], \ \  \zeta_1' \le 0 \ \  {\rm on} \ \ [0,\hbar], \label{David-1}\\
&\zeta_2=0 \ \ {\rm on} \ \ [0, \hbar/{4}], \ \  \zeta_2=1 \ \ {\rm on} \ \  [ \hbar/{2},\hbar], \ \  \zeta_2' \ge 0 \ \ {\rm on} \ \  [0,\hbar]. \label{David}
\end{align}\end{subequations}
We will use elliptic estimates to bound $\mathcal{D}_1^{m,i}$ with $i\ge 1$ in Section \ref{s2.3.2}, and energy estimates to bound $\mathcal{D}^{m,0}$, and
$\mathcal{D}_2^{m,i}$ with $i \ge 1$ in Section \ref{s2.3.3}. For this purpose, we rewrite equation \ef{9-25-2-a}  as
\be\label{problem-1d-a}
   \sa^\iota \lt(\pl_t^2 \oa + \pl_t \oa \rt) + \pl_y \lt( \sa^{\iota+1}  \lt( \lt(1+{\pl_y \oa}\rt)^{-\ga}- 1 \rt) \rt)=0   .
  \ee

\subsubsection{Preliminaries}
It follows from \ef{4.21-1} and \ef{4.21-2'} that for $t\in [0,T]$,
\begin{align}\label{4.21-2}
\sum_{0\le m+i \le [\iota]+4}\lt\|\sa^{ \max\lt\{0,\ \frac{m+2i-3}{2} \rt\}} \pl_t^m\pl_y^i \oa \rt\|_{L^\iy}^2
 \les  \mathcal{E}(t) \le  \ea_0^2,
\end{align}
which implies, with the aid of the smallness of $\ea_0$,  that
\bee\label{5.14-3}
|\pl_y \oa(t,y)|\le 1/2 \ \ {\rm for} \ \  (t,y)\in [0,T]\times I.
\eee
Based on this, it is easy to see that for integers $m\ge 0$, $i\ge 0$ and  $l \ge -1$,
\begin{align}
\lt|\pl_t^m   \pl_y^i \lt(1+{\pl_y \oa}\rt)^{-\ga-l}  \rt|
\les  \mathcal{I}^{m,i},  \label{4.21-5}
\end{align}
where $\mathcal{I}^{m,i}$ are defined inductively as follows:
\begin{align*}
& \mathcal{I}^{0, 0}=1, \\
& \mathcal{I}^{m, i} =  \sum_ {\substack{0\le r\le m,   \ 0\le k \le i
\\ 0\le r +k \le m +i-1
  }} \mathcal{I}^{r ,k} \lt|\pl_t^{m-r} \pl_y^{i-k+1}\oa\rt|
.
 \end{align*}

\begin{lem}\label{5.14-4} It holds that for $t\in [0,T]$,
\begin{subequations}\label{4.28}\begin{align}
&\sum_{
1\le m+i \le [\iota]+3}\lt\|\sa^{\frac{m+2i-1}{2}}\mathcal{I}^{m,i}
\rt\|_{L^\iy}^2  \les \mathcal{E}(t) \le \ea_0^2 , \label{4.28-1}\\
& \sum_{1\le
 m+i \le [\iota]+3}\lt\|\sa^{\frac{\iota+i}{2}}\mathcal{I}^{m,i}
\rt\|_{L^2}^2 \les \sum_{1\le r+k\le m+i+1}\mathcal{D}^{r,k}(t),\label{4.28-2}\\
&\sum_{
1\le m+i \le  [\iota]+4}\lt\|\sa^{\frac{\iota+i+1}{2}}\mathcal{I}^{m,i}
\rt\|_{L^2}^2 \les  \mathcal{D}^{m,i}(t)+\sum_{1\le r+k\le m+i-1}\mathcal{D}^{r,k}(t).\label{4.28-3}
\end{align}\end{subequations}
\end{lem}

{\em Proof}. It follows from \ef{4.21-2} that for $1\le m+i \le [\iota]+3$,
\begin{align*}
&\sa^{\frac{m+2i-1}{2}}\mathcal{I}^{m,i}
=\lt|\sa^{\frac{m+2i-1}{2}}\pl_t^m \pl_y^{i+1}\oa\rt|
\\
&\quad +\sa^{\frac{1}{2}}\sum_ {\substack{r\le m,   \ k \le i
\\ 1\le r +k \le m +i-1
  }}  \sa^{\frac{ r+2k-1}{2}} \mathcal{I}^{r ,k} \lt| \sa^{\frac{m- r+2(i- k)-1}{2}}\pl_t^{m-r} \pl_y^{i+1-k}\oa\rt|
  \\
&  \les \sqrt{\mathcal{E}} + \ea_0 \sum_ {\substack{r\le m,   \ k \le i
\\ 1\le r +k \le m +i-1
  }}  \sa^{\frac{ r+2k-1}{2}} \mathcal{I}^{r ,k} ,
\end{align*}
which, together with the mathematical induction, proves \ef{4.28-1}.

In view of \ef{4.28-1} and \ef{4.21-2}, we see that for $1\le m+i\le [\iota]+3$,
\begin{align*}
&\sa^{\frac{\iota+i}{2}}\mathcal{I}^{m,i}
\les \sa^{\frac{\iota+i}{2}}\lt|\pl_t^m \pl_y^{i+1}\oa\rt|
\\
&\quad +\ea_0 \sum_ {\substack{r\le m,   \ k \le i
\\ 1\le r +k \le m +i-2
  }} \sa^{\frac{\iota+i+1-r-2k}{2}}\lt| \pl_t^{m-r} \pl_y^{i+1-k}\oa\rt|\\
& \quad + \ea_0  \sum_ {\substack{r\le m,   \ k \le i
\\ 1\le r +k = m +i-1
  }} \sa^{\frac{\iota+k}{2}} \mathcal{I}^{r,k}.
\end{align*}
It follows from \ef{hard} that
$$
 \lt\|\sa^{\frac{\iota+i}{2}}\pl_t^m \pl_y^{i+1}\oa\rt\|_{L^2}^2 \les
 \sum_{l=0,1} \lt\|\sa^{\frac{\iota+i}{2}+1}\pl_t^m \pl_y^{i+1+l}\oa\rt\|_{L^2}^2
 \les  \sum_{l=0,1} \mathcal{D}^{m,i+l};
$$
and that for $r\le m$,   $ k \le i$, $ 1\le r +k \le m +i-2$,
\begin{align*}
&\lt\| \sa^{\frac{\iota+i+1-r-2k}{2}} \pl_t^{m-r} \pl_y^{i+1-k}\oa\rt\|_{L^2}^2 \\
  \les
  &
  \sum_{0\le l\le r+k} \lt\| \sa^{\frac{\iota+i+1-r-2k}{2}+r+k} \pl_t^{m-r} \pl_y^{i+1-k+l}\oa\rt\|_{L^2}^2\\
  \les &   \sum_{0\le l\le r+k} \lt\| \sa^{\frac{\iota+i+1-k+l}{2}} \pl_t^{m-r} \pl_y^{i+1-k+l}\oa\rt\|_{L^2}^2 \\
  \le & \sum_{0\le l\le r+k} \mathcal{D}^{m-r,i-k+l}
  ,
\end{align*}
because
$\iota+i+1-r-2k \ge \iota+3-(m+i)\ge 0$.
So, one has that
\begin{align*}
\lt\|\sa^{\frac{\iota+i}{2}}\mathcal{I}^{m,i}
\rt\|_{L^2}^2 \les \sum_{1\le r+k\le m+i+1}\mathcal{D}^{r,k}+ \ea_0^2  \sum_ {\substack{r\le m,   \ k \le i
\\ 1\le r +k = m +i-1
  }} \lt\| \sa^{\frac{\iota+k}{2}} \mathcal{I}^{r,k}\rt\|_{L^2}^2
\end{align*}
for $1\le m+i\le [\iota]+3$.
This, together with the mathematical induction, proves \ef{4.28-2}.

Similarly, we can prove \ef{4.28-3} by noticing that
\begin{align*}
&\sa^{\frac{\iota+i+1}{2}}\mathcal{I}^{m,i}
\les \sa^{\frac{\iota+i+1}{2}}\lt|\pl_t^m \pl_y^{i+1}\oa\rt|
\\
&\quad +\ea_0 \sum_ {\substack{r\le m,   \ k \le i
\\ 1\le r +k \le m +i-2
  }} \sa^{\frac{\iota+i+2-r-2k}{2}}\lt| \pl_t^{m-r} \pl_y^{i+1-k}\oa\rt|\\
& \quad + \ea_0  \sum_ {\substack{r\le m,   \ k \le i
\\ 1\le r +k = m +i-1
  }} \sa^{\frac{\iota+k+1}{2}} \mathcal{I}^{r,k},
\end{align*}
and that
for $r\le m$,   $k \le i$, $ 1\le r +k \le m +i-2$,
\begin{align}
&\lt\| \sa^{\frac{\iota+i+2-r-2k}{2}} \pl_t^{m-r} \pl_y^{i+1-k}\oa\rt\|_{L^2}^2
\notag\\
  \les
  &
  \sum_{0\le l\le r+k-1} \lt\| \sa^{\frac{\iota+i+2-r-2k}{2}+r+k-1} \pl_t^{m-r} \pl_y^{i+1-k+l}\oa\rt\|_{L^2}^2
  \notag \\
  \le & \sum_{0\le l\le r+k-1} \mathcal{D}^{m-r,i-k+l}
  .\label{5.5-5}
\end{align}
This finishes the proof of Lemma \ref{5.14-4}. \hfill $\Box$

\subsubsection{Elliptic estimates near the bottom}\label{s2.3.2}

\begin{lem}
It holds that for   $m+i \le 4+[\iota]$ and $i \ge 1$, and $t\in [0,T]$,
\be\label{4.13-1}
\mathcal{ D}^{m, i}_1(t) \lesssim \sum_{0\le r\le m+i}\mathcal{D}_1^{r,0} (t)  \les \sum_{0\le r\le m+i}\mathcal{D}^{r,0} (t) .
\ee

\end{lem}

{\em Proof}.
It follows from the Sobolev embedding and the a priori assumption \ef{4.21-1} that
\begin{align}\label{4.21-3}
\lt\|\pl_t^m \oa \rt\|_{W^{[\iota]+4-m, \iy}(I_1)}^2 \les
\lt\|\pl_t^m \oa \rt\|_{H^{[\iota]+5-m}(I_1)}^2
\les  \mathcal{E} \le \ea_0^2,
\end{align}
where $I_1=(0,3\hbar/4)$. This implies that for
$1\le m+i \le [\iota]+3$,
\begin{align*}
&\lt\|\mathcal{I}^{m, i}\rt\|_{L^\iy(I_1)} \le
 \lt\|\pl_t^{m } \pl_y^{i +1}\oa\rt\|_{L^\iy(I_1)}
 \\
&\quad  + \sum_ {\substack{ r\le m,   \ k \le i
\\ 1\le r +k \le m +i-1
  }} \lt\|\mathcal{I}^{r ,k}\rt\|_{L^\iy(I_1)} \lt\|\pl_t^{m-r} \pl_y^{i-k+1}\oa\rt\|_{L^\iy(I_1)}\\
&\quad  \les \sqrt{\mathcal{E}}
+ \ea_0 \sum_ {\substack{ r\le m,   \  k \le i
\\ 1\le r +k \le m +i-1
  }} \lt\|\mathcal{I}^{r ,k}\rt\|_{L^\iy(I_1)} ,
\end{align*}
which, together with the mathematical induction, gives that
\begin{align}\label{4.21-4}
\sum_{1\le m+i \le [\iota]+3}\lt\|\mathcal{I}^{m, i}\rt\|_{L^\iy(I_1)}\les \sqrt{\mathcal{E}} \les  \ea_0.
\end{align}

We divide \ef{problem-1d-a} by $\sa^{\iota}$, and take $\pl_t^m \pl_y^{i-1}$ $(i\ge 1)$ onto the resulting equation to obtain
\begin{align*}
\ga \sa \pl_t^m  \pl_y^{i+1} \oa
=\ga \nu (i-1)  \pl_t^m  \pl_y^{i} \oa+\pl_t^{m+2} \pl_y^{i-1}  \oa + \pl_t^{m+1} \pl_y^{i-1}  \oa  +
\pl_t^{m} \pl_y^{i-1} W,
\end{align*}
which means
\begin{align}\label{4.21-6}
\lt\|\za_1 \pl_t^m  \pl_y^{i+1} \oa\rt\|_{L^2}^2
\les \mathcal{D}_1^{m,i-1} + \mathcal{D}_1^{m+1,i-1}+\lt\|\za_1
\pl_t^{m} \pl_y^{i-1} W\rt\|_{L^2}^2.
\end{align}
Here $\za_1$ is a smooth cut-off function satisfying \ef{David-1}, and
\begin{align*}
W=\ga  \sa \lt(1- \lt(1+{\pl_y \oa}\rt)^{-\ga-1} \rt) \pl_y^2 \oa
- \nu (1+\iota)  \lt( \lt(1+{\pl_y \oa}\rt)^{-\ga}- 1 \rt).
\end{align*}
In view of \ef{4.21-5}, \ef{4.21-3} and \ef{4.21-4}, we see that
\begin{align*}
&\lt\|\za_1
\pl_t^{m} \pl_y^{i-1} W\rt\|_{L^2}^2
\les \lt\|\za_1 |\pl_y \oa|
\pl_t^{m} \pl_y^{i+1} \oa \rt\|_{L^2}^2\\
&\quad+\sum_{r\le m, \ k\le i-1, \ r+k \ge 1}\lt\|\za_1 \mathcal{I}^{r,k}
\pl_t^{m-r} \pl_y^{i+1-k} \oa \rt\|_{L^2}^2\\
&\quad+\sum_{r\le m, \ k\le i-1 }\lt\|\za_1 \mathcal{I}^{r,k}
\pl_t^{m-r} \pl_y^{i-k} \oa \rt\|_{L^2}^2\\
&   \les \ea_0^2 \lt\|\za_1 \pl_t^m  \pl_y^{i+1} \oa\rt\|_{L^2}^2
+ \ea_0^2 \sum_{r\le m, \ k\le i-1, \ r+k \ge 1} \mathcal{D}_1^{m-r,i-k}
\\
& \quad+  \sum_{r\le m, \ k\le i-1 } \mathcal{D}_1^{m-r,i-1-k},
\end{align*}
which implies, with the aid of \ef{4.21-6} and the smallness of $\ea_0$, that
\begin{align*}
&\lt\|\za_1 \pl_t^m  \pl_y^{i+1} \oa\rt\|_{L^2}^2
\les  \mathcal{D}_1^{m+1,i-1}+  \sum_{r\le m, \  k\le i ,\  r+k\le m+i-1} \mathcal{D}_1^{r,k}, \end{align*}
and then
\begin{align}
& \mathcal{D}_1^{m ,i }
\les  \mathcal{D}_1^{m+1,i-1}+   \sum_{r\le m, \  k\le i ,\  r+k\le m+i-1} \mathcal{D}_1^{r,k}.\label{4.21-7}
\end{align}

Based on \ef{4.21-7}, we can use the mathematical induction to prove  \ef{4.13-1}. Clearly, \ef{4.13-1} holds for $m+i=1$ and $i\ge 1$, due to
$\mathcal{D}_1^{0 , 1 }
\les  \mathcal{D}_1^{1,0} +\mathcal{D}_1^{0,0}$.
Suppose \ef{4.13-1} holds for $m+i\le l-1$ and $i\ge 1$, then one has for $m+i=l$ and $i\ge 1$,
\begin{align*}
\mathcal{D}_1^{m ,i }
\les  \mathcal{D}_1^{m+1,i-1}+   \sum_{   r+k\le l-1} \mathcal{D}_1^{r,k} \les \mathcal{D}_1^{m+1,i-1}
+ \sum_{   r \le l-1} \mathcal{D}_1^{r,0} ,
\end{align*}
which proves \ef{4.13-1} from $i=1$ to $i=l$ step by step.
\hfill $\Box$

\subsubsection{Energy estimates}\label{s2.3.3}

\begin{lem}
There exists a small positive constant $\da$ which only depends on $\ga$ and $M$ such that
\begin{align}\label{5.10-3}
  e^{\da t} \mathcal{E}(t) + \int_0^t  e^{\da s} \mathcal{E}(s) ds \les   \mathcal{E}(0), \ \ t\in [0,T] .
\end{align}
\end{lem}

{\em Proof}. The proof consists of three steps, among which
the first two steps devote to deriving energy estimates and the last one showing the decay rate.

{\em Step 1}. In this step, we give the basic energy estimates.
Multiply \ef{problem-1d-a} by $\pl_t \oa$ and integrate the resulting equation over $I$ to get
\be\label{75-1}
\frac{1}{2}\frac{d}{dt}\int  \sa^\iota |\pl_t \oa|^2 dy + \frac{d}{dt}E_1^{0,0}
+\int \sa^{\iota} |\pl_t\oa|^2 dy =0,
\ee
where
$$E_1^{0,0}=\iota \int \sa^{\iota+1} \lt((1+\pl_y \oa)^{1-\ga}-1-(1-\ga)\pl_y\oa\rt)dy.$$
Integrate the product of \ef{problem-1d-a} and $\oa$  over $I$ to obtain
\be\label{75-2}
\frac{1}{2}\frac{d}{dt}\int  \sa^{\iota} \lt(  \oa^2 + 2 \oa \pl_t \oa \rt) dy
+\widetilde{E}_1^{0,0} =
 \int \sa^{\iota} |\pl_t\oa|^2 dy,
\ee
where
$$\widetilde{E}_1^{0,0}=\int  \sa^{\iota+1} \lt(1-(1+\pl_y\oa)^{-\ga}\rt)\pl_y\oa dy.$$
It follows from $2\times\ef{75-1}+\ef{75-2}$ that
\be \label{75-3}
\frac{d}{dt}\mathfrak{E}^{0,0}(t)
+\mathfrak{D}^{0,0}(t)  =0,
\ee
where
\begin{align*}
&\mathfrak{E}^{0,0}=\frac{1}{2} \int  \sa^\iota \lt(2 |\pl_t \oa|^2 +   \oa^2 + 2 \oa \pl_t \oa \rt)  dy + 2E_1^{0,0} ,\\
& \mathfrak{D}^{0,0}=\int \sa^{\iota} |\pl_t\oa|^2 dy + \widetilde{E}_1^{0,0}.
\end{align*}
Due to the Taylor expansion and the Cauchy inequality, one has that
\begin{align}\label{5.8-1}
\mathfrak{E}^{0,0} \thicksim
\mathcal{E}^{0,0} \ \ {\rm and} \  \
 \mathfrak{D}^{0,0} \thicksim
\mathcal{D}^{0,0}.
\end{align}

{\em Step 2}. This step dovotes to showing the higher-order energy estimates.  Let
$m$ and $i$ be nonnegative integers satsifying $m+i\ge 1$, then we divide \ef{problem-1d-a} by $\sa^{\iota}$, and take $\pl_t^m \pl_y^{i}$  onto the resulting equation to obtain
 \begin{align}\label{75-4}
\pl_t^{m+2} \pl_y^i \oa   +
\pl_t^{m+1}  \pl_y^i \oa
  +\sa^{-\iota-i}\pl_y \lt( \sa^{\iota+i+1}\pl_t^m \pl_y^i \lt(1+ \pl_y \oa\rt)^{-\ga}\rt) =0.
   \end{align}
Let $\za_2$ be a smooth cut-off function   satisfying \ef{David} and set
\begin{align*}
\chi=1 \ {\rm when} \  i=0, \ {\rm and} \
\chi=\za_2 \ {\rm when} \ i\ge 1.
 \end{align*}
We integrate the product of   \ef{75-4} and $\chi^2 \sa^{\iota+i}\pl_t^{m+1}  \pl_y^i \oa$ over $I$ and use the boundary condition $\pl_t^{m+1}\oa(t,0)=0$ to get
\begin{align}
&\frac{1}{2}\frac{d}{dt} \int \chi^2  \sa^{\iota+i}\lt|\pl_t^{m+1} \pl_y^i \oa \rt|^2 dy
+\frac{1}{2}\frac{d}{dt} \lt( E_1^{m,i} +E_2^{m,i}
\rt)\notag \\
&+\int \chi^2  \sa^{\iota+i}\lt|\pl_t^{m+1} \pl_y^i \oa \rt|^2 dy = \sum_{1\le k\le 3 }L_k^{m,i}, \label{76-1}
\end{align}
where
\begin{align}
&{R}^{m,i}=\pl_t^m \pl_y^i  \lt(1+ \pl_y \oa\rt)^{-\ga}+\ga \lt(1+ \pl_y \oa\rt)^{-\ga-1}\pl_t^{m} \pl_y^{i+1} \oa,\label{4.22}\\
& E_1^{m,i}=\ga \int \chi^2  \sa^{\iota+i+1}(1+\pl_y \oa)^{-\ga-1}\lt|\pl_t^{m} \pl_y^{i+1} \oa \rt|^2 dy,\notag\\
& E_2^{m,i}=-2\int \chi^2  \sa^{\iota+i+1} R^{m,i}\pl_t^{m} \pl_y^{i+1} \oa  dy,\notag\\
&L_1^{m,i}=\frac{\ga}{2} \int \chi^2  \sa^{\iota+i+1} \lt( \pl_t (1+\pl_y \oa)^{-\ga-1} \rt) \lt|\pl_t^{m} \pl_y^{i+1} \oa \rt|^2 dy,\notag\\
&L_2^{m,i}=-\int \chi^2  \sa^{\iota+i+1} (\pl_t R^{m,i})\pl_t^{m} \pl_y^{i+1} \oa  dy,\notag\\
&L_3^{m,i}=2\int \chi \chi'  (\pl_t^{m+1}\pl_y^i \oa)\sa^{\iota+i+1} \pl_t^m \pl_y^i \lt(1+ \pl_y \oa\rt)^{-\ga}dy.\notag
\end{align}
Similarly, we integrate the product  of  \ef{75-4} and $\chi^2 \sa^{\iota+i}\pl_t^{m}  \pl_y^i \oa$ over $I$ to obtain
\begin{align}
&\frac{d}{dt} E_3^{m,i}-\int \chi^2  \sa^{\iota+i}\lt|\pl_t^{m+1} \pl_y^i \oa \rt|^2 dy +E_1^{m,i} =2 \sum_{k=4,5} L_k^{m,i}, \label{76-2}
\end{align}
where
\begin{align*}
&E_3^{m,i}=\frac{1}{2} \int \chi^2  \sa^{\iota+i}\lt(\lt|\pl_t^{m} \pl_y^i \oa \rt|^2 +2 (\pl_t^{m+1} \pl_y^i \oa) \pl_t^{m} \pl_y^i \oa\rt)dy,\\
&L_4^{m,i}=\frac{1}{2}\int \chi^2  \sa^{\iota+i+1}  R^{m,i} \pl_t^{m} \pl_y^{i+1} \oa  dy,\\
&L_5^{m,i}=\int \chi \chi'  (\pl_t^{m}\pl_y^i \oa)\sa^{\iota+i+1} \pl_t^m \pl_y^i  \lt(1+ \pl_y \oa\rt)^{-\ga}dy.
\end{align*}
It follows from $2\times\ef{76-1}+\ef{76-2}$ that
\begin{align}
& \frac{d}{dt}\mathfrak{E}^{m,i}(t)
+\mathfrak{D}^{m,i}(t) = 2\sum_{1\le k\le 5} L_k^{m,i}, \label{76-3}
\end{align}
where
\begin{align*}
&\mathfrak{E}^{m,i}
=\int \chi^2  \sa^{\iota+i}\lt|\pl_t^{m+1} \pl_y^i \oa \rt|^2 dy+\sum_{1\le k\le 3} E_k^{m,i}, \\
& \mathfrak{D}^{m,i}=\int \chi^2  \sa^{\iota+i}\lt|\pl_t^{m+1} \pl_y^i \oa \rt|^2 dy
+E_1^{m,i} .
\end{align*}

We analyze the terms on the right-hand side of \ef{76-3}. When $m+i=1$, it is easy to see that
\begin{subequations}\label{5.5-1}\begin{align}
& 0=R^{m,i}=E_2^{m,i}=L_2^{m,i}=L_4^{m,i},\\
& |L_1^{m,i}| \les \ea_0 \int \chi^2  \sa^{\iota+i+1}\lt|\pl_t^{m } \pl_y^{i+1} \oa \rt|^2 dy,
\end{align}\end{subequations}
due to \ef{4.21-2}.
When $m+i\ge 2$, note that
\begin{align*}
&|R^{m,i}|\les \sum_{r\le m,\ k\le i, \
1\le r+k\le m+i-1} \mathcal{I}^{r,k} |\pl_t^{m-r}\pl_y^{i+1-k} \oa |,\\
&|\pl_t R^{m,i}| \les \sum_{r\le m,\ k\le i, \
1\le r+k} \mathcal{I}^{r,k} |\pl_t^{m+1-r}\pl_y^{i+1-k} \oa |,
\end{align*}
because of \ef{4.21-5}.
Then, it yields from \ef{4.21-2}, \ef{4.28} and \ef{5.5-5} that for $2\le m+i \le 4+[\iota]$,
\begin{align}
& \int    \sa^{\iota+i+1}|R^{m,i}|^2 dy
\notag \\
& \les   \ea_0^2 \sum_{\substack{r\le m,\ k\le i, \\
1\le r+k\le m+i-2}}  \int \sa^{\iota+i+2-r-2k}|\pl_t^{m-r}\pl_y^{i+1-k} \oa |^2 dy\notag\\
& + \ea_0^2  \sum_{\substack{r\le m,\ k\le i, \\
 r+k = m+i-1}} \int  \sa^{\iota+k+1}  |\mathcal{I}^{r,k}|^2 dy
\les \ea_0^2  \sum_{1\le r+k\le m+i-1} \mathcal{D}^{r,k} ,\label{5.5-3}
\end{align}
and
\begin{align}
& \int  \sa^{\iota+i+1}|\pl_t R^{m,i}|^2 dy
 \notag
 \\
& \les   \ea_0^2 \sum_{\substack{r\le m,\ k\le i, \\
1\le r+k\le m+i-2}}  \int \sa^{\iota+i+2-r-2k}|\pl_t^{m+1-r}\pl_y^{i+1-k} \oa |^2 dy  \notag\\
& + \ea_0^2 \sum_{{r\le m,\ k\le i, \
r+k=m+i-1 }} \int  \sa^{\iota+k}  |\mathcal{I}^{r,k}|^2 dy\notag \\
&+ \ea_0^2  \int  \sa^{\iota+i+1}  |\mathcal{I}^{m,i}|^2 dy
\les \ea_0^2  \sum_{1\le r+k\le m+i} \mathcal{D}^{r,k}. \label{5.5-2}
\end{align}
Indeed,
the following estimate has been used to derive \ef{5.5-2}, which is due to \ef{hard}.
\begin{align*}
&\lt\| \sa^{\frac{\iota+i+2-r-2k}{2}} \pl_t^{m+1-r} \pl_y^{i+1-k}\oa\rt\|_{L^2}^2\\
  \les
  &
  \sum_{0\le l\le r+k-1} \lt\| \sa^{\frac{\iota+i+2-r-2k}{2}+r+k-1} \pl_t^{m+1-r} \pl_y^{i+1-k+l}\oa\rt\|_{L^2}^2 \\
  \le & \sum_{0\le l\le r+k-1} \mathcal{D}^{m+1-r,i-k+l}
\end{align*}
for $r\le m$, $ k\le i$,
$1\le r+k\le m+i-2$. So, we derive from the Cauchy inequality, and \ef{5.5-1}-\ef{5.5-2} that
for $1\le m+i \le 4+[\iota]$,
\begin{subequations}\label{5.6}\begin{align}
&|E_2^{m,i}|\le \ea_0 \int \chi^2  \sa^{\iota+i+1}\lt|\pl_t^{m} \pl_y^{i+1} \oa \rt|^2 dy \notag\\
&\qquad + \ea_0\sum_{1\le r+k\le m+i-1} \mathcal{D}^{r,k},\label{5.6-1}\\
&\sum_{k=1,2,4}|L_k^{m,i}|\les \ea_0\sum_{1\le r+k\le m+i} \mathcal{D}^{r,k}. \label{5.6-2}
\end{align}\end{subequations}

In view of  \ef{4.21-4}, we see that
for $1\le m+i\le [\iota]+4$,
\begin{align*}
& \int_{\hbar /4}^{\hbar/2} |\mathcal{I}^{m,i}|^2  dy
 \les \sum_ {\substack{ r\le m,   \  k \le i
\\  r +k \le m +i-1
  }}\int_{\hbar/4}^{\hbar/2}   \lt|\pl_t^{m-r} \pl_y^{i-k+1}\oa\rt|^2 dy \notag\\
&  \le \sum_ {\substack{ r\le m,   \  k \le i
\\  r +k \le m +i-1
  }}\int  \lt|\za_1 \pl_t^{m-r} \pl_y^{i-k+1}\oa\rt|^2 dy
  \le \sum_ {\substack{ r\le m,   \  k \le i
\\  r +k \le m +i-1
  }}\mathcal{D}_1^{m-r,i-k},
\end{align*}
which, together with the Cauchy inequality, \ef{4.21-5} and \ef{4.13-1}, implies that for $m+i\le [\iota]+4$ and $i\ge 1$,
\begin{align}
 \sum_{k=3,5}|L_k^{m,i}| & \les \int_{\hbar/4}^{\hbar/2} \lt(\lt|\pl_t^{m+1} \pl_y^{i} \oa \rt| +\lt|\pl_t^{m} \pl_y^{i} \oa \rt|\rt)\mathcal{I}^{m,i} dy\notag\\
&\les \sum_ {r\le m, \ k\le i
  }\mathcal{D}_1^{r,k} \les \sum_ {r\le m+i
  }\mathcal{D}^{r,0}  . \label{5.6-3}
\end{align}
As a conclusion of \ef{76-3}, \ef{5.6-2}, \ef{5.6-3} and
$L_3^{m,0}=L_5^{m,0}=0$, we obtain that for $1\le m+i\le [\iota]+4$,
\begin{align}\label{5.6-4}
& \frac{d}{dt}\mathfrak{E}^{m,i}(t)
+\mathfrak{D}^{m,i}(t)\les \mathcal{M}^{m,i},
\end{align}
where
\begin{align*}
& \mathcal{M}^{m,0}=\ea_0\sum_{1\le r+k\le m} \mathcal{D}^{r,k},\\
& \mathcal{M}^{m,i}=\ea_0\sum_{1\le r+k\le m+i} \mathcal{D}^{r,k}
 +  \sum_{r\le m+i} \mathcal{D}^{r,0}, \ \ i\ge 1.
\end{align*}

{\em Step 3}. In this step, we prove \ef{5.10-3}.
Clearly, it holds that
\begin{subequations}\label{5.7}\begin{align}
&\mathcal{D}^{m,0}\les \mathfrak{D}^{m,0}, \ \
\mathcal{D}_2^{m,i}\les  \mathfrak{D}^{m,i} \  \ {\rm for} \ \  i\ge 1, \label{5.7-a}\\
& \mathcal{D}^{m,i} \les \sum_{r\le m+i} \mathfrak{D}^{r,0}+ \mathfrak{D}^{m,i} \  \ {\rm for} \ \  i\ge 1, \label{5.7-b}
\end{align}\end{subequations}
where \ef{5.7-b} follows from \ef{4.13-1} and \ef{5.7-a}.
In view of \ef{5.6-4} and
  \ef{5.7-a}, we see that there exists a positive constant $K$  only  depending on $\ga$ and $M$ such that
\begin{align*}
& \sum_{1\le i, \ m+i\le [\iota]+4} \lt\{ \frac{d}{dt}
\mathfrak{E}^{m,i}(t)
+\mathfrak{D}^{m,i}(t) \rt\} \notag\\
& \le
K \ea_0\sum_{1\le r+k\le [\iota]+4} \mathcal{D}^{r,k} + K
\sum_{r\le [\iota]+4}  \mathfrak{D}^{r,0},
\end{align*}
which, together with \ef{75-3}, \ef{5.6-4}, \ef{5.7} and the smallness of $\ea_0$, implies that
\begin{align}\label{5.7-1}
& \frac{d}{dt}\mathfrak{E}(t)+ \frac{1}{2}
\mathfrak{D}(t) \le  0 ,
\end{align}
where
\begin{align*}
&\mathfrak{E}=\sum_{m\le [\iota]+4}(K+1)\mathfrak{E}^{m,0}
+\sum_{1\le i, \ m+i\le [\iota]+4}\mathfrak{E}^{m,i},\\
& \mathfrak{D}=\sum_{m+i\le [\iota]+4}\mathfrak{D}^{m,i}.
\end{align*}

It follows from the Cauchy inequality,   \ef{5.6-1} and the smallness of $\ea_0$ that for $m+i\ge 1$,
\begin{align*}
&\mathfrak{E}^{m,i} \ge 6^{-1}\lt( \mathfrak{D}^{m,i} + \lt\|\chi \sa^{\frac{\iota+i}{2}}\pl_t^{m} \pl_y^i \oa \rt\|_{L^2}^2 \rt)
-C\ea_0\sum_{1\le r+k\le m+i-1} \mathcal{D}^{r,k},\\
&\mathfrak{E}^{m,i} \le 2 \mathfrak{D}^{m,i} + \lt\|\chi \sa^{\frac{\iota+i}{2}}\pl_t^{m} \pl_y^i \oa \rt\|_{L^2}^2
+C\ea_0\sum_{1\le r+k\le m+i-1} \mathcal{D}^{r,k},
\end{align*}
for a certain positive constant $C$ only depending on $\ga$ and $M$. This,
 together with
\ef{5.8-1}, \ef{5.7} and the smallness of $\ea_0$, implies that
\begin{align}\label{5.10}
\mathfrak{E} \thicksim \mathfrak{D}
+ \lt\|\sa^{\frac{\iota}{2}}\oa \rt\|_{L^2}^2.
\end{align}
Due to the boundary condition $\oa (t,0)=0$, and the Sobolev embedding  $H^{\iota+[\iota]+4,[\iota]+3}(I)\hookrightarrow H^{(2+[\iota]-\iota)/2}(I)$ $ \hookrightarrow L^\iy(I)$, we have
\begin{align*}
\lt|\oa(t,y)\rt|=&\lt|\int_0^y \pl_z\oa(t,z)dz\rt|\les \lt\|\pl_y\oa(t,\cdot)\rt\|_{L^\iy}
\les    \sum_{i\le [\iota]+3} \mathcal{D}^{0,i}
\les  \mathfrak{D},
\end{align*}
where \ef{5.8-1} and \ef{5.7} have been used to derive the last inequality. Substitute this into \ef{5.10} to obtain
\be
\label{5.10-1}
\mathfrak{D}
+ \lt\|\sa^{\frac{\iota}{2}}\oa \rt\|_{L^2}^2\les \mathfrak{E}\les \mathfrak{D}.
\ee

Let $\da\in (0,1/4)$ be a small constant to be determined later, then we multiply \ef{5.7-1} by $e^{\da t}$, and use   \ef{5.10-1} to get
\begin{align*}
& \frac{d}{dt}\lt(e^{\da t}\mathfrak{E}(t)\rt)+ \frac{1}{2} e^{\da t}
\mathfrak{D}(t) \le  \da e^{\da t}\mathfrak{E}(t)\les \da e^{\da t}\mathfrak{D}(t)  .
\end{align*}
Hence, there exists a constant $\da\le 1/4$ only depending on $\ga$ and $M$ such that
\begin{align*}
& \frac{d}{dt}\lt(e^{\da t}\mathfrak{E}(t)\rt)+ \frac{1}{4} e^{\da t}
\mathfrak{D}(t) \le 0.
\end{align*}
Integrate the above equation over $[0,t]$ to get
\be\label{5.10-2}
e^{\da t}\mathfrak{E}(t) + \frac{1}{4} \int_0^t  e^{\da s}
\mathfrak{D}( s ) ds \le \mathfrak{E}(0).
\ee
It follows from \ef{5.8-1}, \ef{5.7} and \ef{5.10-1} that
\begin{align*}
\mathcal{E}= & \sum_{m+i\le [\iota]+4} \lt(\mathcal{D}^{m,i}+ \lt\| \sa^{\frac{\iota+i}{2}} \pl_t^{m}   \pl_y^i \oa \rt\|^2_{L^2} \rt) \\
\le &  3\sum_{m+i\le [\iota]+4} \mathcal{D}^{m,i} +\lt\| \sa^{\frac{\iota}{2}} \oa \rt\|^2_{L^2}
\les \mathfrak{D} \les \mathfrak{E},
\end{align*}
which, together with \ef{5.10-2}, proves \ef{5.10-3}. \hfill $\Box$

\section{The two- and three-dimensional motions }
We set $\Oa$ as the reference domain, and
define $x$ as the Lagrangian flow of the velocity $u$ by
\be\label{lag}
\pl_t x(t,y)= u(t, x(t,y))  \  {\rm for}  \  t>0,  \ {\rm and}  \  x(0,y)=x_0(y)  \ {\rm for}  \ y\in \Omega.
\ee
We define $\varrho(t,y)$ as the Lagrangian density, $v(t,y)$ the Lagrangian velocity, $A(t,y)$  the inverse of the Jacobian matrix, and $J(t,y)$ the Jacobian determinant  by
\begin{align*}
&\varrho(t,y)=\rho(t, x(t,y)),  \ \ v(t,y)=  u(t, x(t,y)), \\
& {A}(t,y)=\lt(\frac{\pl x}{\pl y}\rt)^{-1}, \ \
{J}(t,y)={\rm det} \lt(\frac{\pl x}{\pl y}\rt).
\end{align*}
Thus, system \ef{eq-ed} can be written in Lagrangian coordinates as
\begin{subequations}\label{2.1n}\begin{align}
& \pl_t  \varrho  +   \varrho {A}_i^k \pl_k v^i = 0 &  {\rm in}&  \ (0,T]\times \Omega, \label{2.1na}\\
 &   \varrho \pl_{t} v_i   +  {A}_i^k \pl_k ( { \varrho^\ga} )  = -\varrho v_i -  \varrho  g \da_{in}  & {\rm in}&  \ (0,T]\times \Omega,\label{2.1nb}\\
 & \varrho =0  &   {\rm on}& \ (0,T]\times  \{y_n=\hbar\}, \label{2.1nc}\\
 &  v_n=0    &    {\rm on}& \  (0,T]\times \{y_n=0\}, \label{2.1nd}\\
&( \varrho,   x, v)=(\rho_0(x_0), x_0, u_0(x_0)) & {\rm on} &  \ \{t=0\} \times  \Omega, \label{2.1nf}
 \end{align} \end{subequations}
where $v^i=v_i$, $\pl_k=\frac{\pl}{\pl y_k}$.
It follows from  \ef{2.1na} and
$\pl_t {J} = {J} {A}_i^k \pl_k v^i   $ that
$$ \varrho(t,y){J}(t,y)= \varrho (0,y){J}(0,y)=\rho_0\lt(x_0(y)\rt) {\rm det} \lt(\frac{\pl x_0(y)}{\pl y}\rt) .$$
We choose $x_0(y)$ such that
$\rho_0\lt(x_0(y)\rt) {\rm det} \lt(\frac{\pl x_0(y)}{\pl y}\rt)=  \bar\rho (y)$.
The existence of such an $x_0$ follows from the Dacorogna-Moser theorem (cf. \cite{DM}) and \ef{initial density}. So, the Lagrangian density can be expressed as
\be\label{1.6'}
\varrho(t,y)=\bar\rho(y) {J}^{-1}(t,y), \ \ {\rm where} \ \
\bar\rho(y) =\lt(\nu(\hbar-y_n)\rt)^{{1}/({\ga-1})},
   \ee
and problem \ef{2.1n}  reduces to
\begin{subequations}\label{system}\begin{align}
&   \bar\rho \pl_{t} v_i   +  J{A}_i^k \pl_k \lt( \bar\rho^\ga J^{-\ga} \rt)  = -\bar\rho v_i -  \bar\rho  g \da_{in}  & {\rm in}&  \ (0,T]\times \Omega,\label{system-a}\\
 &  v_n=0    &    {\rm on}& \  (0,T]\times \{y_n=0\}, \label{system-c}\\
&(    x, v)=( x_0, u_0(x_0)) & {\rm on} &  \ \{t=0\} \times  \Omega. \label{system-d}
\end{align}\end{subequations}
In the setting, the moving vacuum boundary for problem \ef{eq-ed} is given by
\be
\label{mvb-3}
\Gamma(t,x_*(t,y))=x_n(t,y_*,\hbar).
\ee

\subsection{Notation and main results}

We let $\pl_k=\frac{\pl}{\pl y_k}$,
$\bar\pl^{\al_*}=\pl_1^{\al_1}\pl_2^{\al_2}\cdots\pl_{n-1}^{\al_{n-1}}$ for multi-index $\al_*=(\al_1,\al_2,\cdots, \al_{n-1})$,
$\pl^\al=\bar\pl^{\al_*}\pl_n^{\al_n}$ for multi-index $\al=(\al_*, \al_n)$, $\bar\pl^j=\lt\{\bar\pl^{\al_*}: \  |\al_*|=j\rt\}$ and
$\pl^j=\lt\{\pl^{\al}: \  |\al|=j\rt\}$ for nonnegative integer $j$.
The divergence of a vector filed $F$ is ${\rm div} F=\da^{k}_i \pl_k F^i$,  the curl of a vector filed $F$ for $n=2$ is ${\rm curl} F=\pl_1 F_2- \pl_2 F_1$,
and the $i$-th component of the curl of a vector filed $F$ for $n=3$ is
$   [ {\rm curl} F ]_i = \epsilon^{ijk} \pl_j F_k$,
where  $\epsilon^{ijk}$ is the standard permutation symbol given by
\begin{align*}
\epsilon^{ijk}=
\begin{cases}
&1,    \ \  \ \  \textrm{even permutation of} \ \ \{1,2,3\},   \\
&-1,   \ \  \  \textrm{odd permutation of} \ \ \{1,2,3\}, \\
&0,    \ \ \  \  \textrm{otherwise}.
\end{cases}
\end{align*}
Along the flow map $x$, the $i$-th component of the gradient  of a function $f$ is
$ \lt[\nabla_x f \rt]_i = A^k_i \pl_k f$,
the divergence of a vector filed $F$ is ${\rm div}_x F = A^k_i \pl_k F^i$, the curl of a vector filed $F$ for $n=2$ is ${\rm curl}_x F=\lt[\nabla_x F_2 \rt]_1-\lt[\nabla_x F_1 \rt]_2=A^r_1\pl_r F_2-  A^r_2 \pl_r F_1$,
and the $i$-th component of the curl of a vector filed $F$ for $n=3$ is
 $
[ {\rm curl}_x F ]_{i}=\epsilon^{ijk} \lt[\nabla_x F_k \rt]_j   =\epsilon^{ijk} A^r_j \pl_r F_k$.

We define the perturbation $\oa$ by
$\oa(t,y) = x(t,y)- y $,
and set $\sigma(y) = \bar\rho^{\ga-1}(y)=\nu(\hbar-y_n)$ and  $\iota= ({\ga-1})^{-1}$.
We introduce that for nonnegative integers $m,j,i$,
\begin{align*}
\mathcal{E}^{m,j,i}(t)= & \lt\| \sa^{\frac{\iota+i}{2}}  \pl_t^{m+1}  \bar\pl^j \pl_n^i \oa \rt\|_{L^2}^2+ \lt\| \sa^{\frac{\iota+i}{2}} \pl_t^{m}  \bar\pl^j \pl_n^i \oa \rt\|^2_{L^2} \notag  \\
& +   \lt\| \sa^{\frac{\iota+i+1}{2}}
  \pl_t^m \pl \bar\pl^j\pl_n^i \oa \rt\|^2_{L^2} ,
\end{align*}
and define the higher-order weighted Sobolev norm
$\mathcal{E}(t)$ by
\begin{align}\label{1123}
\mathcal{E}(t) =
\sum_{0\le m+j+i\le [\iota]+n+4} \mathcal{E}^{m,j,i}(t)
=\mathcal{E}_I(t)+\mathcal{E}_{II}(t),
\end{align}
where
$
\mathcal{E}_{I}(t)=
\sum_{ m+j+i\le [\iota]+n+3} \mathcal{E}^{m+1,j,i}(t)$ and $
\mathcal{E}_{II}(t)=
\sum_{ j+i\le [\iota]+n+4} \mathcal{E}^{0,j,i}(t).
$
In addition to \ef{1123}, we need the following Sobolev norm for the curl:
$$
 \mathcal{V}_{a}(t)=
 \sum_{j+i=[\iota]+n+4}\mathcal{V}^{0,j,i}(t), \ \ {\rm where} \ \ \mathcal{V}^{m,j,i}= \lt\|   \sa^{\frac{\iota+i+1}{2}}
 \pl_t^m  \bar\pl^j \pl_n^i {\rm curl}_x v \rt\|_{L^2}^2. $$

Now, it is ready to state the main result.

\begin{thm}\label{5.19-1}
Let $n=2,3$.
There exists a positive constant $\bar \ea>0$ such that if $\mathcal{E}(0)+\mathcal{V}_{a}(0) \le \bar\ea$, then problem \ef{system} admits a global smooth solution in $[0,\iy)\times \Oa$ satisfying
\begin{align}\label{21Ju18}
e^{\da t}\mathcal{E}_{I}(t)+ \mathcal{E}_{II}(t)
+e^{\da t} \mathcal{V}_{a}(t) \le C (\mathcal{E}(0)+\mathcal{V}_{a}(0)) \ \ {\rm for} \ t \ge 0,
\end{align}
where $C$ and $\da$ are positive constants which only depend on the adiabatic exponent $\ga$ and the initial total mass $M$, but do not depend on the time $t$.
\end{thm}

As a corollary of Theorem \ref{5.19-1}, we have the following theorem for solutions to the original vacuum free boundary problem \ef{eq-ed} concerning the convergence of  the velocity $u$ to that of the stationary solution, and the deviation  of the vacuum boundary $\Gamma
$ and density $\rho$ from those of the stationary solution.

\begin{thm}\label{5.19-2}
Let $n=2,3$.
There exists a positive constant $\bar \ea>0$ such that if $\mathcal{E}(0)+\mathcal{V}_{a}(0) \le \bar\ea$, then problem \ef{eq-ed} admits a global  smooth solution $(\rho,u,\Oa(t))$ for $t\in [0,\iy)$ satisfying
\begin{subequations}\label{21.6.18}
\begin{align}
&\lt|u(t,x(t,y))\rt|\le C \sqrt{e^{-\da t }(\mathcal{E}+\mathcal{V}_{a})(0)}, \\
& \lt|\rho(t,x(t,y))-\bar\rho (y)\rt|\le C (\hbar-y_n)^{{1}/({\ga-1})}\sqrt{(\mathcal{E}+\mathcal{V}_{a})(0)} , \\
& \lt|\Gamma(t,x_*(t,y))- \hbar \rt|\le C \sqrt{(\mathcal{E}+\mathcal{V}_{a})(0)},
 \end{align}\end{subequations}
 for all $y\in \Oa$ and $t\ge 0$. Here
 $C$ and $\da$ are positive constants which only depend on the adiabatic exponent $\ga$ and the initial total mass $M$, but do not depend on the time $t$.
\end{thm}

\subsection{Proof of Theorems \ref{5.19-1} and \ref{5.19-2}}
We first prove
Theorem \ref{5.19-1}. The proof of the global existence of smooth solutions is based on the local existence theory (cf. \cite{10',16'}), together with the following a priori estimates stated in Proposition \ref{5.19-3} whose proof will be given in the next two subsections.

\begin{prop}\label{5.19-3}Let $x(t,y)=y+\oa(t,y)$ be a solution to problem \ef{system} in the time interval $[0,T]$ satisfying  the a priori assumption
\begin{align}\label{1120}
\mathcal{E}(t)\le \ea_0^2, \ \  t\in [0,T],
\end{align}
for some suitably small fixed positive number $\ea_0$ independent of $t$. Then there exist positive constants $C$ and $\da$ independent of $t$, which only depend on $\ga$ and $M$, such that
for $t\in[0,T]$,
\begin{align*}
\mathcal{E}_I(t) \le C e^{-\da t}  \mathcal{E}_I(0), \ \  \mathcal{E}_{II}(t) +e^{\da t}\mathcal{V}_{a}(t)  \le C (\mathcal{E}+\mathcal{V}_{a})(0).
\end{align*}
\end{prop}

To show Theorem \ref{5.19-2}, we recall the following estimate obtained in Lemma 3.3 of \cite{HZeng}.
\begin{align}
&\sum_{
  m+j+2i\le 4} \lt\|\pl_t^{m}  \bar\pl^j \pl_n^i  \oa \rt\|_{L^\iy}^2  + \sum_{\substack{
  m+j+2i= 5 }}  \lt\| \pl_t^{m}  \bar\pl^j \pl_n^i \oa \rt\|_{H^{\frac{n+[\iota]-\iota}{2}}}^2     \notag\\
&     + \sum_{\substack{
6 \le m+j+2i \\
m+j+i\le  [\iota]+ n +3 }}  \lt\|\sa^{\frac{ m+j+2i-4}{2}}\pl_t^{m}  \bar\pl^j \pl_n^i \oa \rt\|_{L^\iy}^2 \les \mathcal{E}(t), \label{5.19-4}
\end{align}
provided that $\mathcal{E}(t)$ is finite. The proof of \ef{5.19-4} for $n=2$ is the same as that for $n=3$ shown in  \cite{HZeng}, which is based on \ef{wsv} and \ef{hard}. Indeed, it holds that for $ m+j+2i= 5$,
\begin{align}
 &\lt\| \pl_t^{m}  \bar\pl^j \pl_n^i \oa \rt\|_{H^{(n+[\iota]-\iota)/{2}}}
 \les  \lt\| \pl_t^{m}  \bar\pl^j \pl_n^i \oa \rt\|_{H^{\iota+i+b,b}}
 \les \mathcal{E}(t),\notag\\
 &\lt\|\sa \pl_t^{m}  \bar\pl^j \pl_n^i \oa \rt\|_{H^{(n+2+[\iota]-\iota)/{2}}}
 \les  \lt\| \sa \pl_t^{m}  \bar\pl^j \pl_n^i \oa \rt\|_{H^{\iota+i+b-2,b}}
 \les \mathcal{E}(t), \label{6.3.1}
\end{align}
where $b=[\iota]+n+5-(m+j+i)$.
It follows from \ef{lag}, \ef{1.6'} and \ef{mvb-3} that
\begin{align*}
& \rho(t,x(t,y))-\bar\rho (y)= \bar\rho(y) J^{-1}(t,y)\lt(1-J(t,y)\rt)   ,  \\
&u(t,x(t,y))=\pl_t \oa(t,y) , \ \ \Gamma(t,x_*(t,y))- \hbar = \oa_n(t,y_*,\hbar) .
 \end{align*}
This, together with \ef{21Ju18}, \ef{Jaco}, and \ef{5.19-4} which implies $\|\oa\|_{L^\iy}+\|\pl\oa\|_{L^\iy}\les \sqrt{\mathcal{E}_{II}}$ and $\|\pl_t\oa\|_{L^\iy}\les \sqrt{\mathcal{E}_{I}}$,  proves \ef{21.6.18}.

\subsection{A priori estimates for $\mathcal{E}^{m,j,i}$ with $m\ge 1$ }

We introduce the following Sobolev norms for the  dissipation:
\begin{align*}
&\mathcal{D}^{m,j,i}(t)=\lt\| \sa^{\frac{\iota+i}{2}}  \pl_t^{m+1}   \bar\pl^j \pl_n^i v \rt\|_{L^2}^2 +\lt\|   \sa^{\frac{\iota+i+1}{2}}
 \pl_t^m\pl \bar\pl^j \pl_n^i v \rt\|_{L^2}^2,
\end{align*}
where $m,j,i$ are nonnegative integers.
To specify  the behavior of solutions near the bottom and the top,  we divide $\mathcal{D}^{m,j,i}$
into two parts $\mathcal{D}_1^{m,j,i}$ and $\mathcal{D}_2^{m,j,i}$ as follows:
\begin{align*}
 &  \mathcal{D}_1^{m,j,i}(t)=\lt\| \za_1 \pl_t^{m+1}   \bar\pl^j \pl_n^i v \rt\|_{L^2}^2 +\lt\|\za_1
 \pl_t^m\pl \bar\pl^j \pl_n^i v \rt\|_{L^2}^2,\\
 & \mathcal{D}_2^{m,j,i}(t)=\lt\| \zeta_2\sa^{\frac{\iota+i}{2}}  \pl_t^{m+1}   \bar\pl^j \pl_n^i v \rt\|_{L^2}^2 +\lt\|  \zeta_2 \sa^{\frac{\iota+i+1}{2}}
 \pl_t^m\pl \bar\pl^j \pl_n^i v \rt\|_{L^2}^2,
\end{align*}
where $\zeta_1=\zeta_1(y_n)$ and $\zeta_2=\zeta_2(y_n)$ are  smooth cut-off functions satisfying \ef{4.13}.
We will use elliptic estimates to bound $\mathcal{D}_1^{m,j,i}$ for $i\ge 1$ in Section \ref{s3.3.3}, and energy estimates to bound $\mathcal{D}_2^{m,j,0}$, and
$\mathcal{D}_2^{m,j,i}$ for $i\ge 1$ in Section \ref{s3.3.4}.

\subsubsection{Preliminaries}
It follows from \ef{1120}-\ef{6.3.1} and $H^{(n+1)/2}\hookrightarrow L^\iy$ that for $t\in [0, T]$,
\begin{align}
&\sum_{
  m+j+2i\le 4} \lt\|\pl_t^{m}  \bar\pl^j \pl_n^i  \oa \rt\|_{L^\iy}^2  + \sum_{\substack{
  m+j+2i= 5 }} \lt( \lt\| \pl_t^{m}  \bar\pl^j \pl_n^i \oa \rt\|_{H^{\frac{n+[\iota]-\iota}{2}}}^2
  +\lt\| \sa \pl_t^{m}  \bar\pl^j \pl_n^i \oa \rt\|_{L^\iy}^2  \rt)    \notag\\
&     + \sum_{\substack{
6 \le m+j+2i \\
m+j+i\le  [\iota]+ n +3 }}  \lt\|\sa^{\frac{ m+j+2i-4}{2}}\pl_t^{m}  \bar\pl^j \pl_n^i \oa \rt\|_{L^\iy}^2 \les \mathcal{E}(t) \le \ea_0^2, \label{12-8}
\end{align}
which means that
\bee
|\pl\oa(t,y)| \les \ea_0   \ \ {\rm for} \ \  (t,y)\in [0,T]\times \Oa.
\eee
Since $JA$ is the adjugate matrix of $(\frac{\pl x}{\pl y})$ and $x(t,y)=y+\oa(t,y)$, then
\bee
J A = \lt(\frac{\pl x}{\pl y} \rt)^*  =   \lt(1 +{\rm div}\oa\rt)\mathbb{I} -  \lt(\frac{\pl\oa}{\pl y}\rt) + \da_{n3} B,
\eee
where $\mathbb{I}$ is the identity matrix, and $B$  is the adjugate matrix of $(\frac{\pl \oa}{\pl y})$ given by
\bee
B= \lt(\frac{\pl\oa}{\pl y} \rt)^* = \lt[\begin{split} \pl_2 \omega  \times \pl_3 \omega \\
\pl_3 \omega  \times  \pl_1\omega  \\
\pl_1 \omega  \times \pl_2 \omega \end{split}\rt]  .
\eee
This, together with the fact that $(\frac{\pl x}{\pl y})(\frac{\pl x}{\pl y})^*=J \mathbb{I} $, implies that
 \begin{align}\label{Jaco}
 J=1+{\rm div} \oa + 2^{-1}\lt(|{\rm div} \oa|^2 + |{\rm curl} \oa|^2- |\pl \oa|^2\rt) + 3^{-1}  \da_{n3} B^s_r \pl_s \oa^r .
\end{align}
So, we have  that  for $t\in [0,T]$,
\be\label{7.9}
\|J-1\|_{L^\iy} \les \|\pl \oa\|_{L^\iy}  \les \ea_0  \ {\rm and}  \   \|A-  \mathbb{I}\|_{L^\iy} \les \|\pl \oa\|_{L^\iy}  \les \ea_0,
\ee
which gives, with the aid of the smallness of $\ea_0$, that  for $(t,y)\in [0,T]\times \Oa$,
 \be\label{6.7-1a}
  2^{-1}\le J \le 2     \ \  {\rm and} \ \
\max_{1\le i\le n, \
 1\le j\le n}  \lt| A^i_j \rt|\le 2 .
    \ee
 Moreover,  we have
$
   |[\na_x f]_i -\pl_i f| =|( A^r_i - \da^r_i) \pl_r f| \les \ea_0 |\pl f|
$ for any function $f$,
which means
 \be\label{6.7-1c}
   2^{-1} |\pl f| \le |\na_x f| \le 2 |\pl f|.
\ee

The differentiation formulae for $J$ and $A$ are
\begin{align*}
&\pl_j J=J A^s_r \pl_{j }\pl_{s } \oa^{r},     \  \  \ \  \pl_t J=J A^s_r   \pl_s v^r,  \\
&\pl_j A^k_i = - A^k_r   A^s_i  \pl_{j }\pl_{s } \oa^r,  \ \  \ \
  \pl_t  A^k_i = - A^k_r A^s_i  \pl_s v^r,
\end{align*}
which, together with  the mathematical induction and \ef{6.7-1a}, implies that for  any polynomial function $\mathscr{P}$ and nonnegative integers $m,j,i$,
\begin{align*}
\lt|\pl_t^m  \bar\pl^j \pl_n^i \mathscr{P} (A )  \rt| + \lt|\pl_t^m  \bar\pl^j \pl_n^i \mathscr{P}(J) \rt|
\les  \mathcal{I}^{m,j,i},
\end{align*}
where $\mathcal{I}^{m,j,i}$ are defined inductively as follows:
\begin{subequations}\label{n5.30}\begin{align}
& \mathcal{I}^{0,0,0}=1, \\
& \mathcal{I}^{m,j,i} =
\sum_{\substack{0\le r\le m, \ 0\le l \le j, \ 0\le k \le i
\\ 0\le r+l+k \le m+j+i-1
  }} \mathcal{I}^{r,l,k} \lt|\pl_t^{m-r}\bar\pl^{j-l}\pl_n^{i-k}\pl\oa\rt|
. \label{5.30b}
 \end{align}\end{subequations}

\begin{lem}
It holds that for any nonnegative integers $m,j,i$, and $t\in [0,T]$,
\begin{subequations}\begin{align}
&\sum_{\substack{ r\le m, \  l \le j, \  k \le i
\\ 1\le r+l+k \le m+j+i-1
  }}\lt\|\sa^{\frac{\iota+i+1}{2}}\mathcal{I}^{r,l,k}
\pl_t^{m-r}\bar\pl^{j-l}
\pl_n^{i-k}\pl\oa\rt\|_{L^2}^2
\notag\\
& \les \mathcal{E}(t)\sum_{r+l+k\le m+j+i-1}\mathcal{E}^{r,l,k}(t),\ \ \ \   2\le m+j+i \le [\iota]+n+5,\label{1.27-1}\\
&\sum_ {\substack{ r\le m, \  l \le j, \  k \le i
\\ 1 \le r+l+k
  }}\lt\|  \sa^{\frac{\iota+i+1}{2}} \mathcal{I}^{r,l,k} \pl_t^{m-r}\bar\pl^{j-l}\pl_n^{i-k}\pl v\rt\|_{L^2}^2 \notag\\
&  \les \mathcal{E}(t)\sum_{r+l+k\le m+j+i-1}\mathcal{D}^{r,l,k} (t), \  \ 1\le m+j+i \le [\iota]+n+4  ,\label{1.27-2} \\
& \sum_ {\substack{ (l,k)\in \mathfrak{S} \setminus \mathfrak{S}_1
  }}\lt\|  \sa^{\frac{\iota+i+1}{2}} \mathcal{I}^{1,l,k} \bar\pl^{j-l}\pl_n^{i-k}\pl v\rt\|_{L^2}^2\notag\\
&  \les \mathcal{E}(t)\sum_{l+k\le j+i-1}\mathcal{D}^{0,l,k} (t), \  \ 1\le j+i \le [\iota]+n+4, \label{5-30-1}
 \end{align}\end{subequations}
where $\mathfrak{S}=\{(l,k)\in \mathbb{Z}^2\big|0\le l\le j, \  0\le k\le i\}$ and  $\mathfrak{S}_1=\{(0,0), \ (0,1),  \ (j,i-1), \ (j,i) \}$.
\end{lem}

{\em Proof}. The proof of \ef{1.27-1}  for $n=3$ can be found in Lemma 4.7 of \cite{HZeng}, and that for $n=2$ can be shown in a similar way, so we omit the details. In fact, the derivation of \ef{1.27-1} is the same as that of \ef{1.27-2} which will be given in what follows.
To simplify the presentation, we set
\begin{align*}
&S=\lt\{(r,l,k)\in \mathbb{Z}^3 \ \big|\ 0\le r\le m, \  0\le l \le j, \  0\le k \le i
  \rt\},\\
 & S_1=\lt\{(r,l,k)\in S \ \big|\ 1\le  r+l+k \le m+j+i-1
  \rt\} ,\\
&   S_2=\lt\{(r,l,k)\in S \ \big|\ 1\le  r+l+k
  \rt\} .
\end{align*}

{\em Step 1}. In this step, we prove that
\begin{align}
&\sum_{
1\le   m+j+2i\le 2} \lt\|\mathcal{I}^{m,j,i} \rt\|_{L^\iy}^2
  +\sum_{
  m+j+2i= 3}\lt(  \lt\|\mathcal{I}^{m,j,i} \rt\|_{{L^{q^*}}}^2
  +\lt\|\sa\mathcal{I}^{m,j,i} \rt\|_{{L^{\iy}}}^2\rt) \notag\\
&   + \sum_{\substack{
4 \le m+j+2i \\
m+j+i\le  [\iota]+ n +2 }}  \lt\|\sa^{\frac{m+j+2i-2}{2}}
\mathcal{I}^{m,j,i}\rt\|_{L^\iy}^2
\les \mathcal{E}, \label{12.10}
\end{align}
where $q^*={2n}/({\iota-[\iota]})$ for
$\iota\neq [\iota]$, and $q^*\in (2, \iy)$ for $\iota = [\iota]$.

The proof is based on \ef{12-8}, the Sobolev embedding and the mathematical induction.
When $1\le m+j+2i\le 2$, it follows from \ef{n5.30} and \ef{12-8}   that
\begin{align*}
& \lt\|\mathcal{I}^{m,j,i} \rt\|_{L^\iy}
 \le  \lt\|\pl_t^{m}\bar\pl^{j}
 \pl_n^{i}\pl\oa\rt\|_{L^\iy}
 +\sum_ {(r,l,k)\in S_1}\lt\|\mathcal{I}^{r,l,k}\rt\|_{L^\iy}\\
& \times
\lt\|\pl_t^{m-r}\bar\pl^{j-l}\pl_n^{i-k}
\pl\oa\rt\|_{L^\iy}
\les \sqrt{\mathcal{E}}+ \ea_0 \sum_ {(r,l,k)\in S_1}\lt\|\mathcal{I}^{r,l,k}\rt\|_{L^\iy},
\end{align*}
which, together with the mathematical induction, implies that
\be\label{1.27-3}
\sum_{
1\le  m+j+2i\le 2} \lt\|\mathcal{I}^{m,j,i} \rt\|_{L^\iy}
  \les \sqrt{\mathcal{E}}.
\ee
When
$ m+j+2i= 3$, it yields from \ef{n5.30}, \ef{12-8},   $H^{(n+[\iota]-\iota)/{2}}\hookrightarrow L^{q^*}$, and \ef{1.27-3} that
\begin{align}
& \lt\|\mathcal{I}^{m,j,i} \rt\|_{L^{q^*}}
 \le \lt\|\pl_t^{m}\bar\pl^{j}
 \pl_n^{i}\pl\oa\rt\|_{L^{q^*}}
 +\sum_ {(r,l,k)\in S_1}\lt\|\mathcal{I}^{r,l,k}\rt\|_{L^\iy}\notag\\
& \ \ \times
\lt\|\pl_t^{m-r}\bar\pl^{j-l}\pl_n^{i-k}
\pl\oa\rt\|_{L^{q^*}} \les \lt\|\pl_t^{m}\bar\pl^{j}
 \pl_n^{i}\pl\oa\rt\|_{H^{(n+[\iota]-\iota)/{2}}}
\notag\\
 & \ \  +\sum_ {(r,l,k)\in S_1}\lt\|\mathcal{I}^{r,l,k}\rt\|_{L^\iy}
\lt\|\pl_t^{m-r}\bar\pl^{j-l}\pl_n^{i-k}
\pl\oa\rt\|_{L^\iy}
\les \sqrt{\mathcal{E}}, \notag
\\
& \lt\|\sa \mathcal{I}^{m,j,i} \rt\|_{L^\iy}
 \le \lt\|\sa\pl_t^{m}\bar\pl^{j}
 \pl_n^{i}\pl\oa\rt\|_{L^\iy}
 +\sum_ {(r,l,k)\in S_1}\lt\|\mathcal{I}^{r,l,k}\rt\|_{L^\iy}\notag\\
 & \ \ \times
\lt\|\pl_t^{m-r}\bar\pl^{j-l}\pl_n^{i-k}
\pl\oa\rt\|_{L^\iy}
\les \sqrt{\mathcal{E}}, \label{Feb2-1}
\end{align}
due to  $r+l+2k\le  2$ and
$m-r+j-l+2(i-k)\le 2$.
When
$ m+j+2i \ge 4$ and $m+j+i\le  [\iota]+ n +2$,  it
produces from \ef{n5.30} and \ef{12-8} that
\begin{align*}
& \lt\|\sa^{\frac{m+j+2i-2}{2}} \mathcal{I}^{m,j,i} \rt\|_{L^\iy}
 \le \lt\|\sa^{\frac{m+j+2i-2}{2}} \pl_t^{m}\bar\pl^{j}
 \pl_n^{i}\pl\oa\rt\|_{L^\iy}
 \\
& +\sum_ {(r,l,k)\in S_1}\lt\|\sa^{\frac{r+l+2k-1}{2}} \mathcal{I}^{r,l,k}\rt\|_{L^\iy}
\lt\|\sa^{\frac{m-r+j-l+2(i-k)-1}{2}} \pl_t^{m-r}\bar\pl^{j-l}\pl_n^{i-k}
\pl\oa\rt\|_{L^\iy}\\
& \les  \sqrt{\mathcal{E}}+ \ea_0\sum_ {(r,l,k)\in S_1}\lt\|\sa^{\frac{r+l+2k-1}{2}} \mathcal{I}^{r,l,k}\rt\|_{L^\iy},
\end{align*}
which, together with \ef{1.27-3}, \ef{Feb2-1} and the mathematical induction, gives the bound for the last term on the left-hand side of \ef{12.10}.

{\em Step 2}. In this step, we prove that
\begin{subequations}\label{Feb1}\begin{align}
&\sum_{\substack{5-n\le z\le 3, \ 3-z\le h\le \min\{2,\ 3-z/2\}\\
1\le m+j+i \le [\iota]+n+z+h-1 } }\lt\|\sa^{\frac{\iota+i+h}{2}} \mathcal{I}^{m,j,i}\rt\|_{L^{q_z}}^2 \les \mathcal{E}, \ {\rm if} \
\iota\neq[\iota], \label{Feb1-a}\\
&\sum_{\substack{4-n\le z\le 3 \\ \max\{0,\ 2-z\}\le h\le \min\{2,\ (n+3-z)/2\}\\
1\le m+j+i \le [\iota]+2+z+h} }\lt\|\sa^{\frac{\iota+i+h}{2}} \mathcal{I}^{m,j,i}\rt\|_{L^{q_z}}^2 \les \mathcal{E}, \ {\rm if} \
\iota=[\iota], \label{Feb1-b}\\
&\sum_{\substack{
1\le m+j+i \le [\iota]+n+3 } }\lt\|\sa^{\frac{\iota+i+2}{2}} \mathcal{I}^{m,j,i}\rt\|_{L^{q_2}}^2 \les \mathcal{E}, \ {\rm if}  \ n=2, \label{Feb1-c}\\
&\sum_{\substack{z=2,3,\ 3-z\le h \le 1\\
1\le m+j+i \le [\iota]+n+z+h} }\lt\|\sa^{\frac{\iota+i+h}{2}} \mathcal{I}^{m,j,i}\rt\|_{L^{q_z}}^2 \les \mathcal{E}, \ {\rm if} \
\iota=[\iota], \ n=3, \label{Feb1-d}\\
 &\sum_{h=0,1, \ 1\le m+j+i \le [\iota]+n+3+h }\lt\|\sa^{\frac{\iota+i+h}{2}} \mathcal{I}^{m,j,i}\rt\|_{L^{2}}^2
\les \mathcal{E},\label{Feb1-e}
\end{align}\end{subequations}
where $q_z=2n/(n+z-3+[\iota]-\iota)$. To prove \ef{Feb1-a}-\ef{Feb1-d}, we set
$$\mathcal{Q}^{m,j,i}_{r,l,k}
=\lt\|  \sa^{\frac{\iota+i+h}{2}} \mathcal{I}^{r,l,k} \pl_t^{m-r}\bar\pl^{j-l}\pl_n^{i-k}\pl \oa \rt\|_{L^{q_z}}^2, \ \  (r,l,k)\in S.$$

It follows from $dist(y, \pl \Oa)\les \sa(y) $, \ef{hard} and $a+\iota+i+h-2b>-1$ that for $(r,l,k)\in S$,
\begin{align*}
&\lt\|\sa^{\frac{\iota+i+h}{2}}
\pl_t^{m-r}\bar\pl^{j-l}\pl_n^{i-k}\pl\oa
\rt\|_{H^{a,b}}^2\\
&\les \sum_{ h_1 \le b} \sum_{ h_2 \le h_1}
\lt\|\sa^{\frac{a+\iota+i+h}{2}-h_2}
\pl_t^{m-r}\bar\pl^{j-l}\pl_n^{i-k }\pl^{1+h_1-h_2}\oa
\rt\|_{L^2}^2 \\
& \les \sum_{ h_1 \le b} \sum_{ h_2 \le h_1} \sum_{ h_3 \le h_2}
\lt\|\sa^{\frac{a+\iota+i+h}{2}}
\pl_t^{m-r}\bar\pl^{j-l}\pl_n^{i-k+h_3 }\pl^{1+h_1-h_2}\oa
\rt\|_{L^2}^2  \les \mathcal{E} ,
\end{align*}
where
$b=[\iota]+n +4-(m+j+i)+r+l+k$ and
$a =\max\{b-k+1, \ 2b-1-\iota-i+\hat\da\}-h$ for any $\hat\da>0$.  This, together with
\ef{wsv}, implies that for $(r,l,k)\in S$,
\be\label{21.6.9-1}
\lt\|\sa^{\frac{\iota+i+h}{2}}
\pl_t^{m-r}\bar\pl^{j-l}\pl_n^{i-k}\pl\oa
\rt\|_{H^{s_1}}   \les  \sqrt{\mathcal{E}},
\ee
where $$s_1=\min\{[\iota]+n+3+(r+l+2k)+h-(m+j+i), \ \iota+i+1+h-\hat\da\}/2$$
 for any $\hat\da>0$.
Similarly, we have that for $(r,l,k)\in S$ and $r+l+2k\ge 4$,
\be\label{21.6.9-2}
\lt\|\sa^{\frac{\iota+i+2+h-(r+l+2k)}{2}}
\pl_t^{m-r}\bar\pl^{j-l}\pl_n^{i-k}\pl\oa
\rt\|_{H^{s_2}}\les \sqrt{\mathcal{E}},
\ee
where
$$s_2=\min\{[\iota]+n+5+h-(m+j+i), \ \iota+3+i+h-(r+l+2k)-\hat\da\}/2$$ for any $\hat\da>0$, because of $H^{a,b}\hookrightarrow H^{s_2}$
with $b=[\iota]+n +4-(m+j+i)+r+l+k$ and
$a =\max\{b-1+r+l+k, \ 2b-3-\iota-i+(r+l+2k)+\hat\da\}-h$ for any $\hat\da>0$.

We first prove \ef{Feb1-a}. In this case, \ef{21.6.9-1} hold  with
$
s_1=\min\{4-z+r+l+2k, \ \iota+1+h-(1+[\iota]-\iota)/2\}/2.
$
When $r+l+2k\le 2$, we use \ef{12.10}, \ef{21.6.9-1}, and $H^{ s_1} \hookrightarrow H^{(3-z+\iota-[\iota])/2}\hookrightarrow L^{q_z}$ to get
$$
\mathcal{Q}^{m,j,i}_{r,l,k} \le  \lt\| \mathcal{I}^{r,l,k}
\rt\|_{L^\iy}^2\lt\|\sa^{\frac{\iota+i+h}{2}}
\pl_t^{m-r}\bar\pl^{j-l}\pl_n^{i-k}\pl\oa
\rt\|_{L^{q_z}}^2\les \mathcal{E}. $$
When $r+l+2k=3$, it follows from \ef{12.10}, \ef{21.6.9-1}, and $H^{ s_1} \hookrightarrow H^{(3-z+2\iota-2[\iota])/2}\hookrightarrow L^{q}$ that
$$
\mathcal{Q}^{m,j,i}_{r,l,k}
 \le \lt\| \mathcal{I}^{r,l,k}
\rt\|_{L^{q^*}}^2\lt\|\sa^{\frac{\iota+i+h}{2}}
\pl_t^{m-r}\bar\pl^{j-l}\pl_n^{i-k}\pl\oa
\rt\|_{L^q}^2 \les \mathcal{E}^2,
$$
where $q=2n/(n+z-3+2[\iota]-2\iota))$.
When $4\le r+l+2k\le m+j+2i-3$,  we use \ef{12.10}, \ef{21.6.9-2} with
$
s_2 = (6.5-n-z+\iota-[\iota] )/2,
$
and $ H^{ s_2} \hookrightarrow H^{(3-z+\iota-[\iota])/2}\hookrightarrow L^{q_z}$   to obtain
$$
\mathcal{Q}^{m,j,i}_{r,l,k}\le \lt\|\sa^{\frac{r+l+2k-2}{2}} \mathcal{I}^{r,l,k}
\rt\|_{L^\iy}^2\lt\|\sa^{\frac{\iota+i+2+h-(r+l+2k)}{2}}
\pl_t^{m-r}\bar\pl^{j-l}\pl_n^{i-k}\pl\oa
\rt\|_{L^{q_z}}^2
\les\mathcal{E}^2.$$
When $r+l+2k\ge \max\{4, \ m+j+2i-2\}$, it follows from \ef{12-8} that
$$
\mathcal{Q}^{m,j,i}_{r,l,k}\le \lt\|\sa^{\frac{\iota+i+h}{2}} \mathcal{I}^{r,l,k}
\rt\|_{L^{q_z}}^2\lt\|
\pl_t^{m-r}\bar\pl^{j-l}\pl_n^{i-k}\pl\oa
\rt\|_{L^\iy}^2
\les \ea_0^2 \lt\|\sa^{\frac{\iota+k+h}{2}} \mathcal{I}^{r,l,k}
\rt\|_{L^{q_z}}^2.$$
So, \ef{Feb1-a} is a consequence of  the mathematical induction and the estimates obtained above.

Similarly, we can prove \ef{Feb1-b}. Indeed, \ef{21.6.9-1} holds with
$s_1=\min\{n+1-z+r+l+2k, \ 1+h+\iota/2\}/2$ satisfying that $s_1\ge (3-z)/2$ for $r+l+2k\le 2$, and $s_1\ge (3-z)/2+\iota/4 $ for $r+l+2k=3$;
and  \ef{21.6.9-2} holds with $s_2 = ( 3.5-z ) /2$ satisfying $s_2\ge (3-z)/2$ for $4\le r+l+2k\le m+j+2i-3$.

The proof of \ef{Feb1-c} is the same as that of \ef{Feb1-a} except for the case of $r+l+2k=3$. In fact,
\ef{21.6.9-1} holds with
$s_1=1$  for $r+l+2k\le 2$, and \ef{21.6.9-2} holds with
$s_2=1$  for $4\le r+l+2k\le m+j+2i-3$.
When $r+l+2k=3$,  one has
$$
\mathcal{Q}^{m,j,i}_{r,l,k}
 \le \lt\|\sa^{\frac{1}{2}}\mathcal{I}^{r,l,k}
\rt\|_{L^{q_2}}^2\lt\|\sa^{\frac{\iota+i+1}{2}}
\pl_t^{m-r}\bar\pl^{j-l}\pl_n^{i-k}\pl\oa
\rt\|_{L^\iy}^2 \les \mathcal{E}^2,
$$
since $ \|\sa^{1/2}\pl_t^r\bar\pl^l \pl_n^k \pl \oa
\|_{H^{a,b}}\les \sqrt{\mathcal{E}}$ with $b=[\iota]+n+4-(r+l+k)$ and $a=\max\{\iota+k+b, \ 2b-1+[\iota]-\iota\}$ satisfying $b-a/2\ge (1+\iota-[\iota])/2$, and
$\|\sa^{({\iota+i+1})/{2}}
\pl_t^{m-r}\bar\pl^{j-l}\pl_n^{i-k}\pl\oa\|_{H^{s_1}}
\les \sqrt{\mathcal{E}}
$ with $s_1\ge 1+ \min\{\iota/4, 1\}$ by use of  \ef{21.6.9-1}.

Next, we prove \ef{Feb1-d}. Notice that \ef{21.6.9-1} holds with
$s_1=\min\{3-z+r+l+2k, \ 1+h+\iota/2\}/2$ satisfying that $s_1\ge (3-z)/2$ for $r+l+2k\le 2$, and $s_1\ge (3-z)/2+\iota/4 $ for $r+l+2k=3$; and \ef{21.6.9-2} holds with $s_2 = ( 3.5-z ) /2$ satisfying $s_2\ge (3-z)/2$ for $4\le r+l+2k\le m+j+2i-4$. When $r+l+2k= m+j+2i-3$, it follows from \ef{12.10}, \ef{Feb1-b}, and $ r+l+k = m+j+i-2$ or $m+j+i-3$ that
$$
\mathcal{Q}^{m,j,i}_{r,l,k}\le \lt\|\sa^{\frac{\iota+k+h}{2}} \mathcal{I}^{r,l,k}
\rt\|_{L^{q_{z-1}}}^2\lt\|
\pl_t^{m-r}\bar\pl^{j-l}\pl_n^{i-k}\pl\oa
\rt\|_{L^{q^*}}^2 \les \mathcal{E}^2 .
$$

When $[\iota]=\iota$, \ef{Feb1-e} follows from $q_3=2$, \ef{Feb1-b} and \ef{Feb1-d}. When $[\iota]\neq \iota$, we can use a similar way to the derivation of \ef{Feb1-a} to prove  \ef{Feb1-e}.
\ef{21.6.9-1} holds with
$s_1=\min\{r+l+2k, \ \iota+1/2\}/2$ satisfying $s_1\ge 0$ for $r+l+2k\le 2$, and $s_1\ge (\iota-[\iota])/2$ for $r+l+2k=3$; and \ef{21.6.9-2} holds with
$s_2\ge (\iota-[\iota])/4>0$  for $4\le r+l+2k\le m+j+2i-3$.  Indeed, when $r+l+2k=3$, it holds that $H^{ s_1} \hookrightarrow H^{(\iota-[\iota])/2}\hookrightarrow L^{q_3}$   and
\begin{align*}
&\lt\|  \sa^{\frac{\iota+i+h}{2}} \mathcal{I}^{r,l,k} \pl_t^{m-r}\bar\pl^{j-l}\pl_n^{i-k}\pl \oa \rt\|_{L^{2}}\\
&\le \lt\|  \mathcal{I}^{r,l,k} \rt\|_{L^{q^*}} \lt\|  \sa^{\frac{\iota+i+h}{2}} \pl_t^{m-r}\bar\pl^{j-l}\pl_n^{i-k}\pl \oa \rt\|_{L^{q_3}}\les \mathcal{E}^2.
\end{align*}

{\em Step 3}. Based on \ef{12.10} and \ef{Feb1},  we can prove \ef{1.27-2} in this step. We set
\begin{align*}
&\mathcal{P}^{m,j,i}_{r,l,k}=\lt\|  \sa^{\frac{\iota+i+1}{2}} \mathcal{I}^{r,l,k} \pl_t^{m-r}\bar\pl^{j-l}\pl_n^{i-k}\pl v\rt\|_{L^2}^2, \ \  (r,l,k)\in S_2,\\
& \mathfrak{T}=\mathcal{E}(t)\sum_{r+l+k\le m+j+i-1}\mathcal{D}^{r,l,k} (t).
\end{align*}
When $r+l+2k \le 2$, it yields from \ef{12.10} that
\begin{align}
&\mathcal{P}^{m,j,i}_{r,l,k}
\le \lt\| \mathcal{I}^{r,l,k} \rt\|_{L^\iy}^2  \lt\|  \sa^{\frac{\iota+i+1}{2}}  \pl_t^{m-r}\bar\pl^{j-l}\pl_n^{i-k}\pl v\rt\|_{L^2}^2 \les \mathfrak{T}.\label{1-27-1}
\end{align}
When $r+l+2k =3$, it produces from \ef{12.10},
$H^{1,1}\hookrightarrow H^{1/{2}}\hookrightarrow L^{2n/(n-1)}$, \ef{hard}
and  $r+l+k\ge 2$
that
\begin{align}
&\mathcal{P}^{m,j,i}_{r,l,k}
\le \lt\| \mathcal{I}^{r,l,k} \rt\|_{L^{2n}}^2  \lt\|  \sa^{\frac{\iota+i+1}{2}}  \pl_t^{m-r}\bar\pl^{j-l}\pl_n^{i-k}\pl v\rt\|_{L^{2n/(n-1)}}^2
\les \mathfrak{T}. \label{Feb8}
\end{align}

To deal with
the case of $r+l+2k \ge 4$, we set
\begin{align*}
&S_3=\lt\{(r,l,k)\in S_2 \ \big|\  4\le r+l+2k, \ \ r+l+k \le [\iota]+n+2
  \rt\},\\
&S_{31}=\lt\{(r,l,k)\in S_3 \ \big|\  r+l+2k <\iota+i+4
  \rt\},\\
&   S_{32}=S_3\setminus S_{31}=\lt\{(r,l,k)\in S_3 \ \big|\  r+l+2k \ge \iota+i+4
  \rt\}.
\end{align*}
When $(r,l,k)\in S_{31}$, we have
\begin{align}
&\mathcal{P}^{m,j,i}_{r,l,k}\le \lt\|\sa^{\frac{r+l+2k-2}{2}} \mathcal{I}^{r,l,k}\rt\|_{L^\iy}^2
\lt\|\sa^{\frac{\iota+i+3-r-l-2k}{2}} \pl_t^{m-r}\bar\pl^{j-l}\pl_n^{i-k}\pl v\rt\|_{L^2}^2\notag\\
& \les  \mathcal{E} \sum_{h_1 \le r+l+k-2} \lt\|\sa^{\frac{\iota+i+r+l-1}{2}} \pl_t^{m-r}\bar\pl^{j-l}\pl_n^{i-k+h_1}\pl v\rt\|_{L^2}^2\notag\\
& \les  \mathcal{E}\sum_{a+b+c\le m+j+i-2}\mathcal{D}^{a,b,c}, \label{June7}
\end{align}
where the second inequality follows from
\ef{12.10}, \ef{hard} and $\iota+i+3-r-l-2k>-1$. When $(r,l,k)\in S_{32}$, one gets
\be\label{6.5-1}
 \lt\| \pl_t^{m-r}\bar\pl^{j-l}\pl_n^{i-k}\pl v\rt\|_{L^{2n}}^2 \les
\sum_{r+l+k\le m+j+i-1}\mathcal{D}^{r,l,k}.
\ee
If $(r,l,k)\in S_{32}$,  then $r+l+k =[\iota]+n+1$ or $[\iota]+n+2$ for $n=3$, and $r+l+k =[\iota]+n+2$ for $n=2$.
Consider the case of $n=3$ and $r+l+k =[\iota]+n+1$, in which $[\iota]=\iota$ and $k=i$, then $H^{\iota+[\iota]+n+1,\ [\iota]+n} \hookrightarrow H^{(n-1)/{2}}\hookrightarrow L^{2n}$.
Consider the case of $n=3$ and $r+l+k =[\iota]+n+2$, in which $k=i$, or $[\iota]=\iota$ and $k=i-1$,  thus $H^{\iota+[\iota]+n+2,\ [\iota]+n+1} \hookrightarrow H^{(n+[\iota]-\iota)/{2}}\hookrightarrow L^{2n}$  or
 $H^{\iota+[\iota]+n+3,\ [\iota]+n+1} \hookrightarrow H^{(n-1)/{2}}\hookrightarrow L^{2n}$.
Consider the case of $n=2$ and $r+l+k =[\iota]+n+2$, in which $[\iota]=\iota$ and $k=i$, then $H^{\iota+[\iota]+n+2,\ [\iota]+n+1} \hookrightarrow H^{n/{2}}\hookrightarrow L^{2n}$. So, it follows from \ef{6.5-1} and \ef{Feb1} that for $(r,l,k)\in S_{32}$,
\begin{align}
\mathcal{P}^{m,j,i}_{r,l,k} \le \lt\|  \sa^{\frac{\iota+i+1}{2}} \mathcal{I}^{r,l,k}\rt\|_{L^{2n/(n-1)}}^2 \lt\| \pl_t^{m-r}\bar\pl^{j-l}\pl_n^{i-k}\pl v\rt\|_{L^{2n}}^2
\les \mathfrak{T} .\label{Feb8-6}
\end{align}

Note that
\begin{subequations}\label{Feb3}\begin{align}
& \lt\|\pl v \rt\|_{L^\iy}^2 \les \sum_{l+k\le [\iota]+n+2} \mathcal{D}^{0,l,k},  \label{Feb3-1} \\
 & \lt\| \pl_t \pl v \rt\|_{L^\iy}^2 +\lt\| \bar \pl \pl v \rt\|_{L^\iy}^2  +\lt\|  \pl_n \pl v \rt\|_{L^{q^*}}^2 \les \sum_{l+k\le [\iota]+n+3} \mathcal{D}^{0,l,k} , \label{Feb3-2}
\end{align}\end{subequations}
where the bounds for $\pl v$, $\pl_t \pl v $ and $\bar \pl \pl v$ follow from  $H^{\iota+[\iota]+n+3, \ [\iota]+n+2}\hookrightarrow H^{({n+1+[\iota]-\iota})/{2}}
\hookrightarrow L^\iy$, and that for $\pl_n \pl v $
 from  $H^{\iota+[\iota]+n+4, \ [\iota]+n+2}\hookrightarrow H^{({n+[\iota]-\iota})/{2}}
\hookrightarrow L^{q^*}$.
When $r+l+k =[\iota]+n+3$, we get for
$m+j+i=[\iota]+n+3 $,
\begin{align}
\mathcal{P}^{m,j,i}_{r,l,k}\le \lt\|  \sa^{\frac{\iota+i+1}{2}} \mathcal{I}^{m,j,i} \rt\|_{L^2}^2\lt\|\pl v\rt\|_{L^\iy}^2
\les \mathfrak{T}, \label{2-8-1}
\end{align}
and for $m+j+i=[\iota]+n+4 $,
\begin{align}
\mathcal{P}^{m,j,i}_{r,l,k}
\les & \lt\|  \sa^{\frac{\iota+i+1}{2}} \mathcal{I}^{m-1,j,i} \rt\|_{L^2}^2\lt\| \pl_t \pl v\rt\|_{L^\iy}^2+
 \lt\|  \sa^{\frac{\iota+i+1}{2}} \mathcal{I}^{m,j-1,i}  \rt\|_{L^2}^2\lt\| \bar\pl \pl v\rt\|_{L^\iy}^2
 \notag \\
& +
 \lt\|  \sa^{\frac{\iota+i+1}{2}} \mathcal{I}^{m,j,i-1}  \rt\|_{L^{2n/ (n-1)}}^2\lt\| \pl_n \pl v\rt\|_{L^{2n}}^2
\les  \mathfrak{T},
  \label{2-8-2}
\end{align}
where  \ef{Feb1} and  \ef{Feb3} have been used  to derive the estimates above.
When $r+l+k =[\iota]+n+4$, it follows from
$m+j+i =[\iota]+n+4$,
\ef{Feb1} and \ef{Feb3-1} that
\begin{align}\label{2-8-3}
\mathcal{P}^{m,j,i}_{r,l,k}\le \lt\|  \sa^{\frac{\iota+i+1}{2}} \mathcal{I}^{m,j,i} \rt\|_{L^2}^2\lt\|\pl v\rt\|_{L^\iy}^2
\les \mathcal{E} \sum_{a+b+c\le m+j+i-2}\mathcal{D}^{a,b,c}.
\end{align}
Finally, \ef{1.27-2} is a conclusion of \ef{1-27-1}-\ef{June7}, \ef{Feb8-6} and
\ef{2-8-1}-\ef{2-8-3}.

{\em Step 4}. In a similar way to deriving \ef{1.27-2}, we prove \ef{5-30-1} in this step.
 In what follows, we assume $(l,k)\in \mathfrak{S} \setminus \mathfrak{S}_1$, and set
$$X^{j,i}_{l,k}=\lt\|  \sa^{\frac{\iota+i+1}{2}} \mathcal{I}^{1,l,k} \bar\pl^{j-l}\pl_n^{i-k}\pl v\rt\|_{L^2}^2, \  \
 \mathscr{T}=\mathcal{E}(t)\sum_{l+k\le j+i-1}\mathcal{D}^{0,l,k} (t).$$
When $l+2k\le 1$, then $(l,k)=(1,0)$ and
\begin{align}
&X^{j,i}_{l,k}\le \lt\|  \mathcal{I}^{1,1,0} \rt\|_{L^\iy}^2 \lt\|  \sa^{\frac{\iota+i+1}{2}} \bar\pl^{j-1}\pl_n^{i}\pl v\rt\|_{L^2}^2
\les \mathscr{T},\label{5.31-a}
\end{align}
due to \ef{12.10}.
When $l+2k=2$, then $(l,k)=(2,0)$ and
\begin{align}
&X^{j,i}_{l,k}\le \lt\|  \mathcal{I}^{1,2,0} \rt\|_{L^{2n}}^2 \lt\|  \sa^{\frac{\iota+i+1}{2}} \bar\pl^{j-2}\pl_n^{i}\pl v\rt\|_{L^{{2n}/({n-1})}}^2
\les \mathscr{T},\label{5.31-b}
\end{align}
because of   \ef{12.10} and $H^{1,1}\hookrightarrow H^{1/2}\hookrightarrow L^{{2n}/({n-1})}$. To deal with the case of $l+2k\ge 3$, we set
\begin{align*}
&\mathfrak{S}_2=\{
(l,k)\in \mathfrak{S} \setminus \mathfrak{S}_1\big| 3\le l+2k, \ l+k\le [\iota]+n+1\},\\
&\mathfrak{S}_{21}=\{
(l,k)\in \mathfrak{S}_2 \big|  l+2k<\iota+i+3\},\\
&\mathfrak{S}_{22}=\mathfrak{S}_2\setminus\mathfrak{S}_{21}=\{
(l,k)\in \mathfrak{S}_2 \big|  l+2k \ge \iota+i+3\}.
\end{align*}
When $(l,k)\in \mathfrak{S}_{21}$, we use \ef{12.10}, and \ef{hard} $l+k-1$ times to achieve
$$X^{j,i}_{l,k}\le \lt\|  \sa^{\frac{l+2k-1}{2}} \mathcal{I}^{1,l,k} \rt\|_{L^\iy}^2 \lt\|  \sa^{\frac{\iota+i+2-l-2k}{2}} \bar\pl^{j-l}\pl_n^{i-k}\pl v\rt\|_{L^2}^2\les  \mathscr{T} .$$
If $(l,k)\in \mathfrak{S}_{22}$,  then $l+k =[\iota]+n$ or $[\iota]+n+1$ for $n=3$, and $l+k =[\iota]+n+1$ for $n=2$.
Consider the case of $n=3$ and $l+k =[\iota]+n$, in which $[\iota]=\iota$ and $k=i$, then one gets
$$X^{j,i}_{l,k}
\les \lt\|  \sa^{\frac{\iota+i+1}{2}} \mathcal{I}^{1,l,i} \rt\|_{L^{q_1}}^2 \lt\|\bar\pl^{j-l}\pl v\rt\|_{L^{n}}^2
\les \mathscr{T}, $$
noticing \ef{Feb1} and $H^{\iota+[\iota]+n,\ [\iota]+n-1}\hookrightarrow H^{({n-2})/{2}}
\hookrightarrow L^{n}$.
Consider the case of $n=3$ and $l+k =[\iota]+n+1$, in which $k=i$, or $[\iota]=\iota$ and $k=i-1$,  then we obtain
$$X^{j,i}_{l,k}
\les \lt\|  \sa^{\frac{\iota+i+1}{2}} \mathcal{I}^{1,l,i} \rt\|_{L^{q_2}}^2 \lt\|\bar\pl^{j-l}\pl v\rt\|_{L^{{2n}/({1+\iota-[\iota]})}}^2
\les \mathscr{T} , $$
with the help of \ef{Feb1} and $H^{\iota+[\iota]+n+1,\ [\iota]+n}\hookrightarrow H^{({n-1+[\iota]-\iota})/{2}}
\hookrightarrow L^{{2n}/({1+\iota-[\iota]})}$;
or
$$X^{j,i}_{l,k}\les \lt\|  \sa^{\frac{\iota+i+1}{2}} \mathcal{I}^{1,l,i-1} \rt\|_{L^{q_1}}^2 \lt\|\bar\pl^{j-l}\pl_n \pl v\rt\|_{L^{n}}^2
\les  \mathscr{T} , $$
due to \ef{Feb1} and $H^{\iota+[\iota]+n+2,\ [\iota]+n}\hookrightarrow H^{({n-2})/{2}}
\hookrightarrow L^{n}$.
Consider the case of $n=2$ and $l+k =[\iota]+n+1$, in which $[\iota]=\iota$ and $k=i$, then one has
$$X^{j,i}_{l,k}\les \lt\|  \sa^{\frac{\iota+i+1}{2}} \mathcal{I}^{1,l,i} \rt\|_{L^{2n}}^2 \lt\|\bar\pl^{j-l}\pl v\rt\|_{L^{{2n}}}^2
 \les \mathscr{T} , $$
noting \ef{Feb1} and $H^{\iota+[\iota]+n+1,\ [\iota]+n}\hookrightarrow H^{({n-1})/{2}}
\hookrightarrow L^{2n}$. So,
\begin{align}
X^{j,i}_{l,k}\les  \mathscr{T} ,  \ \  (l,k)\in \mathfrak{S}_2. \label{6-3-1}
\end{align}

When $l+k =[\iota]+n+2$, it follows from \ef{Feb1}, $H^{\iota+[\iota]+n+2, \ [\iota]+n+1}\hookrightarrow H^{(n+[\iota]-\iota)/2}\hookrightarrow L^{q^*}$ and $H^{\iota+[\iota]+n+3, \ [\iota]+n+1}\hookrightarrow H^{(n-1+[\iota]-\iota)/2}\hookrightarrow L^{{2n}/(1+\iota-[\iota])}$ that for
$j+i=[\iota]+n+3 $,
\begin{align}\label{6.5.1}
X^{j,i}_{l,k}\les \lt\|  \sa^{\frac{\iota+i+1}{2}} \mathcal{I}^{1,j-1,i} \rt\|_{L^q}^2 \lt\|  \bar\pl\pl v\rt\|_{L^{q^*}}^2\les  \mathscr{T},
\end{align}
and for $j+i=[\iota]+n+4 $,
\begin{align}
&X^{j,i}_{l,k}\les \lt\|  \sa^{\frac{\iota+i+1}{2}} \mathcal{I}^{1,j-2,i} \rt\|_{L^q}^2\lt\| \bar\pl^{2}\pl v\rt\|_{L^{q^*}}^2\notag\\
&+\lt\|  \sa^{\frac{\iota+i+1}{2}} \mathcal{I}^{1,j-1,i-1} \rt\|_{L^{q_2}}^2\lt\| \bar\pl\pl_n\pl v\rt\|_{L^{{2n}/({1+\iota-[\iota]})}}^2\notag\\
&+\lt\|  \sa^{\frac{\iota+i-1}{2}} \mathcal{I}^{1,j,i-2} \rt\|_{L^q}^2\lt\| \sa \pl_n^{2}\pl v\rt\|_{L^{q^*}}^2\les  \mathscr{T}, \label{6.5.2}
\end{align}
where $q=q_3$ for $[\iota]\neq \iota$, and $q=q_2$ for $[\iota]=\iota$.
When $l+k =[\iota]+n+3$, we get for $j+i =[\iota]+n+4$,
\begin{align}
X^{j,i}_{l,k}\le \lt\|  \sa^{\frac{\iota+i+1}{2}} \mathcal{I}^{1,j-1,i} \rt\|_{L^2}^2 \lt\| \bar\pl\pl v\rt\|_{L^\iy}^2\les  \mathscr{T}, \label{6.5.3}
\end{align}
due to \ef{Feb1} and \ef{Feb3-2}. So, \ef{5-30-1} is a conclusion of \ef{5.31-a}-\ef{6.5.3}.
\hfill $\Box$

\subsubsection{Curl estimates}

\begin{lem}It holds that for $m+j+i\le [\iota]+n+3$, and $t\in [0,T]$,
\begin{align}
\frac{d}{dt} \mathcal{V}^{m,j,i}(t) +  \mathcal{V}^{m,j,i}(t) \les
\mathcal{E} (t)\sum_{r+l+k\le m+j+i} \mathcal{D}^{r,l,k}(t). \label{Feb15}
\end{align}
\end{lem}

{\em Proof}. It follows from \ef{system-a} that
 \be\label{11-29}
{\rm curl}_x \pl_t v + {\rm curl}_x v =0.
\ee
For any nonnegative integers $m,\al$, and multi-index $\ba$, we take
$\pl_t^m\bar\pl^\ba \pl_n^\al$ onto \ef{11-29} to obtain
\begin{align}\label{1.19}
  \pl_t^{m+1}\bar\pl^\ba \pl_n^\al {\rm curl}_x v +  \pl_t^m\bar\pl^\ba \pl_n^\al {\rm curl}_x v =\pl_t^m\bar\pl^\ba \pl_n^\al  \lt[\pl_t, {\rm curl}_x\rt] v.
\end{align}
Integrate the inner product of \ef{1.19} and  $\sa^{\iota+\al+1} \pl_t^m\bar\pl^\ba \pl_n^\al {\rm curl}_x v$ over $\Oa$ and use  the Cauchy inequality to get
\begin{align}
&\frac{d}{dt} \int \sa^{\iota+\al+1}  \lt|\pl_t^{m}\bar\pl^\ba \pl_n^\al {\rm curl}_x v \rt|^2 dy +  \int \sa^{\iota+\al+1}  \lt|\pl_t^{m}\bar\pl^\ba \pl_n^\al {\rm curl}_x v \rt|^2 dy \notag \\
&\le \int \sa^{\iota+\al+1} \lt| \pl_t^m\bar\pl^\ba \pl_n^\al  \lt[\pl_t, {\rm curl}_x\rt] v \rt|^2 dy . \label{2-15}
\end{align}
Note that
\begin{align*}
&\lt| \pl_t^m\bar\pl^\ba \pl_n^\al  \lt[\pl_t, {\rm curl}_x\rt] v \rt|\\
\les & \sum_ {\substack{r\le m, \ l \le |\ba|, \  k \le \al }}
\mathcal{I}^{r+1,l,k}  \lt| \pl_t^{m-r}\bar\pl^{|\ba|-l} \pl_n^{\al -k} \pl v \rt|\\
= & \sum_ {\substack{1\le r\le m+1, \ l \le |\ba|, \  k \le \al}}
\mathcal{I}^{r,l,k}  \lt| \pl_t^{m+1-r}\bar\pl^{|\ba|-l} \pl_n^{\al-k} \pl v \rt|,
\end{align*}
then \ef{Feb15} is a conclusion of \ef{2-15} and \ef{1.27-2}. \hfill $\Box$

\subsubsection{Elliptic estimates near the bottom}\label{s3.3.3}

\begin{lem} It holds that for $\ m+j+i\le [\iota]+n+3$ and $i\ge 1$, and $t\in [0, T]$,
\begin{align}\label{div-est}
&\mathcal{D}^{m,j,i}_{1}(t)\les \sum_{ r+l\le m+j+i} \mathcal{D}^{r,l,0}(t) + \sum_{ r+l+k\le m+j+i-1} \mathcal{D}_2^{r,l,k}(t)\notag\\
&\qquad + \sum_{r+l+k\le m+j+i, \  k \ge 1} \lt\|\zeta_1{\rm curl}  \pl_t^r  \bar\pl^l \pl_n^k v \rt\|_{L^2}^2(t).
\end{align}

\end{lem}

{\em Proof}. Note that
\begin{align}
& \lt\|\za_1
 \pl_t^m\pl \bar\pl^j \pl_n^i v \rt\|_{L^2} \le \lt\|\pl(\za_1
 \pl_t^m \bar\pl^j \pl_n^i v) \rt\|_{L^2}+\lt\|\za_1'
 \pl_t^m \bar\pl^j \pl_n^i v \rt\|_{L^2}
 \notag \\
& \quad \les
\lt\|\za_1
 \pl_t^m \bar\pl^j \pl_n^i v \rt\|_{H^1}+ \lt\|
 \pl_t^m \bar\pl^j \pl_n^i v \rt\|_{L^2(\widetilde{\Oa})},
\notag\\
& \lt\|\zeta_1 \pl_t^m  \bar\pl^j \pl_n^i v \rt\|_{H^1}
\les  \lt\|\zeta_1 \pl_t^m  \bar\pl^j \pl_n^i v \rt\|_{L^2}
+\lt\|{\rm curl} \lt(\zeta_1 \pl_t^m  \bar\pl^j \pl_n^i v \rt) \rt\|_{L^2}
 \notag\\
& \quad +\lt\|{\rm div}\lt(\zeta_1 \pl_t^m  \bar\pl^j \pl_n^i v \rt) \rt\|_{L^2}+\lt\| \bar\pl \lt(\zeta_1 \pl_t^m  \bar\pl^j \pl_n^i v \rt) \rt\|_{L^2},
\label{1.15}\\
&  \lt\|
 \pl_t^m \bar\pl^j \pl_n^i v \rt\|_{L^2(\widetilde{\Oa})}
 \les \lt\|\za_2\sa^{\frac{\iota+i}{2}}
 \pl_t^m \pl \bar\pl^j \pl_n^{i-1} v \rt\|_{L^2(\widetilde{\Oa})}
 \le  \sqrt{ \mathcal{D}_2^{m,j,i-1} },  \notag
\end{align}
where
$\widetilde{\Omega}=\mathbb{T}^{n-1} \times ( \hbar/2, 3\hbar/4)$, and \ef{1.15}
follows from  the Hodge Decomposition \ef{Hodge}.
Then
\begin{align}
&\mathcal{D}_1^{m,j,i}
 \les \mathcal{D}_1^{m+1,j,i-1}
+\mathcal{D}_1^{m,j+1,i-1}
+\mathcal{D}_1^{m,j,i-1}
+\mathcal{D}^{m,j,i-1}_{2}
 \notag\\
& \quad +\lt\|\zeta_1{\rm curl}  \pl_t^m  \bar\pl^j \pl_n^i v \rt\|_{L^2}^2+\lt\|\zeta_1{\rm div} \pl_t^m  \bar\pl^j \pl_n^i v \rt\|_{L^2}^2 .
\label{1-15}
\end{align}
To control the $L^2$-norm of $\zeta_1{\rm div} \pl_t^m  \bar\pl^j \pl_n^i v $, we
divide \ef{system-a} by $\sa^\iota$, take
 $\pl_t$ onto the resulting equation  and  use $\pl_k  \sa=-\nu\da_{kn}$ to get
$$
 \ga \sa \pl_i {\rm div} v =\pl_{t}^2  v_i   + \pl_t  v_i+ \sa W_i   -(\iota+1)\nu\pl_t ( A^n_i J^{1-\ga}),
$$
where
\begin{align*}
&W_i=\pl_t \lt( J A^k_i\rt)  \pl_k J^{-\ga} -\ga J A^k_i \pl_k \lt( J^{-\ga} \lt(A^r_s-\da^r_s\rt) \pl_r v^s \rt)
\\
& \ \ -\ga J A^k_i \lt( \pl_k  J^{-\ga} \rt) {\rm div} v -\ga (J^{1-\ga} A^k_i -\da^k_i) \pl_k {\rm div} v.
\end{align*}
We take $\pl_t^m \bar\pl^j \pl_n^{i-1}$ ($i\ge 1$) onto the equation above to obtain
\begin{align*}
 & \ga  \sa \pl_t^m  \bar\pl^j \pl_n^{i}  {\rm div} v  = \ga (i-1) \nu \pl_t^m  \bar\pl^j \pl_n^{i-1}  {\rm div} v\\
 & \ \ + \pl_t^{m+2}  \bar\pl^j \pl_n^{i-1}   v_n   + \pl_t^{m+1}  \bar\pl^j \pl_n^{i-1}  v_n\\
 & \ \  + \pl_t^m  \bar\pl^j \pl_n^{i-1} \lt(\sa W_n \rt) -(\iota+1) \nu  \pl_t^{m+1}  \bar\pl^j \pl_n^{i-1}  ( A^n_n J^{1-\ga}),
\end{align*}
which implies
\begin{align}
&\lt\|\zeta_1{\rm div} \pl_t^m  \bar\pl^j \pl_n^i v  \rt\|_{L^2}^2=\lt\|\zeta_1 \pl_t^m  \bar\pl^j \pl_n^i {\rm div}  v  \rt\|_{L^2}^2\notag\\
& \les   \mathcal{D}_1^{m,j,i-1}+\mathcal{D}_1^{m+1,j,i-1}
+\lt\|\zeta_1\pl_t^{m+1}  \bar\pl^j \pl_n^{i-1}  ( A^n_n J^{1-\ga})\rt\|_{L^2}^2
\notag\\
& +\lt\|\zeta_1\pl_t^m  \bar\pl^j \pl_n^{i-1} W_n \rt\|_{L^2}^2+(i-1)^2\lt\|\zeta_1\pl_t^m  \bar\pl^j \pl_n^{i-2} W_n \rt\|_{L^2}^2.\label{1-16}
\end{align}

It follows from the Sobolev embedding and the a priori assumption \ef{1120} that
\be\label{1119}
\lt\|\pl_t^m \oa  \rt\|_{W^{[\iota]+n+3-m,\iy}(\Oa_1)}\les \lt\|\pl_t^m \oa \rt\|_{H^{[\iota]+n+5-m}(\Oa_1)}
\les  \sqrt{\mathcal{E}}
\le \ea_0,
\ee
where
$ \Oa_1=\mathbb{T}^{n-1} \times (0,3\hbar/4)$.
This implies that
for $1\le m+j+i\le [\iota]+n+2$,
\begin{align*}
&\lt\|\mathcal{I}^{m,j,i}\rt\|_{L^\iy(\Oa_1)}
\le  \lt\|\pl_t^{m}\bar\pl^{j}\pl_n^{i}
\pl\oa\rt\|_{L^\iy(\Oa_1)}
\\
& + \sum_ {\substack{ r\le m, \  l \le j, \  k \le i
\\ 1\le r+l+k \le m+j+i-1
  }} \lt\|\mathcal{I}^{r,l,k}\rt\|_{L^\iy(\Oa_1)}
\lt\|\pl_t^{m-r}\bar\pl^{j-l}\pl_n^{i-k}\pl\oa
\rt\|_{L^\iy(\Oa_1)}\\
& \les  \sqrt{\mathcal{E}}+ \ea_0 \sum_ {\substack{ r\le m, \ l \le j, \  k \le i
\\ 1\le r+l+k \le m+j+i-1
  }} \lt\|\mathcal{I}^{r,l,k}\rt\|_{L^\iy(\Oa_1)},
\end{align*}
which, together with the mathematical induction, gives that
\begin{align}\label{12.15}
&\sum_{1\le m+j+i\le [\iota]+n+2}\lt\|\mathcal{I}^{m,j,i}\rt\|_{L^\iy(\Oa_1)} \les \sqrt{\mathcal{E}}\les \ea_0.
\end{align}

In view of  \ef{12.15}, we see that  for $m+j+i \le [\iota]+n+3$,
\begin{align}
&\lt\|\zeta_1\pl_t^{m+1}  \bar\pl^j \pl_n^{i-1}  ( A^n_n J^{1-\ga})\rt\|_{L^2}^2
=\lt\|\zeta_1\pl_t^{m}  \bar\pl^j \pl_n^{i-1} \lt( \pl_t  ( A^n_n J^{1-\ga}) \rt) \rt\|_{L^2}^2
\notag\\
& \les
\sum_ {\substack{ r\le m, \ l \le j, \  k\le i-1 }} \lt\|\zeta_1\mathcal{I}^{r,l,k} \pl_t^{m-r} \bar\pl^{j-l}\pl_n^{i-1-k}\pl v\rt\|_{L^2}^2
\notag\\
& \les  \ea_0^2 \sum_ {\substack{r\le m, \ l \le j, \ k \le i-1 }} \lt\|\zeta_1 \pl_t^{m-r} \bar\pl^{j-l}\pl_n^{i-1-k}\pl v\rt\|_{L^2}^2
\notag\\
& \les  \ea_0^2 \sum_ {\substack{ r\le m, \  l \le j, \  k \le i-1 }}\mathcal{D}_1^{m-r,j-l,i-1-k}.
\label{Jan16}
 \end{align}
Note that
\begin{align*}
 &\lt|\pl_t^m  \bar\pl^j \pl_n^{i-1} W_n \rt| \les  \lt|\pl_t^m  \bar\pl^j \pl_n^{i-1}
\pl^2\oa \rt|\lt|\pl v\rt|
\\
 &\quad + \sum_ {\substack{ r_1+r_2 \le m, \  l_1+l_2 \le j,  \  k_1+k_2 \le i-1 \\
  r_2+l_2+k_2 \le m+j+i-2}} \mathcal{I}^{r_1,l_1,k_1}
\lt|\pl_t^{r_2} \bar\pl^{l_2}\pl_n^{k_2}\pl^2 \oa \rt| \\
 &\quad \times \lt|\pl_t^{m-r_1-r_2} \bar\pl^{j-l_1-l_2}\pl_n^{i-1-k_1-k_2}\pl v\rt| +\lt|\pl \oa \rt|\lt|\pl_t^{m} \bar\pl^{j}\pl_n^{i-1}\pl^2 v\rt|\\
 &\quad +\sum_ {\substack{ r\le m, \ l \le j, \  k \le i-1 \\
1\le r+l+k }} \mathcal{I}^{r,l,k} \lt|\pl_t^{m-r} \bar\pl^{j-l}\pl_n^{i-1-k}\pl^2 v\rt|
 ,
\end{align*}
then it  follows from \ef{1119} and \ef{12.15} that for $m+j+i \le [\iota]+n+3$,
\begin{align}\label{1.17}
&\lt\|\zeta_1\pl_t^m  \bar\pl^j \pl_n^{i-1}   W_n  \rt\|_{L^2}^2
\les \lt\| \zeta_1\pl_t^m  \bar\pl^j \pl_n^{i-1}
\pl^2\oa \lt|\pl v\rt|\rt\|_{L^2}^2
\notag\\
& \quad + \ea_0^2 \sum_{r\le m, \ l \le j, \ k \le i-1} \lt(\mathcal{D}_1^{r,l,k}
+\mathcal{D}_1^{r,l+1,k}
+\mathcal{D}_1^{r,l,k+1}\rt).
\end{align}
When $m+j+i \le [\iota]+n+2$, we have
\begin{align}\label{1.17-1}
 \lt\| \zeta_1\pl_t^m  \bar\pl^j \pl_n^{i-1}
\pl^2\oa \lt|\pl v\rt|\rt\|_{L^2}^2
\les \ea_0^2 \lt\| \zeta_1\pl v\rt\|_{L^2}^2
\le \ea_0^2 \mathcal{D}_1^{0,0,0},
\end{align}
due to \ef{1119}. When $m+j+i = [\iota]+n+3$, one gets
\begin{align}
& \lt\| \zeta_1\pl v\rt\|_{L^\iy}^2
 \les \lt\| \zeta_1\pl v\rt\|_{H^2}^2
 \les \sum_{r=1,2,3} \lt\| \zeta_1\pl^r v\rt\|_{L^2}^2
 +\sum_{r=1,2} \lt\| \pl^r v\rt\|_{L^2(\widetilde{\Oa})}^2
\notag \\
& \les \sum_{l+k\le  2} \mathcal{D}_1^{0,l,k}
 +\sum_{r=1,2} \lt\|\za_2 \sa^{\frac{\iota+r}{2}} \pl^r v\rt\|_{L^2(\widetilde{\Oa})}^2
 \les L, \label{Feb21.3}
\end{align}
where
$\widetilde{\Omega}=\mathbb{T}^{n-1} \times (\hbar/2, 3\hbar/4)$ and $L=\sum_{l+k\le  2} \mathcal{D}_1^{0,l,k}
 +\sum_{l+k\le 1} \mathcal{D}_2^{0,l,k}$;
which, together with \ef{1120}, implies that
\begin{align}\label{1.17-2}
& \lt\| \zeta_1\pl_t^m  \bar\pl^j \pl_n^{i-1}
\pl^2\oa \lt|\pl v\rt|\rt\|_{L^2}^2 \notag\\
& \les \lt\| \sa^{\frac{\iota+i+1}{2}} \pl_t^m  \bar\pl^j \pl_n^{i-1}
\pl^2\oa \rt\|_{L^2(\Oa_1)}^2  \lt\| \zeta_1\pl v\rt\|_{L^\iy}^2 \les \ea_0^2 L,
\end{align}
where
$ \Oa_1=\mathbb{T}^{n-1} \times (0,3\hbar/4)$.
Now,  it follows from \ef{1.17}-\ef{1.17-2} that
for $m+j+i \le [\iota]+n+3$,
\begin{align}\label{1-17}
&\lt\|\zeta_1\pl_t^m  \bar\pl^j \pl_n^{i-1}   W_n  \rt\|_{L^2}^2
\les \ea_0^2 \sum_{r\le m, \ l \le j, \ k \le i-1} (\mathcal{D}_1^{r,l,k} \notag \\
&\qquad +\mathcal{D}_1^{r,l+1,k}
+\mathcal{D}_1^{r,l,k+1})+\ea_0^2 \varsigma L,
\end{align}
where
$\varsigma=1$ when $m+j+i = [\iota]+n+3$, and
$\varsigma=0$ otherwise.
Similarly, we obtain
for $m+j+i \le [\iota]+n+3$ and $i\ge 2$,
\begin{align}\label{1-18}
&\lt\|\zeta_1\pl_t^m  \bar\pl^j \pl_n^{i-2}   W_n  \rt\|_{L^2}^2 \notag \\
& \les \ea_0^2 \sum_{r\le m, \ l \le j, \ k \le i-1} \lt(\mathcal{D}_1^{r,l,k}
+\mathcal{D}_1^{r,l+1,k}
+\mathcal{D}_1^{r,l,k+1}\rt).
\end{align}
It yields from \ef{1-15},  \ef{1-16}, \ef{Jan16}, \ef{1-17} and \ef{1-18} that
\begin{align*}
&\mathcal{D}_1^{m,j,i}
 \les \lt\|\zeta_1{\rm curl}  \pl_t^m  \bar\pl^j \pl_n^i v \rt\|_{L^2}^2 + \mathcal{D}_1^{m+1,j,i-1}
+\mathcal{D}_1^{m,j+1,i-1}
+\mathcal{D}_1^{m,j,i-1}
 \notag\\
& +\mathcal{D}^{m,j,i-1}_{2}
 +\ea_0^2 \sum_{r\le m, \ l \le j, \ k \le i-1} \lt(\mathcal{D}_1^{r,l,k}
+\mathcal{D}_1^{r,l+1,k}
+\mathcal{D}_1^{r,l,k+1}\rt)
 +\ea_0^2 \varsigma L,
\end{align*}
which, together with the smallness of  $\ea_0$, implies that
\begin{align}\label{Jan18}
& \mathcal{D}_1^{m,j,i}
 \les  \lt\|\zeta_1{\rm curl}  \pl_t^m  \bar\pl^j \pl_n^i v \rt\|_{L^2}^2 +  \mathcal{D}_1^{m+1,j,i-1}
+\mathcal{D}_1^{m,j+1,i-1}
\notag \\
& \qquad  +\sum_{r+l+k\le m+j+i-1} (\mathcal{D}_1^{r,l,k}
+\mathcal{D}_2^{r,l,k}).
\end{align}

Based on \ef{Jan18}, it is easy to
use  the  mathematical induction to prove \ef{div-est}.  Clearly, \ef{div-est} holds  for
 $m+j+i=1$  and $i\ge 1$,  due to  \ef{Jan18} which reads
\begin{align*}
\mathcal{D}_1^{0,0,1}
 \les \lt\|\zeta_1{\rm curl} \pl_n v \rt\|_{L^2}^2 + \mathcal{D}_1^{1,0,0}
+\mathcal{D}_1^{0,1,0}
+\mathcal{D}_1^{0,0,0}
+\mathcal{D}^{0,0,0}_{2}
 .
\end{align*}
Suppose that \ef{div-est} holds for $m+j+i\le l-1$ and  $i\ge 1$, then one has for  $m+j+i=l$ and  $i\ge 1$,
\begin{align*}
& \mathcal{D}_1^{m,j,i}
 \les   \mathcal{D}_1^{m+1,j,i-1}
+\mathcal{D}_1^{m,j+1,i-1}
 +\sum_{ r+l\le l-1} \mathcal{D}^{r,l,0}
\notag\\
 &\quad
  +\sum_{r+l+k\le l-1}
\mathcal{D}_2^{r,l,k}+\sum_{r+l+k\le l, \  k \ge 1} \lt\|\zeta_1{\rm curl}  \pl_t^r  \bar\pl^l \pl_n^k v \rt\|_{L^2}^2,
\end{align*}
which proves \ef{div-est} from $i=1$ to $i=l$ step by step.
\hfill $\Box$

\subsubsection{Energy estimates}\label{s3.3.4}

\begin{lem}There exists a  positive constant $\da\le 1/4$ which only depends on $\ga$ and $M$  such that
\begin{align}\label{n21-4-8}
e^{\da t} \mathcal{E}_I(t)+
\int_0^t e^{\da s} \mathcal{E}_I(s) ds
\les \mathcal{E}_I(0), \ \  t\in [0,T].
\end{align}
\end{lem}

{\em Proof}. Take $\pl_t$ onto \ef{system-a}, divide the resulting equation by $\sa^\iota$ and use the Piola identity $\pl_k(JA^k_i)=0$ $(1\le i\le n)$ to get
\be\label{11-27}
 \pl_{t}^2 v_i   +  \pl_t v_i -   \sa \pl_k {H}^k_i -(\iota+1)(\pl_k \sa){H}^k_i   =0,
\ee
where
\begin{align*}
{H}^k_i=J^{1-\ga} \lt( A^k_r A^s_i \pl_s v^r + \iota^{-1} A^k_i {\rm div}_x v \rt).
\end{align*}
For any nonnegative integers $m,\al$, and multi-index $\ba$, we take
$\pl_t^m\bar\pl^\ba \pl_n^\al$ onto \ef{11-27} to obtain
\begin{align}\label{81-1}
&\pl_{t}^{m+2} \bar\pl^\ba \pl_n^\al v_i   +   \pl_t^{m+1}  \bar\pl^\ba \pl_n^\al v_i -  \sa^{-\iota-\al} \pl_k \lt( \sa^{\iota+\al+1} \pl_t^m\bar\pl^\ba \pl_n^\al {H}^k_i \rt)\notag\\
& =\al \lt((\pl_n \sa) \pl_k - (\pl_k \sa) \pl_n\rt) \pl_t^m\bar\pl^\ba \pl_n^{\al-1}{H}^k_i={R}_{1;i}^{m,\ba,\al}.
\end{align}
Let $\zeta_2=\zeta_2(y_n)$ be a smooth cut-off function satisfying \ef{David} and set
\be\label{1-25}
\chi=1 \ \  {\rm when} \ \  \al=0, \ \ {\rm and} \ \  \chi=\zeta_2 \ \ {\rm when} \ \ \al\ge 1.
\ee

We integrate the product of \ef{81-1}  and $\chi^2 \sa^{\iota+\al} \pl_t^{m+1}  \bar\pl^\ba \pl_n^\al v^i$ over $\Oa$ and use the boundary condition $\sa=0$ on $\{y_n=\hbar\}$ to get
\begin{align*}
&\frac{1}{2}\frac{d}{dt}\int \chi^2 \sa^{\iota+\al}  \lt| \pl_t^{m+1}  \bar\pl^\ba \pl_n^\al v \rt|^2 dy   +  \int \chi^2 \sa^{\iota+\al}  \lt| \pl_t^{m+1}   \bar\pl^\ba \pl_n^\al v \rt|^2 dy\\
&+\sum_{1\le j \le 3} L_j^{m,\ba,\al}=0,
\end{align*}
where
\begin{align*}
 L_1^{m,\ba,\al}=& \int_{\{y_n=0\}}  \chi^2  \sa^{\iota+\al+1} \lt( \pl_t^m \bar\pl^\ba \pl_n^\al {H}^n_i \rt)  \pl_t^{m+1}  \bar\pl^\ba \pl_n^\al v^i  dy_* ,\\
 L_2^{m,\ba,\al}=& \int  \chi^2  \sa^{\iota+\al+1} \lt( \pl_t^m\bar\pl^\ba \pl_n^\al {H}^k_i \rt)  \pl_k \pl_t^{m+1}  \bar\pl^\ba \pl_n^\al v^i   dy,\\
 L_3^{m,\ba,\al}=&2\int \chi( \pl_k \chi )  \sa^{\iota+\al+1} \lt( \pl_t^m\bar\pl^\ba \pl_n^\al {H}^k_i \rt)   \pl_t^{m+1}  \bar\pl^\ba \pl_n^\al v^i   dy \\
&- \int \chi^2 \sa^{\iota+\al} {R}_{1;i}^{m,\ba,\al} \pl_t^{m+1}  \bar\pl^\ba \pl_n^\al v^i  dy.
\end{align*}
Clearly, the boundary term $L_1^{m,\ba,\al}=0$ for $\al \ge 1$. When $\al =0$,
it follows from the boundary condition \ef{system-c}:  $v_n |_{\{y_n=0\}}=0$ that $ x_n  |_{\{y_n=0\}}=0$, which implies $JA^n_i |_{\{y_n=0\}}=0$ and hence ${H}^n_i |_{\{y_n=0\}}=0$ for $1\le i\le n-1$, so that $L_1^{m,\ba,0}=0$.
We set
\begin{align*}
{R}_{2;i}^{m,\ba,\al;k}=
\pl_t^m \bar\pl^\ba \pl_n^\al \lt(J^{1-\ga} \lt( A^k_r A^s_i \pl_s v^r + \iota^{-1} A^k_i {\rm div}_x v \rt) \rt)\\
 -J^{1-\ga} \lt( A^k_r A^s_i \pl_s \pl_t^m \bar\pl^\ba \pl_n^\al v^r + \iota^{-1} A^k_i {\rm div}_x \pl_t^m \bar\pl^\ba \pl_n^\al v \rt),
\end{align*}
and use \ef{nabt} to obtain
\begin{align*}
L_2^{m,\ba,\al}=&\frac{1}{2}\frac{d}{dt}\sum_{1\le j \le 3} E_j^{m,\ba,\al}
+L_4^{m,\ba,\al},
\end{align*}
where
\begin{align*}
&E_1^{m,\ba,\al}= \int \chi^2 \sa^{\iota+\al+1}  J^{1-\ga}
\lt( \lt|\na_x \pl_t^m\bar\pl^\ba \pl_n^\al v \rt|^2  + \iota^{-1}  \lt|{\rm div}_x \pl_t^m\bar\pl^\ba \pl_n^\al v \rt|^2   \rt) dy,\\
&E_2^{m,\ba,\al}= -\int \chi^2 \sa^{\iota+\al+1}  J^{1-\ga} \lt|{\rm curl}_x \pl_t^m\bar\pl^\ba \pl_n^\al v \rt|^2  dy,
\\
&E_3^{m,\ba,\al}= 2\int \chi^2 \sa^{\iota+\al +1} {R}_{2;i}^{m,\ba,\al;k}   \pl_k \pl_t^{m}  \bar\pl^\ba \pl_n^\al v^i  dy,\\
& L_4^{m,\ba,\al}=  \int \chi^2 \sa^{\iota+\al+1}  J^{1-\ga}  \lt[ \nabla_x \pl_t^m\bar\pl^\ba \pl_n^\al v^r \rt]_i \lt[\nabla_x v^s\rt]_r\lt[\nabla_x \pl_t^m \bar\pl^\ba \pl_n^\al v^i\rt]_s dy \\
&\ -\iota^{-1} \int \chi^2  \sa^{\iota+\al+1}  J^{1-\ga}
\lt(\pl_t A^k_i \rt) \lt({\rm div}_x \pl_t^m\bar\pl^\ba \pl_n^\al v \rt) \pl_k  \pl_t^m\bar\pl^\ba \pl_n^\al v^i dy \\
&\ -\frac{1}{2}\int \chi^2  \sa^{\iota+\al+1} \lt(\pl_t  J^{1-\ga} \rt) \lt( \lt|\na_x \pl_t^m\bar\pl^\ba \pl_n^\al v \rt|^2  + \iota^{-1}  \lt|{\rm div}_x \pl_t^m\bar\pl^\ba \pl_n^\al v \rt|^2   \rt)
 dy\\
 &\  +\frac{1}{2}\int \chi^2  \sa^{\iota+\al+1}\lt(\pl_t  J^{1-\ga} \rt)  \lt|{\rm curl}_x \pl_t^m\bar\pl^\ba \pl_n^\al v \rt|^2 dy\\
 &\ -\int \chi^2 \sa^{\iota+\al+1} \lt(\pl_t {R}_{2;i}^{m,\ba,\al;k} \rt) \pl_k  \pl_t^{m}  \bar\pl^\ba \pl_n^\al v^i dy.
\end{align*}
So, we have
\begin{align}
&\frac{1}{2}\frac{d}{dt} \int \chi^2 \sa^{\iota+\al}  \lt| \pl_t^{m+1}  \bar\pl^\ba \pl_n^\al v \rt|^2 dy+  \int \chi^2 \sa^{\iota+\al}  \lt| \pl_t^{m+1}   \bar\pl^\ba \pl_n^\al v \rt|^2 dy \notag  \\
& +\frac{1}{2}\frac{d}{dt}\sum_{1\le j\le 3} E_j^{m,\ba,\al}+\sum_{j=3, 4} L_j^{m,\ba,\al}=0. \label{86-1}
\end{align}

Similarly, we
integrate the product of \ef{81-1}  and $\chi^2 \sa^{\iota+\al} \pl_t^{m}  \bar\pl^\ba \pl_n^\al v^i$ over $\Oa$, and use  the boundary conditions and \ef{nab} to achieve
\begin{align}
&\frac{1}{2}\frac{d}{dt}\int \chi^2 \sa^{\iota+\al}  \lt(\lt| \pl_t^{m}  \bar\pl^\ba \pl_n^\al v \rt|^2 + 2 \lt(\pl_t^{m+1}  \bar\pl^\ba \pl_n^\al v_i\rt)   \pl_t^{m}  \bar\pl^\ba \pl_n^\al v^i \rt)dy \notag\\
& - \int \chi^2 \sa^{\iota+\al}  \lt| \pl_t^{m+1}   \bar\pl^\ba \pl_n^\al v \rt|^2 dy+ \sum_{j=1,2} E_j^{m,\ba,\al} =-2\sum_{j=5,6}L_j^{m,\ba,\al},  \label{86-2}
\end{align}
where
\begin{align*}
L_5^{m,\ba,\al}=&\frac{1}{2}\int \chi^2 \sa^{\iota+\al+1}   {R}_{2;i}^{m,\ba,\al;k}\pl_k  \pl_t^m \bar\pl^\ba \pl_n^\al v^i  dy,\\
L_6^{m,\ba,\al}=& \int \chi (\pl_k \chi)  \sa^{\iota+\al+1} \lt( \pl_t^m\bar\pl^\ba \pl_n^\al {H}^k_i \rt)  \pl_t^{m}  \bar\pl^\ba \pl_n^\al v^i  dy\\
&- \frac{1}{2} \int \chi^2 \sa^{\iota+\al} {R}_{1;i}^{m,\ba,\al} \pl_t^{m}  \bar\pl^\ba \pl_n^\al v^i  dy .
\end{align*}
It produces from $2\times\ef{86-1}+\ef{86-2}$ that
\begin{align}\label{Feb20}
\frac{d}{dt} \mathfrak{E}^{m,\ba,\al}(t)+ \mathfrak{D}^{m,\ba,\al}(t)
=-E_2^{m,\ba,\al}-2\sum_{3\le j \le 6} L_j^{m,\ba,\al} ,
\end{align}
where
\begin{align*}
&\mathfrak{E}^{m,\ba,\al}
=\frac{1}{2}\int \chi^2 \sa^{\iota+\al}  \lt(2 \lt| \pl_t^{m+1}  \bar\pl^\ba \pl_n^\al v \rt|^2 + \lt| \pl_t^{m}  \bar\pl^\ba \pl_n^\al v \rt|^2 \rt. \notag\\
&\quad \lt.+ 2 \lt(\pl_t^{m+1}  \bar\pl^\ba \pl_n^\al v_i \rt)   \pl_t^{m}  \bar\pl^\ba \pl_n^\al v^i\rt)dy
+\sum_{1\le j\le 3} E_j^{m,\ba,\al},\\
&\mathfrak{D}^{m,\ba,\al}=\int \chi^2 \sa^{\iota+\al}  \lt| \pl_t^{m+1}   \bar\pl^\ba \pl_n^\al v \rt|^2 dy+
E_1^{m,\ba,\al}.
\end{align*}

{\em Step 1}. In this step, we prove that
for $m+|\ba|+\al\le [\iota]+n+3$,
\begin{align}\label{Feb21}
\frac{d}{dt} \mathfrak{E}^{m,\ba,\al}+ \mathfrak{D}^{m,\ba,\al}
\les \mathcal{M}^{m,|\ba|,\al} ,
\end{align}
where
\begin{align*}
& \mathcal{M}^{m, |\ba|, 0} =  \mathcal{V}^{m,|\ba|,0} + \ea_0 \sum_{r+l+k\le m+|\ba| }\mathcal{D}^{r,l,k},\\
& \mathcal{M}^{m,|\ba|,\al}= \mathcal{D}_2^{m,|\ba|+1,\al-1} +\mathcal{D}_2^{m+1,|\ba|,\al-1}+\sum_{r+l\le m+|\ba|+\al}
 \mathcal{D}^{r,l,0}\\
&\qquad  + \sum_{r+l+k\le m+|\ba|+\al }(\mathcal{V}^{r,l,k} + \ea_0\mathcal{D}^{r,l,k})
 \\
 &\qquad   +\sum_{\substack{ r+l+k\le m+|\ba|+\al-1, \  k \ge 1 } } \mathcal{D}_2^{r,l,k} , \ \ \al \ge 1.
\end{align*}

It needs to analyze the terms on the right-hand side of \ef{Feb20}.
Note that
\begin{align*}
&\lt| {R}_{2;i}^{m,\ba,\al;k} \rt|\les
 \sum_ {\substack{ r\le m, \ l \le |\ba|, \ b \le \al
\\ 1 \le r+l+b
  }} \mathcal{I}^{r,l,b} \lt|\pl_t^{m-r}\bar\pl^{|\ba|-l}\pl_n^{\al-b}\pl v\rt|
,\\
  &\lt|\pl_t {R}_{2;i}^{m,\ba,\al;k}\rt|\les
 \sum_ {\substack{ r\le m+1, \  l \le |\ba|, \  b \le \al
\\ 1 \le r+l+b
  }} \mathcal{I}^{r,l,b} \lt|\pl_t^{m+1-r}\bar\pl^{|\ba|-l}\pl_n^{\al-b}\pl v\rt|,
\end{align*}
then it follows from the Cauchy inequality, \ef{12-8} and \ef{1.27-2} that
\begin{align}
&\lt|E_3^{m,\ba,\al}\rt|
\les  \ea_0 \int \chi^2 \sa^{\iota+\al+1}  \lt|\pl \pl_t^m\bar\pl^\ba \pl_n^\al v \rt|^2
\notag\\
&\qquad\qquad +  \ea_0 \sum_{r+l+k\le m+|\ba|+\al -1}\mathcal{D}^{r,l,k},
\label{Feb20.1}\\
&\sum_{j=4,5} \lt|L_j^{m,\ba,\al}\rt|
\les  \ea_0 \sum_{r+l+k\le m+|\ba|+\al }\mathcal{D}^{r,l,k}.
\label{Feb20.2}
\end{align}

When
 $\al \ge 1$, notice that
\begin{align*}
&  \lt| \pl_t^m\bar\pl^\ba \pl_n^\al {H}^k_i \rt|
\les
\sum_ {\substack{ r\le m, \ l \le |\ba|, \ b \le \al
  }} \mathcal{I}^{r,l,b}  \lt|\pl_t^{m-r}\bar\pl^{|\ba|-l}\pl_n^{\al-b}\pl v\rt|,\\
&\lt|{R}_{1;i}^{m,\ba,\al} \rt|
\les
\sum_ {\substack{ r\le m, \ l \le |\ba|+1, \ b \le \al-1
  }} \mathcal{I}^{r,l,b} \lt|\pl_t^{m-r}\bar\pl^{|\ba|+1-l}\pl_n^{\al-1-b}\pl v\rt|,
\end{align*}
then it yields from \ef{12.15} and \ef{1.27-2} that
\begin{align}
& \int_{\Oa_2} \lt| \pl_t^m\bar\pl^\ba \pl_n^\al {H}^k_i \rt|^2  dy
\le \int \zeta_1^2  \lt| \pl_t^m\bar\pl^\ba \pl_n^\al {H}^k_i \rt|^2  dy \notag \\
&  \les \mathcal{D}_1^{m, |\ba|, \al}
 +\ea_0^2 \sum_{r+l+b\le m+|\ba|+\al-1}
 \mathcal{D}^{r,l,b}, \label{Feb21.1}\\
 &\int \chi^2 \sa^{\iota+\al}\lt|{R}_{1;i}^{m,\ba,\al} \rt|^2 dy
\les
  \mathcal{D}_2^{m,|\ba|+1,\al-1}
\notag  \\
 & \qquad +\ea_0^2 \sum_{r+l+k\le m+|\ba|+\al -1}\mathcal{D}^{r,l,k}, \label{Feb21.2}
\end{align}
where
$ \Oa_2=\mathbb{T}^{n-1} \times (\hbar/4, \hbar/2)$. Indeed, the following estimate has been used to derive \ef{Feb21.1}.
\begin{align*}
& \int \lt|\zeta_1\mathcal{I}^{m,|\ba|,\al}\pl v\rt|^2 dy
\les  \lt\|\za_1 \pl v \rt\|_{L^\iy}^2 \int \sa^{\iota+\al+1} \lt|\pl_t^m \bar\pl^{|\ba|} \pl_n^\al \pl \oa \rt|^2 dy\\
& \quad +
\sum_ {\substack{r\le m, \ l \le |\ba|, k \le \al
\\ 1\le r+l+k \le m+j+i-1
  }} \int \zeta_1^2  \lt|\mathcal{I}^{r,l,k}\pl_t^{m-r}
  \bar\pl^{|\ba|-l}\pl_n^{\al-k}\pl\oa\rt|^2 |\pl v|^2 dy\\
& \les \ea_0^2  \lt\|\za_1 \pl v \rt\|_{L^\iy}^2
+ \ea_0^4 \mathcal{D}_1^{0,0,0} \les
\ea_0^2 \sum_{r+l+k\le m+|\ba|+\al-1}
 \mathcal{D}^{r,l,k}
\end{align*}
for $m+|\ba|+\al=[\iota]+n+3$, where the first inequality follows from \ef{n5.30}, the second from \ef{1120}, \ef{1119} and \ef{12.15}, and the last from \ef{Feb21.3}. It produces from the Cauchy inequality, \ef{Feb21.1} and \ef{Feb21.2} that for $\al \ge 1$,
\begin{align}
& \sum_{j=3,6} \lt| L_j^{m,\ba,\al}\rt|
\les \mathcal{D}_1^{m, |\ba|, \al}
+\mathcal{D}_2^{m,|\ba|+1,\al-1} +\mathcal{D}_2^{m+1,|\ba|,\al-1}\notag\\
& \qquad  +\mathcal{D}^{m,|\ba|,\al-1}
 +\ea_0^2 \sum_{r+l+k\le m+|\ba|+\al-1}
 \mathcal{D}^{r,l,k}. \label{Feb21.4}
\end{align}

In view of \ef{1.27-2} and
\begin{align}
&|  {\rm curl}_x \pl_t^m\bar\pl^\ba \pl_n^\al v - \pl_t^m\bar\pl^\ba \pl_n^\al {\rm curl}_x v|\notag\\
&\les \sum_ {\substack{ r\le m, \ l \le |\ba|, \ k \le \al,
\  1 \le r+l+k
  }} \mathcal{I}^{r,l,k}| \pl_t^{m-r}\bar\pl^{|\ba|-l}\pl_n^{\al-k}\pl v|,
  \label{21June16-1}
\end{align}
we see that
\begin{align}
 \lt\| \sa^{\frac{\iota+\al+1}{2}} {\rm curl}_x \pl_t^m\bar\pl^\ba \pl_n^\al v \rt\|_{L^2}^2\les \mathcal{V}^{m,|\ba|, \al} +
 \ea_0^2 \sum_{r+l+k\le m+|\ba|+\al -1}\mathcal{D}^{r,l,k},\label{Feb25-1}
\end{align}
Due to \ef{7.9}, we have
\begin{align*}
|  {\rm curl}_x \pl_t^m\bar\pl^\ba \pl_n^\al v - {\rm curl} \pl_t^m\bar\pl^\ba \pl_n^\al  v|
\les \ea_0 |\pl \pl_t^m\bar\pl^\ba \pl_n^\al  v|,
\end{align*}
which, together with \ef{Feb25-1}, implies
\begin{align}\label{Feb25-2}
& \lt\| \sa^{\frac{\iota+\al+1}{2}} {\rm curl} \pl_t^m\bar\pl^\ba \pl_n^\al v \rt\|_{L^2}^2\notag \\
& \les \mathcal{V}^{m,|\ba|, \al} + \ea_0^2 \mathcal{D}^{m,|\ba|, \al}
 +
 \ea_0^2 \sum_{r+l+k\le m+|\ba|+\al-1 }\mathcal{D}^{r,l,k}.
\end{align}
So, \ef{Feb21} is a conclusion of \ef{Feb20},  \ef{Feb20.2}, \ef{Feb21.4}, \ef{Feb25-1}, \ef{Feb25-2}, \ef{div-est} and $L_3^{m,\ba,0}=L_6^{m,\ba,0}=0$.

{\em Step 2}. In this step, we prove that there exist positive constants $K, K_b,  K_{b,\al}$ only depending  on $\ga$ and $ M$ such that
\begin{align}\label{Feb27-3}
\frac{d}{dt} \mathfrak{E}(t) +  \frac{1}{2} \mathfrak{D}(t) + \mathcal{V}(t)
 \le 0 ,
\end{align}
where
\begin{align*}
&\mathfrak{D}=\sum_{m+|\ba|+\al\le [\iota]+n+3} \mathfrak{D}^{m,\ba,\al}, \ \  \mathcal{V}= \sum_{m+j+i\le [\iota]+n+3  } \mathcal{V}^{m,j,i},
\\
&\mathfrak{E}= \sum_{b\le [\iota]+n+3}K_b\sum_{\al\le b} K_{b,\al} \sum_{m+|\ba|=b-\al} \mathfrak{E}^{m,\ba,\al} +K\mathcal{V}
 \end{align*}
satisfy the following estimates:
\begin{align} \label{Feb28-1}
& \mathfrak{E} \les
\mathfrak{D}+ \widetilde{\mathfrak{D}}  + \mathcal{V}\ \ {\rm and} \ \
\mathfrak{E}  \ge
6^{-1}\lt(2^{-1}\mathfrak{D}+ \widetilde{\mathfrak{D}} \rt)
+ \mathcal{V}\\
& {\rm with} \ \  \widetilde{\mathfrak{D}}=\sum_{m+|\ba|+\al \le [\iota]+n+3}\lt\|\chi\sa^{\frac{\iota+\al}{2}}\pl_t^m \bar\pl^\ba \pl_n^\al v\rt\|_{L^2}^2. \notag
\end{align}

Let $b\le [\iota]+n+3$ and $\al\le b$ be  non-negative integer, and $K_{b,\al}\ge 1$ be constants to be determined later.
It follows from \ef{Feb21} that
\begin{align}\label{Feb26}
&\frac{d}{dt}\sum_{\al\le b} K_{b,\al} \sum_{m+|\ba|=b-\al} \mathfrak{E}^{m,\ba,\al}+ D_b
 \les N_b^I+ N_b^{II},
\end{align}
where
\begin{align*}
D_b= &\sum_{\al\le b} K_{b,\al} \sum_{m+|\ba|=b-\al} \mathfrak{D}^{m,\ba,\al} -C \sum_{1\le \al \le b} K_{b, \al} \sum_{r+l= b} \mathcal{D}^{r,l,0}
\\
&-C\sum_{2 \le \al \le b} K_{b, \al} \sum_{m+|\ba|=b+1-\al} \mathcal{D}_2^{m,|\ba|, \al-1},\\
N_b^{I}= & \sum_{\al \le b} K_{b,\al} \sum_{r+l+k\le b }(\mathcal{V}^{r,l,k} + \ea_0\mathcal{D}^{r,l,k}),\\
N_b^{II} =&  \sum_{1\le \al \le b} K_{b,\al} \lt(\sum_{\substack{ r+l \le b-1 } } \mathcal{D}^{r,l,0} +\sum_{\substack{ r+l+k\le b-1, \  k\ge 1 } } \mathcal{D}_2^{r,l,k}
   \rt).
\end{align*}
Due to \ef{6.7-1c}, we have
\begin{align}\label{Feb27}
\mathcal{D}^{m,j,0}\les \sum_{|\ba|=j}\mathfrak{D}^{m,\ba,0}, \ \
\mathcal{D}_2^{m,j,i}\les \sum_{|\ba|=j}\mathfrak{D}^{m,\ba,i} \ \ {\rm for}  \ \ i \ge 1,
\end{align}
and then
\begin{align*}
D_b \ge & \lt(K_{b,0}-C\sum_{1\le \al \le b}K_{b,\al}\rt)\sum_{m+|\ba|=b} \mathfrak{D}^{m,\ba,0}
\\
&+\sum_{1\le \al\le b}(K_{b,\al}-C K_{b, \al+1 })  \sum_{m+|\ba|=b-\al} \mathfrak{D}^{m,\ba,\al}
,\\
N_b^{II}\les & \sum_{1\le \al \le b} K_{b,\al}
\sum_{m+|\ba|+\al\le b-1} \mathfrak{D}^{m,\ba,\al},
\end{align*}
where $K_{b,b+1}=0$.
We may choose $K_{b,b}=1$ and suitable large constants $K_{b,\al}$ $(\al\le b-1)$ only depending on $\ga$ and $M$ such that
$$D_b \ge \sum_{m+|\ba|+\al=b} \mathfrak{D}^{m,\ba,\al},  \ \ N_b^{II}\les
\sum_{m+|\ba|+\al\le b-1} \mathfrak{D}^{m,\ba,\al},$$
which, together with \ef{Feb26}, implies
\begin{align*}
&\frac{d}{dt}\sum_{\al\le b} K_{b,\al} \sum_{m+|\ba|=b-\al} \mathfrak{E}^{m,\ba,\al}+ \sum_{m+|\ba|+\al=b} \mathfrak{D}^{m,\ba,\al}\notag \\
& \les  N_b^{I}+ \sum_{m+|\ba|+\al\le b-1} \mathfrak{D}^{m,\ba,\al}.
\end{align*}
Similarly, we can choose $K_{[\iota]+n+3}=1$ and suitable large constants $K_b\ge 1$ ($b \le [\iota]+n+2$) only  depending  on $\ga$ and $M$ such that
\begin{align}\label{Feb27-1}
&\frac{d}{dt}\sum_{b\le [\iota]+n+3}K_b\sum_{\al\le b} K_{b,\al} \sum_{m+|\ba|=b-\al} \mathfrak{E}^{m,\ba,\al}\notag\\
& + \sum_{m+|\ba|+\al\le [\iota]+n+3} \mathfrak{D}^{m,\ba,\al}
 \les  \sum_{b\le [\iota]+n+3} N_b^{I} .
\end{align}

It follows from the Cauchy inequality, \ef{Feb20.1}, \ef{Feb25-1} and the smallness of $\ea_0$ that
\begin{subequations}\label{Feb28}\begin{align}
\mathfrak{E}^{m,\ba,\al} & \le
2\mathfrak{D}^{m,\ba,\al}+ \lt\|\chi\sa^{\frac{\iota+\al}{2}}\pl_t^m \bar\pl^\ba \pl_n^\al v\rt\|_{L^2}^2
\notag \\
& +C\ea_0\sum_{r+l+k\le m+|\ba|+\al-1} \mathcal{D}^{r,l,k},\\
\mathfrak{E}^{m,\ba,\al} & \ge
6^{-1}\lt(\mathfrak{D}^{m,\ba,\al}+ \lt\|\chi\sa^{\frac{\iota+\al}{2}}\pl_t^m \bar\pl^\ba \pl_n^\al v\rt\|_{L^2}^2\rt)
\notag \\
& -C \mathcal{V}^{m,|\ba|,\al}
-C \ea_0\sum_{r+l+k\le m+|\ba|+\al-1} \mathcal{D}^{r,l,k}.
\end{align}\end{subequations}
In view of \ef{div-est} and \ef{Feb25-2}, we see that
\begin{align*}
 \mathcal{D}^{m,j,i}
 \les  \sum_{ r+l\le m+j+i} \mathcal{D}^{r,l,0} + \sum_{ r+l+k\le m+j+i, \ k\ge 1 } \mathcal{D}_2^{r,l,k}\notag\\
  + \sum_{r+l+k\le m+j+i }   \mathcal{V}^{r,l,k} + \ea_0^2 \sum_{r+l+k\le m+j+i } \mathcal{D}^{r,l,k},
\end{align*}
which, together with the smallness of $\ea_0$ and \ef{Feb27}, implies that
\begin{align}\label{Feb27-2}
 \mathcal{D}^{m,j,i}
 \les  \sum_{ r+|\ba|+k \le m+j+i} \mathfrak{D}^{r,\ba,k}+ \sum_{r+l+k\le m+j+i }   \mathcal{V}^{r,l,k}   .
\end{align}
So, it produces from
\ef{Feb15}, \ef{Feb27-1}-\ef{Feb27-2} and the smallness of $\ea_0$ that there exists a suitable large constant $K\ge 1$ only depending on $\ga$ and $M$ such that \ef{Feb28-1} and \ef{Feb27-3} hold.

{\em  Step 3}. This step devotes to proving \ef{n21-4-8}.  Clearly, we have
\begin{align*}
 \widetilde{\mathfrak{D}}\les & \sum_{\substack{m+|\ba|+\al \le [\iota]+n+3,\
 1\le m   }}
 \mathfrak{D}^{m-1,\ba,\al}
+\sum_{\substack{|\ba|+\al \le [\iota]+n+3,\
 1\le \al  }}
 \mathfrak{D}^{0,\ba,\al-1}
 \notag \\
& +\sum_{\substack{1\le |\ba| \le [\iota]+n+3  }}
 \lt\|\sa^{ {\iota}/{2}}\bar\pl^\ba v\rt\|_{L^2}^2
 +\lt\|\sa^{ {\iota}/{2}} v\rt\|_{L^2}^2.
\end{align*}
Notice that for $|\ba|\ge 1$,
\begin{align*}
 \lt\|\sa^{ {\iota}/{2}}\bar\pl^\ba v\rt\|_{L^2}^2
 \les & \lt\|\sa^{ {\iota }/{2}+1}\bar\pl^\ba v\rt\|_{L^2}^2
 + \lt\|\sa^{ {\iota}/{2}+1}\bar\pl^\ba \pl_n v\rt\|_{L^2}^2
 \\
 \les & \sum_{|\tilde\ba|=|\ba|-1} \mathfrak{D}^{0,\tilde\ba,0}
 + \mathfrak{D}^{0, \ba,0},
 \end{align*}
where the first inequality follows from \eqref{hard}, then
\begin{align} \label{Mar-29}
 \widetilde{\mathfrak{D}}(t)\les \mathfrak{D}(t)  +\lt\|\sa^{ {\iota}/{2}} v\rt\|_{L^2}^2.
\end{align}

It follows from \ef{system-a} and the Piola identity $\pl_k(JA^k_i)=0$ $(1\le i\le n)$ that
$$\bar\rho \pl_{t} v_i   +   \pl_k \lt({A}_i^k \bar\rho^\ga J^{1-\ga} \rt)  = - \bar\rho  v_i, \ \ 1\le i\le n-1, $$
which means that for $1\le i \le n-1$,
\begin{align*}
\frac{d}{dt} \int_\Oa \bar\rho   v_i  dy  + \int_\Oa \bar\rho   v_i dy  =- \int_{\pl \Oa}  \bar\rho^\ga {A}_i^k J^{1-\ga} N_k dS =0.
 \end{align*}
Because on the boundary $\{y_n=0\}$, we have $N=(0,\cdots,0, -1)$, and the boundary condition $v_n |_{\{y_n=0\}}=0$ which implies $ x_n  |_{\{y_n=0\}}=0$ and $JA^n_i |_{\{y_n=0\}}=0$.
Thus,
$$
 \int_\Oa \bar\rho(y)  v_i(t,y) dy = e^{-t}  \int_\Oa \bar\rho(y)  v_i(0,y) dy
$$
 and then
$$\int_\Oa \bar\rho(y) \lt(  v_i(t,y)  -  M^{-1}  e^{-t}  \int_\Oa \bar\rho(y) v_i(0,y) dy\rt)  dy
  =0 $$
for $1\le i\le  n-1$.
So, there exists $b(t)\in \Oa$ such that
for  $1\le i \le n-1$,
\begin{align}
& v_i(t,b(t))  = M^{-1}  e^{-t}  \int_\Oa \bar\rho(y)  v_i(0,y) dy,\notag\\
& \lt|  v_i(t,y) \rt|\le \lt|  v_i(t,y) - v_i(t,b(t))  \rt| +| v_i(t,b(t))|\notag\\
&\le  \lt\| \pl  v_i(t,\cdot) \rt\|_{L^\iy}|y-b(t)| +  M^{-1}  e^{-t}    \int_\Oa \lt|\bar\rho(y)  v_i(0,y) \rt|dy\notag\\
& \les  \lt\| \pl  v_i(t,\cdot) \rt\|_{L^\iy}
+ e^{-t} \lt\|    v_i(0,\cdot) \rt\|_{L^\iy},\notag
\end{align}
which, with the aid of $H^{n+3+[\iota]+\iota, \  n+2+[\iota]} \hookrightarrow H^{({n+1+[\iota]-\iota})/{2}}\hookrightarrow
L^\iy$ and \ef{Feb27-2}, implies that
\begin{align}
&\lt|  v_i(t,y) \rt|^2 \les \sum_{j+i\le [\iota]+n+2 } \mathcal{D}^{0,j,i}(t)
+ e^{-2t} \lt\|    v_i(0,\cdot) \rt\|_{L^\iy}^2 \notag \\
& \les \sum_{\substack{m+|\ba|+\al \le [\iota]+n+2 }}
 \lt(\mathfrak{D}^{m,\ba,\al}
 +
 \mathcal{V}^{m,|\ba|,\al}\rt)(t) +  e^{-2t}  \lt\|    v_i(0,\cdot) \rt\|_{L^\iy}^2\notag
\end{align}
for $1\le i\le n-1$. Due to the boundary condition $v_n |_{\{y_n=0\}}=0$,
we can   bound $v_n(t,y)$ similarly and obtain that
\begin{align}
& \lt|  v(t,y) \rt|^2 \les   e^{-2t} \lt\|    v (0,\cdot) \rt\|_{L^\iy}^2\notag
\\
& \ \ +\sum_{\substack{m+|\ba|+\al \le [\iota]+n+2 }}
 \lt(\mathfrak{D}^{m,\ba,\al}
 +
 \mathcal{V}^{m,|\ba|,\al}\rt)(t).\notag
\end{align}
This, together with
 \ef{Mar-29} and \ef{Feb28-1}, gives
\begin{align} \label{3.30}
 \widetilde{\mathfrak{D}}(t)\les \mathfrak{E}(t) \les \mathfrak{D}(t)  +\mathcal{V}(t)
 +  e^{-2t}  \lt\|    v (0,\cdot) \rt\|_{L^\iy}^2.
\end{align}

Let $\da\in (0,1/4)$ be a small constant to be determined later, then we multiply \ef{Feb27-3} by $e^{\da t}$, and use   \ef{3.30} to get
\begin{align*}
&\frac{d}{dt} \lt(e^{\da t}\mathfrak{E}(t)\rt) +\frac{1}{2} e^{\da t}\lt( \mathfrak{D} + 2\mathcal{V} \rt)(t)
 \le \da    e^{\da t}\mathfrak{E}(t)\\
& \les \da e^{\da t}\lt( \mathfrak{D} + \mathcal{V} \rt)(t)
+\da e^{-(2-\da)t}  \lt\|    v (0,\cdot) \rt\|_{L^\iy}^2,
\end{align*}
which means that there exists a constant $\da\le 1/4$ only depending on $\ga$ and $M$ such that
\begin{align*}
\frac{d}{dt} \lt(e^{\da t}\mathfrak{E}(t)\rt) +\frac{1}{4} e^{\da t}\lt( \mathfrak{D} + 2\mathcal{V} \rt)(t)
  \les \da e^{-(2-\da)t}  \lt\|    v (0,\cdot) \rt\|_{L^\iy}^2.
\end{align*}
Integrate the equation above over $[0,t]$ to obtain
\begin{align}\label{3-30}
e^{\da t}\mathfrak{E}(t) +\frac{1}{4}\int_0^t  e^{\da s} \lt(  \mathfrak{D} + \mathcal{V} \rt)(s) ds
  \les    \mathfrak{E}(0)+ \da \lt\|    v (0,\cdot) \rt\|_{L^\iy}^2.
\end{align}
Note that
\begin{align}
& \sum_{m+j+i\le [\iota]+n+3}\mathcal{D}^{m,j,i} \les \mathfrak{D} + \mathcal{V},\label{3-30-1}\\
& \sum_{m+j+i\le [\iota]+n+3}\mathcal{E}^{m+1,j,i}
\les   \sum_{m+j+i\le [\iota]+n+3}\mathcal{D}^{m,j,i}  + \lt\|\sa^{ {\iota}/{2}} v\rt\|_{L^2}^2 \les \mathfrak{E} ,
  \label{3-30-2}
\end{align}
where \ef{3-30-1} follows from \ef{Feb27-2}, and the first inequality of \ef{3-30-2} from a similar way to deriving \ef{Mar-29}, and the second from
  \ef{3-30-1} and \ef{Feb28-1}.
As a conclusion of \ef{Feb28-1}, \ef{3-30-2}, \ef{3.30} and \ef{3-30}, we have
\begin{align}\label{21.6.16}
e^{\da t} \mathscr{E}_1(t)+
\int_0^t e^{\da s} \mathscr{E}_1(s) ds
\les \mathscr{E}_1(0),
\end{align}
where
$
\mathscr{E}_1(t)=\sum_{m+j+i\le [\iota]+n+3} (\mathcal{E}^{m+1,j,i}+\mathcal{V}^{m,j,i}   )(t)
$.
In view of \ef{21June16-1} and \ef{1.27-2}, we see
\begin{align*}
\mathcal{V}^{m,j,i} & \les   \lt\| \sa^{\frac{\iota+i+1}{2}} {\rm curl}_x \pl_t^m\bar\pl^j \pl_n^i v \rt\|_{L^2}^2 +
 \ea_0^2 \sum_{r+l+k\le m+j+i -1}\mathcal{D}^{r,l,k}\notag\\
& \les  \lt\| \sa^{\frac{\iota+i+1}{2}}  \pl_t^m\pl \bar\pl^j \pl_n^i v \rt\|_{L^2}^2 +
 \ea_0^2 \mathcal{E}_I \les \mathcal{E}_I
\end{align*}
for $m+j+i\le [\iota]+n+3$.
This, together with \ef{21.6.16}, proves \ef{n21-4-8}.
\hfill $\Box$

\subsection{A priori estimates for $\mathcal{E}^{0,j,i}$}
It follows from \ef{n21-4-8}  that for $t\in [0, T]$,
\begin{align}\label{5.27-1}
\mathscr{E}_2(t) \les \mathscr{E}_2(0)+ \da^{-1} \mathcal{E}_I(0), \ {\rm where} \ \mathscr{E}_2(t)=\sum_{j+i\le [\iota]+n+3}
\mathcal{E}^{0,j,i}(t),
\end{align}
where $\da\le 1/4$ is a positive constant only depending on $\ga$ and $M$. Clearly,
\begin{subequations}\label{June8}\begin{align}
 & \lt\|\pl v \rt\|_{L^\iy}+\lt\| \pl_t \pl v \rt\|_{L^\iy}  +\lt\|  \sa \pl_n \pl v \rt\|_{L^{\iy}} +\lt\|  \sa \pl_t \pl_n \pl v \rt\|_{L^{\iy}}
 + \lt\|\mathcal{I}^{1,0,0} \rt\|_{L^\iy}\notag\\
 &\qquad + \lt\|\mathcal{I}^{2,0,0} \rt\|_{L^\iy} + \lt\|\sa\mathcal{I}^{1,0,1} \rt\|_{L^\iy} + \lt\|\sa\mathcal{I}^{2,0,1} \rt\|_{L^\iy}
 \les \sqrt{\mathcal{E}_{I}}, \label{June8.a}\\
&  \lt\|\pl \oa \rt\|_{L^\iy} \les \sqrt{\mathscr{E}_2}.\label{June8.b}
\end{align}\end{subequations}
In fact, \ef{June8} follows from \ef{12-8} and \ef{12.10}.

\subsubsection{Curl estimates}

\begin{lem}
It holds that for $j+i= [\iota]+n+4$, and $t\in [0,T]$,
\begin{subequations}\begin{align}
&\lt\|\sa^{\frac{\iota+i+1}{2}} \bar\pl^j\pl_n^i {\rm curl}_x \oa \rt\|_{L^2}^2(t)
\les \lt\|\sa^{\frac{\iota+i+1}{2}} \bar\pl^j\pl_n^i {\rm curl}_x \oa \rt\|_{L^2}^2(0)\notag\\
&\qquad \qquad + \mathcal{V}^{0,j,i}(0)+\ea_0^2 \mathscr{E}_2(0)+ \ea_0^2 \da^{-2}\mathcal{E}_{I}(0),\label{21.6.11-2}\\
&\mathcal{V}^{0,j,i}(t)\les e^{-2t}\mathcal{V}^{0,j,i}(0)+ e^{-\da t} \ea_0^2 \mathcal{E}_I(0), \label{21.6.11-3}
 \end{align}\end{subequations}
where $\da\le 1/4$ is a positive constant  only depending on $\ga$ and $M$.

\end{lem}

{\em Proof}.
It follows from \ef{11-29} that
\begin{align}\label{21.6.11-4}
{\rm curl}_x v =  \lt\{{\rm curl}_x v\big|_{t=0}+ \int_0^t e^{\tau} [\pl_\tau, {\rm curl}_x] v d\tau\rt\}e^{-t},
\end{align}
which implies
\begin{align}
{\rm curl}_x \oa = & {\rm curl}_x \oa \big|_{t=0}
+\lt(1-e^{-t}\rt){\rm curl}_x v\big|_{t=0}
+\int_0^t [\pl_s, {\rm curl}_x] \oa ds
\notag \\
& + \int_0^t e^{-s} \int_0^s e^{\tau} [\pl_\tau, {\rm curl}_x] v d\tau ds.  \notag
\end{align}
For any multi-index $\ba$ and nonnegative integer $\al$, we take $\bar\pl^\ba \pl_n^\al$ onto the equation above to get
\begin{align}
&\bar\pl^\ba \pl_n^\al {\rm curl}_x \oa =  \bar\pl^\ba \pl_n^\al {\rm curl}_x \oa \big|_{t=0}
+\lt(1-e^{-t}\rt)\bar\pl^\ba \pl_n^\al
{\rm curl}_x v\big|_{t=0}\notag
\\
& \ +\int_0^t \bar\pl^\ba \pl_n^\al [\pl_s, {\rm curl}_x] \oa ds
 + \int_0^t e^{-s} \int_0^s e^{\tau} \bar\pl^\ba \pl_n^\al [\pl_\tau, {\rm curl}_x] v d\tau ds.\label{1.25-1}
\end{align}

Notice that for  $  |\ba|+\al=[\iota]+n+4$,
\begin{align}
\bar\pl^\ba \pl_n^\al \lt(\pl_t[{\rm curl}_x v ]_l-[{\rm curl}_x \pl_t v ]_l   \rt)
=\bar\pl^\ba \pl_n^\al \lt(\ea^{ljk}(\pl_t A^r_j)\pl_r v_k \rt)= \sum_{h=1,2,3}
Z_{h,l}^{\ba, \al},\label{21.6.12-1}
\end{align}
where
\begin{align*}
&Z_{1,l}^{\ba, \al}=\ea^{ljk} \sum_{0\le h\le \min\{1,\ \al\}} C(h) \lt\{ \pl_t\lt( (\pl_t   \pl_n^h A^r_j)\bar\pl^{\ba} \pl_n^{\al-h} \pl_r \oa_k\rt)\rt. \\
&\qquad\qquad\qquad \qquad \lt.-(\pl_t^2   \pl_n^h A^r_j)\bar\pl^{\ba } \pl_n^{\al-h} \pl_r \oa_k \rt\},\\
&Z_{2,l}^{\ba, \al}=\ea^{ljk} \sum_{\max\{\al-1, \ 0\}\le h\le \al} C( h)\lt\{ \pl_t\lt( ( \bar\pl^{\ba} \pl_n^h A^r_j) \pl_n^{\al-h} \pl_r v_k\rt)\rt. \\
&\qquad\qquad\qquad \qquad \lt.-( \bar\pl^{\ba} \pl_n^h A^r_j)\pl_t \pl_n^{\al-h} \pl_r v_k \rt\},\\
&\lt|Z_{3,l}^{\ba, \al}\rt|\les \sum_{(r,h)\in \mathcal{S}\setminus \mathcal{S}_1} \mathcal{I}^{1,r,h} \lt|\bar\pl^{|\ba|-r} \pl_n^{\al-h} \pl  v \rt|.
\end{align*}
Here $C( h)=1$ for $h=0$ and $\al$, $C(h)=\al$ for $h=1$ and $\al-1$,  $\mathcal{S}=\{(r,h)\in \mathbb{Z}^2\big|0\le  r \le |\ba|, \  0\le h\le \al\}$,   $\mathcal{S}_1=\{(0,0), \ (0,1),\  (|\ba|,\al -1), \ (|\ba|,\al) \}$. It follows from \ef{5-30-1} and \ef{n21-4-8} that
\begin{align}\label{21.6.12-2}
\lt\|\sa^{\frac{\iota+|\al|+1}{2}} \int_0^s e^{\tau} Z_{3,l}^{\ba, \al} d\tau  \rt\|_{L^2}
\les  \ea_0  \int_0^s e^{\tau} \sqrt{\mathcal{E}_I(\tau)} d\tau
\les \ea_0 \sqrt{\mathcal{E}_I(0)} e^{(1-\da/2)s}.
\end{align}
We integrate by parts over time to achieve
\begin{align*}
&\lt\|\sa^{\frac{\iota+|\al|+1}{2}} \int_0^s e^{\tau} Z_{2,l}^{\ba, \al} d\tau  \rt\|_{L^2}
\\
&\les
\lt\|\sa^{\frac{\iota+|\al|+1}{2}}
\mathcal{I}^{0,|\ba|,\al-1} \pl_n\pl v \rt\|_{L^2}(0)
+\lt\|\sa^{\frac{\iota+|\al|+1}{2}}
\mathcal{I}^{0,|\ba|,\al} \pl v \rt\|_{L^2}(0)
\\
&+ e^s\lt\|\sa^{\frac{\iota+|\al|+1}{2}}
\mathcal{I}^{0,|\ba|,\al-1} \pl_n\pl v \rt\|_{L^2}(s)
+e^s\lt\|\sa^{\frac{\iota+|\al|+1}{2}}
\mathcal{I}^{0,|\ba|,\al} \pl v \rt\|_{L^2}(s)\\
&+\int_0^s e^\tau\lt\|\sa^{\frac{\iota+|\al|+1}{2}}
\mathcal{I}^{0,|\ba|,\al-1} \pl_n\pl v \rt\|_{L^2} d\tau +\int_0^s e^\tau
\lt\|\sa^{\frac{\iota+|\al|+1}{2}}
\mathcal{I}^{0,|\ba|,\al} \pl v \rt\|_{L^2} d\tau\\
&+\int_0^s e^\tau\lt\|\sa^{\frac{\iota+|\al|+1}{2}}
\mathcal{I}^{0,|\ba|,\al-1} \pl_t \pl_n\pl v \rt\|_{L^2} d\tau +\int_0^s e^\tau
\lt\|\sa^{\frac{\iota+|\al|+1}{2}}
\mathcal{I}^{0,|\ba|,\al}\pl_t \pl v \rt\|_{L^2} d\tau.
\end{align*}
In view of \ef{Feb1-e} and \ef{June8.a}, we see that
\begin{align*}
&\lt\|\sa^{\frac{\iota+|\al|+1}{2}}
\mathcal{I}^{0,|\ba|,\al-1} \pl_n\pl v \rt\|_{L^2}\le \lt\|\sa^{\frac{\iota+|\al|-1}{2}}
\mathcal{I}^{0,|\ba|,\al-1}\rt\|_{L^2}  \lt\|\sa\pl_n\pl v \rt\|_{L^\iy}\les \ea_0\sqrt{\mathcal{E}_I} ,\\
&\lt\|\sa^{\frac{\iota+|\al|+1}{2}}
\mathcal{I}^{0,|\ba|,\al} \pl v \rt\|_{L^2}\le \lt\|\sa^{\frac{\iota+|\al|+1}{2}}
\mathcal{I}^{0,|\ba|,\al} \rt\|_{L^2}  \lt\|\pl v \rt\|_{L^\iy}\les \ea_0\sqrt{\mathcal{E}_I},\\
& \lt\|\sa^{\frac{\iota+|\al|+1}{2}}
\mathcal{I}^{0,|\ba|,\al-1} \pl_t \pl_n\pl v \rt\|_{L^2} \le \lt\|\sa^{\frac{\iota+|\al|-1}{2}}
\mathcal{I}^{0,|\ba|,\al-1}\rt\|_{L^2}  \lt\| \sa \pl_t \pl_n\pl v \rt\|_{L^\iy}\les \ea_0\sqrt{\mathcal{E}_I}, \\
&
\lt\|\sa^{\frac{\iota+|\al|+1}{2}}
\mathcal{I}^{0,|\ba|,\al}\pl_t \pl v \rt\|_{L^2} \le \lt\|\sa^{\frac{\iota+|\al|+1}{2}}
\mathcal{I}^{0,|\ba|,\al}\rt\|_{L^2}  \lt\|\pl_t \pl v \rt\|_{L^\iy} \les \ea_0\sqrt{\mathcal{E}_I},
\end{align*}
which, together with \ef{n21-4-8}, mean
\begin{align}\label{21.6.12-3}
\lt\|\sa^{\frac{\iota+|\al|+1}{2}} \int_0^s e^{\tau} Z_{2,l}^{\ba, \al} d\tau  \rt\|_{L^2}
\les \ea_0 \sqrt{\mathcal{E}_I(0)} e^{(1-\da/2)s}.
\end{align}
We can bound $Z_{1,l}^{\ba, \al}$ similarly, and then obtain that
 for $|\ba|+\al = [\iota]+n+4$,
\begin{align}
\lt\|\sa^{\frac{\iota+|\al|+1}{2}} \int_0^t e^{-s} \int_0^s e^{\tau} \bar\pl^\ba \pl_n^\al [\pl_\tau, {\rm curl}_x] v d\tau ds \rt\|_{L^2}^2 \les
\ea_0^2 \da^{-2}\mathcal{E}_I(0) .\label{21.6.11-1}
\end{align}

Note that for  $ |\ba|+\al= [\iota]+n+4$,
\begin{align*}
\bar\pl^\ba \pl_n^\al \lt(\pl_t[{\rm curl}_x \oa ]_l-[{\rm curl}_x \pl_t \oa ]_l   \rt)
=\bar\pl^\ba \pl_n^\al \lt(\ea^{ljk}(\pl_t A^r_j)\pl_r \oa_k \rt)=
Y_{1,l}^{\ba, \al} +Y_{2,l}^{\ba, \al},
\end{align*}
where
\begin{align*}
&Y_{1,l}^{\ba, \al}=\pl_t \lt(\ea^{ljk} (\pl_r \oa_k)\bar\pl^\ba \pl_n^\al  A^r_j \rt)
- \ea^{ljk} ( \pl_r v_k)\bar\pl^\ba \pl_n^\al  A^r_j, \\
& \lt| Y_{2,l}^{\ba, \al} \rt|\les
\sum_ {\substack{ r \le |\ba|, \ k \le \al\\
1\le r+k
  }} \lt|\bar\pl^{r}\pl_n^{k}\pl \oa \rt| \sum_ {\substack{ r_1 \le |\ba|-r\\ k_1 \le \al-k
  }} \mathcal{I}^{0,r_1,k_1} \lt|\bar\pl^{|\ba|-r-r_1}\pl_n^{\al-k-k_1}\pl v \rt|\\
& \qquad \quad \le \sum_ {\substack{ r \le |\ba|, \ k \le \al, \
1\le r+k
  }}\mathcal{I}^{0,r,k} \lt|\bar\pl^{|\ba|-r}\pl_n^{\al-k}\pl v \rt|.
\end{align*}
We integrate by parts over time, and use
 \ef{Feb1-e} and \ef{June8} to get
\begin{align}
&\lt\|\sa^{\frac{\iota+|\al|+1}{2}}\int_0^t Y_{1,l}^{\ba, \al} ds \rt\|_{L^2}\notag\\
\les &  \lt\|\sa^{\frac{\iota+|\al|+1}{2}}
\mathcal{I}^{0,|\ba|,\al} \pl \oa\rt\|_{L^2} (0)+\lt\|\sa^{\frac{\iota+|\al|+1}{2}}
\mathcal{I}^{0,|\ba|,\al} \pl \oa\rt\|_{L^2} (t)\notag\\
&+ \int_0^t \lt\|\sa^{\frac{\iota+|\al|+1}{2}}
\mathcal{I}^{0,|\ba|,\al} \pl v\rt\|_{L^2} ds \notag\\
\les & \ea_0 \sqrt{ \mathscr{E}_2(0)} +\ea_0  \sqrt{\mathscr{E}_2(t)}+ \ea_0 \int_0^t \sqrt{\mathcal{E}_I   (s)} ds, \notag
\end{align}
which, with the aid of \ef{n21-4-8} and \ef{5.27-1}, gives
\begin{align}
\lt\|\sa^{\frac{\iota+|\al|+1}{2}}\int_0^t Y_{1,l}^{\ba, \al} ds \rt\|_{L^2}^2
\les  \ea_0^2 \mathscr{E}_2(0)+ \ea_0^2 \da^{-1}\mathcal{E}_I(0). \notag
\end{align}
It yields from \ef{1.27-2} and \ef{n21-4-8} that
\begin{align}
\lt\|\sa^{\frac{\iota+|\al|+1}{2}}\int_0^t Y_{2,l}^{\ba, \al} ds \rt\|_{L^2}^2
\les  \ea_0^2 \da^{-1}\mathcal{E}_I(0). \notag
\end{align}
So, we have that for $ |\ba|+\al = [\iota]+n+4$,
\begin{align}
\lt\|\sa^{\frac{\iota+|\al|+1}{2}} \int_0^t \bar\pl^\ba \pl_n^\al [\pl_s, {\rm curl}_x] \oa ds \rt\|_{L^2}^2 \les  \ea_0^2 \mathscr{E}_2(0)+ \ea_0^2 \da^{-1}\mathcal{E}_I(0). \label{June8-1}
\end{align}

Now, \ef{21.6.11-2} is a consequence of \ef{1.25-1}, \ef{21.6.11-1} and \ef{June8-1}. \ef{21.6.11-3} follows from \ef{21.6.11-4}, and \ef{21.6.12-1}-\ef{21.6.12-3}, since $Z_{1,l}^{\ba, \al}$ and $Z_{2,l}^{\ba, \al}$ can be bounded in the same way.
\hfill $\Box$

\subsubsection{Elliptic estimates near the bottom}
\begin{lem} It holds that for $j+i= [\iota]+n+4$ and $i\ge 1$, and $t\in [0, T]$,
\begin{align}
&\lt\|\za_1 \pl \bar\pl^j\pl_n^{i} \oa\rt\|_{L^2}^2(t) \les \lt\|\sa^{\frac{\iota+1}{2}} \bar\pl^{[\iota]+n+4} \pl \oa\rt\|_{L^2}^2(t)+ (\mathscr{E}_2+ \da^{-2}\mathcal{E}_I)(0)\notag\\
&
 \quad +\sum_{l+k=[\iota]+n+4, \ 1\le k\le i}
 \lt(\lt\|\sa^{\frac{\iota+k+1}{2}} \bar\pl^l\pl_n^k {\rm curl}_x \oa \rt\|_{L^2}^2+ \mathcal{V}^{0,l,k}\rt)(0),\label{June15-5}
\end{align}
where $\zeta_1=\zeta_1(y_n)$ is a smooth cut-off function satisfying \ef{David-1}, and $\da\le 1/4$ is a positive constant  only depending on $\ga$ and $M$.
\end{lem}

{\em Proof}. We use \ef{1.25} and  the Piola identity $\pl_k(JA^k_i)=0$ $(1\le i\le n)$ to rewrite \ef{system-a}  as
\be\label{perted}
 \pl_{t} v_i +  v_i   +  \sa^{-\iota}  \pl_k \lt( \sa^{\iota+1} \mathcal{F}^k_i \rt)  = 0 , \ \ {\rm where} \ \   \mathcal{F}^k_i= {A}_i^k J^{1-\ga} -\da_i^k.
\ee
Clearly, \ef{perted} is equivalent to
$$
\ga \sa \pl_i {\rm div} \oa= \pl_{t} v_i +  v_i - \sa \mathcal{W}_i  -\nu (\iota+1) \mathcal{F}^n_i,
$$
where
$$\mathcal{W}_i =\ga A^k_iJ^{1-\ga}(A^s_r-\da^s_r)\pl_k\pl_s \oa^r+\ga (A^k_iJ^{1-\ga}-\da^k_i)\pl_k {\rm div} \oa .$$
We take $\bar\pl^j\pl_n^{i-1}$ ($i\ge 1$) onto the equation above to get
\begin{align*}
\ga \sa \bar\pl^j\pl_n^{i}  {\rm div} \oa=
\ga(i-1) \nu \bar\pl^j\pl_n^{i-1}  {\rm div} \oa+ \pl_{t}\bar\pl^j\pl_n^{i-1} v_n \notag\\+ \bar\pl^j\pl_n^{i-1}  v_n
- \bar\pl^j\pl_n^{i-1}(\sa \mathcal{W}_n)  -\nu (\iota+1) \bar\pl^j\pl_n^{i-1} \mathcal{F}^n_n,
\end{align*}
which implies
\begin{align}\label{5-29-1}
&\lt\|\za_1 \bar\pl^j\pl_n^{i}  {\rm div} \oa\rt\|_{L^2} \les  \sqrt{\mathcal{E}_I}+\sqrt{\mathscr{E}_2}+\lt\|\za_1 \bar\pl^j\pl_n^{i-1} \mathcal{F}^n_n
\rt\|_{L^2}
\notag\\
&\quad  +\lt\|\za_1 \bar\pl^j\pl_n^{i-1} \mathcal{W}_n
\rt\|_{L^2}   +(i-1)\lt\|\za_1 \bar\pl^j\pl_n^{i-2} \mathcal{W}_n
\rt\|_{L^2}.
\end{align}
Here and thereafter $j+i=[\iota]+n+4$ and $i\ge 1$.

It follows from \ef{n5.30}, \ef{12.15} and \ef{1120} that
\begin{align}
&\lt\|\za_1 \bar\pl^j\pl_n^{i-1} \mathcal{F}^n_n
\rt\|_{L^2}\les
\lt\|\za_1  \mathcal{I}^{0,j,i-1}
\rt\|_{L^2}\notag
\\
&\quad \les \sum_{\substack{l \le j, \ k \le i-1
\\ l+k \le j+i-2
  }}\lt\| \mathcal{I}^{0,l,k}
\rt\|_{L^\iy(\Oa_1)}  \lt\|\sa^{\frac{\iota+i-k}{2}} \bar\pl^{j-l}\pl_n^{i-1-k}\pl\oa
\rt\|_{L^2} \les \sqrt{\mathscr{E}_2} , \label{June15-1}
\\
&\lt\|\za_1 \bar\pl^j\pl_n^{i-1} \mathcal{W}_n
\rt\|_{L^2}\les \lt\|\za_1 \mathcal{I}^{0,j,i-1}
\rt\|_{L^2}\lt\| \pl^2 \oa
\rt\|_{L^\iy(\Oa_1)}
  \notag
\\
&\quad + \sum_{\substack{l \le j, \ k \le i-1
\\1\le l+k \le j+i-2
  }}\lt\| \mathcal{I}^{0,l,k}
\rt\|_{L^\iy(\Oa_1)}  \lt\|\sa^{\frac{\iota+i-k+1}{2}} \bar\pl^{j-l}\pl_n^{i-1-k}\pl^2\oa
\rt\|_{L^2}
\notag\\
&\quad +\lt\| \pl\oa
\rt\|_{L^\iy(\Oa_1)}  \lt\|\sa^{\frac{\iota+i+1}{2}} \bar\pl^{j}\pl_n^{i-1}\pl^2\oa
\rt\|_{L^2}
 \les \ea_0 \sqrt{\mathscr{E}_2} ,  \label{June15-2}
\end{align}
where
$ \Oa_1=\mathbb{T}^{n-1} \times (0,3\hbar/4)$. Indeed,  the following estimates have been used to derive \ef{June15-2}.
$$\lt\| \pl^2 \oa
\rt\|_{L^\iy(\Oa_1)}\les \lt\| \pl^2 \oa
\rt\|_{H^2(\Oa_1)}\les \sum_{r=2,3,4}\lt\|\sa^{(\iota+r)/2} \pl^r \oa
\rt\|_{L^2} \les \sqrt{\mathscr{E}_2}\le \ea_0,$$
and
$\lt\| \pl \oa
\rt\|_{L^\iy(\Oa_1)} \les \sqrt{\mathscr{E}_2}$. In view of \ef{5-29-1}-\ef{June15-2}, we see that
$$\lt\|\za_1 \bar\pl^j\pl_n^{i}  {\rm div} \oa\rt\|_{L^2}^2 \les \mathcal{E}_I+\mathscr{E}_2,
$$
which, together with \ef{Hodge}, means
\begin{align}
&\lt\|\za_1 \pl \bar\pl^j\pl_n^{i} \oa\rt\|_{L^2}^2 \les \lt\|\za_1  \bar\pl^j\pl_n^{i} \oa\rt\|_{H^1}^2+ {\mathscr{E}_2}
\notag\\
&\les \lt\|\za_1 {\rm curl} \bar\pl^j\pl_n^{i} \oa\rt\|_{L^2 }^2+ \lt\|\za_1 \pl \bar\pl^{j+1}\pl_n^{i-1} \oa\rt\|_{L^2}^2
  +
\mathcal{E}_I+
{\mathscr{E}_2}.\label{June15-3}
\end{align}
In a similar way to deriving \ef{June15-2}, we obtain
\begin{align*}
&\lt\|\za_1  ({\rm curl} \bar\pl^j\pl_n^{i}\oa-{\rm curl}_x \bar\pl^j\pl_n^{i}\oa)\rt\|_{L^2 }\\
&\quad \les \lt\| \pl\oa
\rt\|_{L^\iy(\Oa_1)}  \lt\|\sa^{\frac{\iota+i+1}{2}}\pl \bar\pl^{j}\pl_n^{i}\oa
\rt\|_{L^2}\les \ea_0 \sqrt{\mathscr{E}_2} ,\\
&\lt\|\za_1  ({\rm curl}_x \bar\pl^j\pl_n^{i}\oa-\bar\pl^j\pl_n^{i}{\rm curl}_x \oa)\rt\|_{L^2 }\\
&\quad\les \sum_{\substack{l \le j, \ k \le i
\\1\le l+k\le j+i-2
  }}\lt\| \mathcal{I}^{0,l,k}
\rt\|_{L^\iy(\Oa_1)}  \lt\|\sa^{\frac{\iota+i-k+1}{2}} \bar\pl^{j-l}\pl_n^{i-k}\pl\oa
\rt\|_{L^2}\\
&\quad +\sum_{\substack{l \le j, \ k \le i
\\ j+i-1\le l+k
  }}\lt\| \sa^{\frac{\iota+k+1}{2}}\mathcal{I}^{0,l,k}
\rt\|_{L^2}  \lt\| \bar\pl^{j-l}\pl_n^{i-k}\pl\oa
\rt\|_{L^\iy(\Oa_1)}\les \ea_0 \sqrt{\mathscr{E}_2},
\end{align*}
where  \ef{Feb1-e} has been used. Thus,
\begin{align}
\lt\|\za_1 {\rm curl} \bar\pl^j\pl_n^{i} \oa\rt\|_{L^2 }\les  \lt\|\za_1 \bar\pl^j\pl_n^{i}{\rm curl}_x \oa\rt\|_{L^2 }
+\ea_0 \sqrt{\mathscr{E}_2}. \label{June15-4}
\end{align}
Now, \ef{June15-5} follows from \ef{June15-3},
\ef{June15-4}, \ef{21.6.11-2}, \ef{5.27-1} and
\ef{n21-4-8}.
\hfill $\Box$

\subsubsection{Energy estimates}
\begin{lem} It holds that for $t\in [0, T]$,
\begin{align}\label{haha16}
\mathcal{E}_{II}(t)\les \mathcal{E}(0)+\mathcal{V}_{a}(0).
\end{align}
\end{lem}

{\em Proof}.
For any multi-index $\ba$ and nonnegative integer $\al$ with $|\ba|+\al \ge 1$, we take
$\bar\pl^\ba \pl_n^\al$ onto \ef{perted} to get
\begin{align*}
&\pl_{t} \bar\pl^\ba \pl_n^\al v_i   +   \bar\pl^\ba \pl_n^\al v_i + \sa^{-\iota-\al} \pl_k \lt( \sa^{\iota+\al+1} \bar\pl^\ba \pl_n^\al \mathcal{F}^k_i \rt)\notag\\
& = \al \lt((\pl_k \sa) \pl_n-(\pl_n \sa) \pl_k  \rt)\bar\pl^\ba \pl_n^{\al-1}\mathcal{F}^k_i
={S}_{1;i}^{\ba,\al},
\end{align*}
where $\mathcal{F}^k_i= {A}_i^k J^{1-\ga} -\da_i^k$.
Let $\chi$ be defined by \ef{1-25} and set
\begin{align*}
& {S}_{2;i}^{\ba,\al;k}=
\bar\pl^\ba \pl_n^\al \mathcal{F}^k_i+ J^{1-\ga} \lt( A^k_r A^s_i \pl_s \bar\pl^\ba \pl_n^\al \oa^r + \iota^{-1} A^k_i {\rm div}_x \bar\pl^\ba \pl_n^\al \oa \rt),
\end{align*}
then we  obtain,
using a similar way to deriving \ef{86-1}, that
\begin{align}
\frac{1}{2}\frac{d}{dt} \int \chi^2 \sa^{\iota+\al}  \lt|   \bar\pl^\ba \pl_n^\al v \rt|^2 dy  +  \int \chi^2 \sa^{\iota+\al}  \lt|    \bar\pl^\ba \pl_n^\al v \rt|^2 dy
 \notag\\
+\frac{1}{2}\frac{d}{dt}\sum_{1\le j\le 3}  E_j^{\ba,\al}= \sum_{1\le j\le 3} L_j^{\ba,\al},\label{Jan25}
\end{align}
where
\begin{align*}
&E_1^{\ba,\al}= \int \chi^2 \sa^{\iota+\al+1}  J^{1-\ga}
\lt( \lt|\na_x \bar\pl^\ba \pl_n^\al \oa \rt|^2 + \iota^{-1}  \lt|{\rm div}_x \bar\pl^\ba \pl_n^\al \oa \rt|^2  \rt) dy,\\
& E_2^{\ba,\al} = -\int \chi^2 \sa^{\iota+\al+1}  J^{1-\ga}\lt|{\rm curl}_x \bar\pl^\ba \pl_n^\al \oa \rt|^2  dy, \\
&E_3^{\ba,\al}=  - 2\int \chi^2 \sa^{\iota+\al +1} {S}_{2;i}^{\ba,\al;k}   \pl_k   \bar\pl^\ba \pl_n^\al \oa^i  dy,\\
& L_1^{\ba,\al}=-\int \chi^2 \sa^{\iota+\al+1} \lt(\pl_t {S}_{2;i}^{\ba,\al;k} \rt) \pl_k   \bar\pl^\ba \pl_n^\al \oa^i dy,\\
& L_2^{\ba,\al}= - \int \chi^2 \sa^{\iota+\al+1}  J^{1-\ga}  \lt( \lt[ \nabla_x \bar\pl^\ba \pl_n^\al \oa^r \rt]_i \lt[\nabla_x v^s\rt]_r\lt[\nabla_x  \bar\pl^\ba \pl_n^\al \oa^i\rt]_s \rt.\\
&\  \lt. -\iota^{-1}
\lt(\pl_t A^k_i \rt) \lt({\rm div}_x \bar\pl^\ba \pl_n^\al \oa \rt) \pl_k \bar\pl^\ba \pl_n^\al \oa^i \rt)dy+\frac{1}{2}\int \chi^2  \sa^{\iota+\al+1} \lt(\pl_t  J^{1-\ga} \rt)\\
& \  \times \lt( \lt|\na_x \bar\pl^\ba \pl_n^\al \oa \rt|^2  + \iota^{-1}  \lt|{\rm div}_x \bar\pl^\ba \pl_n^\al \oa \rt|^2
-\lt|{\rm curl}_x \bar\pl^\ba \pl_n^\al \oa \rt|^2   \rt)
 dy,\\
& L_3^{\ba,\al}= \int \chi^2 \sa^{\iota+\al} {S}_{1;i}^{\ba,\al} \bar\pl^\ba \pl_n^\al v^i  dy+2\int \chi( \pl_k \chi )  \sa^{\iota+\al+1} \lt( \bar\pl^\ba \pl_n^\al \mathcal{F}^k_i \rt)    \bar\pl^\ba \pl_n^\al v^i   dy .
\end{align*}

Note that for $|\ba|+\al \ge 2$,
\begin{align*}
&\lt|{S}_{2;i}^{\ba,\al;k}\rt|\les
\sum_ {\substack{ l \le |\ba|, \ r \le \al\\
 1\le l+r\le |\ba|+\al-1
  }} \mathcal{I}^{0,l,r} \lt|\bar\pl^{|\ba|-l}\pl_n^{\al-r}\pl \oa \rt|,\\
& \lt|\pl_t {S}_{2;i}^{\ba,\al;k}\rt|\les
\sum_ {\substack{ l \le |\ba|, \ r \le \al, \
 1\le l+r
  }}\mathcal{I}^{0,l,r} \lt|\bar\pl^{|\ba|-l}\pl_n^{\al-r}\pl v \rt|,
\end{align*}
then it follows from the Cauchy inequality,
\ef{1.27-1} and \ef{1.27-2} that
\begin{align}
&\lt|E_3^{\ba,\al}\rt|\les  \ea_0 (\mathcal{E}_I+
{\mathscr{E}_2}) + \ea_0     E_1^{\ba,\al}, \label{5.29-1}\\
&\lt| L_1^{\ba,\al}\rt|\les
\ea_0 e^{\da t} \mathcal{E}_I + \ea_0  e^{-\da t}   E_1^{\ba,\al}. \label{5.29-2}
\end{align}
In view of \ef{June8.a}, we see that
\begin{align}
&\lt| L_2^{\ba,\al}\rt|\les  \int \chi^2 \sa^{\iota+\al+1} |\pl v| |\pl \bar\pl^\ba \pl_n^\al \oa|^2 dy
\les \sqrt{\mathcal{E}_I}  E_1^{\ba,\al} . \label{5.29-3}
\end{align}
Notice that  for $\al \ge 1$,
\begin{align*}
&\lt|{S}_{1;i}^{\ba,\al}\rt|\les
\lt|\bar\pl^{|\ba|+1}\pl_n^{\al-1}\pl \oa \rt| +
\sum_ {\substack{ l \le |\ba|, \ r \le \al-1\\
 1\le l+r
  }}  \mathcal{I}^{0,l,r}   \lt|\bar\pl^{|\ba|+1-l}\pl_n^{\al-1-r}\pl \oa \rt|,\\
  & \lt|\bar\pl^\ba \pl_n^\al \mathcal{F}^k_i\rt|
  \les \lt|\bar\pl^{|\ba|}\pl_n^{\al}\pl \oa \rt| + \sum_ {\substack{ l \le |\ba|, \ r \le \al-1, \
 1\le l+r
  }} \mathcal{I}^{0,l,r} \lt|\bar\pl^{|\ba|-l}\pl_n^{\al-r}\pl \oa \rt|,
\end{align*}
thus it yields from  the Cauchy inequality and
\ef{1.27-1}  that
\begin{align}
&\lt| L_3^{\ba,\al}\rt|\les
 e^{\da t}  \mathcal{E}_I+  e^{-\da t} \lt( \sum_{|\tilde\ba|=|\ba|+1}E_1^{\tilde \ba, \al-1}+  \lt\|\za_1 \pl \bar\pl^{|\ba|}\pl_n^{\al} \oa\rt\|_{L^2}^2 + \ea_0^2\mathscr{E}_2\rt).\label{5.29-4}
\end{align}
Due to \ef{1.27-1}, \ef{Feb1-e} and \ef{June8.b}, one has
\begin{align*}
&\lt\|\sa^{\frac{\iota+i+1}{2}}  ({\rm curl}_x \bar\pl^j\pl_n^{i}\oa-\bar\pl^j\pl_n^{i}{\rm curl}_x \oa)\rt\|_{L^2 }\\
&\les \sum_{\substack{l \le j, \ k \le i,
\ 1\le l+k\le j+i-1
  }}\lt\|\sa^{\frac{\iota+i+1}{2}} \mathcal{I}^{0,l,k} \bar\pl^{j-l}\pl_n^{i-k}\pl\oa
\rt\|_{L^2}
\\
&\quad +\lt\|\sa^{\frac{\iota+i+1}{2}} \mathcal{I}^{0,j,i}\rt\|_{L^2} \lt\| \pl\oa
\rt\|_{L^\iy}\les \ea_0\sqrt{\mathcal{E}_I+\mathscr{E}_2},
\end{align*}
which gives
\begin{subequations}\label{21J16}\begin{align}
&\lt\|\sa^{\frac{\iota+i+1}{2}}\bar\pl^j\pl_n^{i}{\rm curl}_x \oa\rt\|_{L^2 }\les \lt\|\sa^{\frac{\iota+i+1}{2}}  {\rm curl}_x \bar\pl^j\pl_n^{i}\oa\rt\|_{L^2 }
 +\ea_0\sqrt{\mathcal{E}}
  \les \sqrt{\mathcal{E}}, \label{21J16.a}\\
&\lt\|\sa^{\frac{\iota+i+1}{2}}  {\rm curl}_x \bar\pl^j\pl_n^{i}\oa\rt\|_{L^2 }
\les \lt\|\sa^{\frac{\iota+i+1}{2}}\bar\pl^j\pl_n^{i}{\rm curl}_x \oa\rt\|_{L^2 } +\ea_0\sqrt{\mathcal{E}_I+\mathscr{E}_2}. \label{21J16.b}
\end{align}\end{subequations}
As a conclusion of \ef{5.29-2}-\ef{21J16.a}, \ef{June15-5}, \ef{21.6.11-2}, \ef{5.27-1}, \ef{n21-4-8} and $L_3^{\ba,0}=0$, we arrive at
\begin{subequations}\label{21H616}
\begin{align}
&\sum_{1\le j\le 3}\int_0^t \lt| L_j^{\ba,0}\rt| ds \les
\ea_0\mathcal{E}_I(0) + \da^{-1}\lt(
\ea_0+\sqrt{\mathcal{E}_I(0)}\rt)   \sup_{[0,t]}E_1^{\ba,0},\\
& \sum_{1\le j\le 3}\int_0^t \lt| L_j^{\ba,\al}\rt| ds\les
\da^{-3}\mathcal{E}(0) + \da^{-1} \mathcal{V}_a(0) + \da^{-1}\lt(
\ea_0+\sqrt{\mathcal{E}_I(0)}\rt)   \sup_{[0,t]}E_1^{\ba,\al}
\notag\\
&\quad +\da^{-1} \sup_{[0,t]}\lt(\sum_{|\tilde\ba|=|\ba|+1}E_1^{\tilde\ba,\al-1}
+ \sum_{|\tilde\ba|=[\iota]+n+4}E_1^{\tilde\ba,0}\rt),
\ \ \al\ge 1,
\end{align}\end{subequations}
for
$|\ba|+\al=[\iota]+n+4$.

We integrate \ef{Jan25} over $[0,t]$ with $|\ba|+\al=[\iota]+n+4$ from $\al=0$ to $[\iota]+n+4$ step by step, and use \ef{21.6.11-2}, \ef{21J16}, \ef{5.29-1}, \ef{21H616}, \ef{5.27-1},  \ef{n21-4-8},  and the smallness of $\ea_0$ and $\mathcal{E}_I(0)$ to obtain
\begin{align*}
&\int \chi^2 \sa^{\iota+\al}  \lt|   \bar\pl^\ba \pl_n^\al v \rt|^2 dy  +
\int \chi^2 \sa^{\iota+\al+1}  \lt|  \pl \bar\pl^\ba \pl_n^\al \oa \rt|^2 dy\notag\\
& +  \int_0^t \int \chi^2 \sa^{\iota+\al}  \lt|    \bar\pl^\ba \pl_n^\al v \rt|^2 dyds
\les  \mathcal{E}(0) +   \mathcal{V}_a(0),
\end{align*}
which, with the aid of \ef{June15-5}, implies
\begin{align}\label{wei-1}
\lt\|   \sa^{\frac{\iota  }{2}}  \bar\pl^{[\iota]+n+4}   v \rt\|_{L^2}^2+   \sum_{j+i=[\iota]+n+4}\lt\| \sa^{\frac{\iota+i +1}{2}}
  \pl \bar\pl^j\pl_n^i \oa \rt\|^2_{L^2} \les  \mathcal{E}(0) +   \mathcal{V}_a(0).
\end{align}
Clearly,
\begin{subequations}\label{wei-2}\begin{align}
&\sum_{j+i=[\iota]+n+4, \ i\ge 1}\lt\| \sa^{\frac{\iota+i }{2}}
\bar\pl^j\pl_n^i v \rt\|^2_{L^2} \le \sum_{j+i=[\iota]+n+4, \ i\ge 1} \mathcal{E}^{1,j,i-1} \le \mathcal{E}_I,\\
& \sum_{j+i=[\iota]+n+4}\lt\| \sa^{\frac{\iota+i }{2}}
\bar\pl^j\pl_n^i \oa \rt\|^2_{L^2} \le  \mathscr{E}_2 + \lt\| \sa^{\frac{\iota }{2}}
\bar\pl^{[\iota]+n+4
} \oa \rt\|^2_{L^2}.
\end{align}\end{subequations}
In view of \ef{hard}, we see that
\begin{align*}
 \lt\| \sa^{\frac{\iota }{2}}
\bar\pl^{[\iota]+n+4
} \oa \rt\|^2_{L^2}\les  \lt\| \sa^{\frac{\iota+2 }{2}}
\bar\pl^{[\iota]+n+4
} \oa \rt\|^2_{L^2} +\lt\| \sa^{\frac{\iota+2 }{2}}
\bar\pl^{[\iota]+n+4
}\pl_n \oa \rt\|^2_{L^2}\notag\\
\les \lt\| \sa^{\frac{\iota+1 }{2}}\pl
\bar\pl^{[\iota]+n+3
} \oa \rt\|^2_{L^2} +\lt\| \sa^{\frac{\iota+1}{2}} \pl
\bar\pl^{[\iota]+n+4
} \oa \rt\|^2_{L^2},
\end{align*}
which, together with \ef{wei-1}, \ef{wei-2}, \ef{5.27-1} and \ef{n21-4-8}, proves \ef{haha16}. \hfill $\Box$

\renewcommand{\theequation}{A-\arabic{equation}}
\renewcommand{\thethm}{A-\arabic{thm}}
\setcounter{equation}{0}
\setcounter{thm}{0}
\section*{Appendix}  

\noindent{\bf The weighted Sobolev embedding}.
Let $\mathfrak{U}$ be a bounded smooth domain in $\mathbb{R}^n  $, and $d=d(y)=dist(y, \partial  \mathfrak{U})$ be a distance function to the boundary.
 For any positive real number $a$  and nonnegative integer $b$,  we define the  weighted Sobolev space  $H^{a, b}(  \mathfrak{U})$   by
$$ H^{a, b}(\mathfrak{U}) = \lt\{   d^{a/2}f \in L^2(\mathfrak{U}): \ \  \int_\mathfrak{U}    d^a|\pl^k f|^2dy<\infty, \ \  0\le k\le b\rt\}$$
  with the norm
$ \|f\|^2_{H^{a, b}(\mathfrak{U})} = \sum_{k=0}^b \int_\mathfrak{U}    d^a|\pl^k f|^2dy$. Let $H^s( \mathfrak{U} )$ $(s\ge 0)$ be the standard Sobolev space, then  for $b\ge  {a}/{2}$, we have the following embedding of weighted Sobolev spaces (cf. \cite{18'}):
 $ H^{a, b}(\mathfrak{U} )\hookrightarrow H^{b- {a}/{2}}( \mathfrak{U})$
    with the estimate
  \be\label{wsv} \|f\|_{H^{b- {a}/{2}}( \mathfrak{U})} \le \bar C \|f\|_{H^{a, b}(\mathfrak{U} )} \ee
for some constant $\bar C$ only depending on $a$, $ b$ and $\mathfrak{U}$.\\

\noindent{\bf The Hardy inequality}. Let $k>-1$ be a given real number, $\bar\da$ be a positive constant,  and $f$ be a function satisfying
$
\int_0^{\bar\da} y^{k+2} \lt(f^2 + |f'|^2\rt) dy < \iy,
$
then it holds that
\be\label{hardy'}
\int_0^{\bar\da} y^{k } f^2 dy \le \bar C \int_0^{\bar\da} y^{k+2} \lt(f^2 + |f'|^2\rt) dy
\ee
for a certain constant $\bar C$ only depending on $\bar\da$ and $k$,
whose proof  can be found  in \cite{18'}. Indeed, \ef{hardy'} is a   general version of the  standard Hardy inequality:
$\int_0^\iy |y^{-1} f|^2 dy \le C \int_0^\iy |f'|^2 dy$.
As a consequence of \ef{hardy'}, we have the following  estimates.
\begin{lem}
Let $\Oa=\mathbb{T}^{n-1}\times (0,\hbar)$ and $\sa(y)=\nu(\hbar-y_n)$ with $\nu$ and $\hbar$ being positive constants.
Let $k>-1$ be a given real number and $f$ be a function satisfying $\int_\Omega \sa^{k+2} (f^2 + |\pl_n f|^2 )dy <\iy$, then
\begin{align}\label{hard}
\int_\Omega \sa^k f^2 dy \le
 \bar C \int_\Omega  \sa^{k+2}\lt( f^2+ |\pl_n f|^2 \rt) dy
\end{align}
for some constant $\bar C$ only depending on $k$, $\nu$ and $\hbar$.
\end{lem}

{\em Proof}.  When $n=1$, \ef{hard} follows from \ef{hardy'} and $\sa(y)=\nu(\hbar-y)$. When  $n\ge 2$,
it follows from $\sa(y)=\sa(y_n)$ that
$$
\int_{\Omega} \sa^k(y) f^2(y) dy=\int_0^\hbar \sa^k(y_n)    F^2(y_n)  dy_n,
$$
where
$$
F^2 (y_n)=  \int_{\mathbb{T}^{n-1}}   f^2(y_*,y_n)dy_*, \ \  y_*=(y_1,y_2,\cdots,y_{n-1}).
$$
This, together with \ef{hardy'}, $\sa(y)=\nu(\hbar-y_n)$ and the  H\"{o}lder  inequality, implies that
\begin{align*}
&\int_{\Omega} \sa^k f^2 dy\le \bar C
\int_0^\hbar \sa^{k+2}(y_n) \lt( F^2(y_n)+ \lt|F'(y_n)\rt|^2\rt)dy_n\\
& \le
  \bar C\int_0^\hbar \sa^{k+2}(y_n) \int_{\mathbb{T}^{n-1}}  \lt(   f^2   +  \lt| \pl_n f\rt|^2  \rt)dy_* dy_n
\\
& = \bar C\int_\Omega  \sa^{k+2}\lt( f^2+ |\pl_n f|^2 \rt) dy
\end{align*}
for some constants $\bar C$ only depending on $k$, $\nu$ and $\hbar$.
\hfill $\Box$\\

\noindent{\bf The Hodge decomposition}.
\begin{lem} Let $F$ be a vector field defined on $\Oa=\mathbb{T}^{n-1}\times (0,\hbar)$ such that $F\in L^2(\Oa)$, ${\rm curl} F\in L^2(\Oa)$, ${\rm div} F\in L^2(\Oa)$ and  $\bar\pl F \in L^2(\Oa)$, then
\begin{align}\label{Hodge}
\|F\|_{H^1(\Oa)}\le \bar{C} \lt(\|F\|_{L^2(\Oa)}
+ \|{\rm curl} F\|_{L^2(\Oa)} + \|{\rm div} F\|_{L^2(\Oa)} + \|\bar\pl F\|_{L^2(\Oa)} \rt)
\end{align}
for some constant $\bar{C}$ only depending on $\Oa$.
\end{lem}
{\em Proof}. The proof of  \ef{Hodge} can be found in Section 6 of \cite{10'}. Indeed, \ef{Hodge} is a conclusion of the following Hodge-type elliptic estimate:
\begin{align*}
\|F\|_{H^1(\Oa)}\le \bar C \lt(
\|F\|_{L^2(\Oa)}
+ \|{\rm curl} F\|_{L^2(\Oa)} + \|{\rm div} F\|_{L^2(\Oa)}
+\|\bar\pl F \cdot N\|_{H^{-0.5}(\bar\Gamma)} \rt),
\end{align*}
 and
 the following  normal trace estimate:
\begin{align*}
\|\bar\pl F \cdot N\|_{H^{-0.5}(\bar\Gamma)}
\le \bar{C} \lt( \|\bar\pl F\|_{L^2(\Oa)}+\|{\rm div} F\|_{L^2(\Oa)}  \rt),
\end{align*}
where $\bar\Gamma=\{y_n=0\}\cup \{y_n=\hbar\}$, $N$ denotes the outward unit-normal to $\bar\Gamma$, and $\bar{C}$ are positive constants only depending on $\Oa$.
\hfill $\Box$\\

\noindent{\bf The basic identity}.
Recall the following identity in Lemma 4.4 of \cite{HZeng}.
\begin{lem}
For any vector field  $F $ with $F^i=F_i$, we have
\begin{align}
& A^k_rA^s_i (\pl_s  F^r)  \pl_t \pl_k F^i  =2^{-1} \pl_t \lt( |\nabla_x F|^2 - |{\rm curl}_x F|^2 \rt) \notag\\
&\qquad+    \lt[ \nabla_x F^r \rt]_i \lt[\nabla_x v^s\rt]_r\lt[\nabla_x F^i\rt]_s, \label{nabt}\\
&A^k_r A^s_i ( \pl_s F^r ) \pl_k F^i  = |\nabla_x F|^2 - |{\rm curl}_x F|^2. \label{nab}
\end{align}
\end{lem}

\noindent{\bf Acknowledgements}.
This research was supported in part by NSFC  Grants  12171267 and 11822107.

\bibliographystyle{plain}

\end{document}